\documentclass[11pt]{article}
\usepackage{amssymb,latexsym}
\usepackage{epsfig}
\usepackage{eufrak}
\usepackage{amsmath}
\usepackage{mathrsfs}
\usepackage{color}

\newtheorem{theorem}{Theorem}
\newtheorem{lemma}{Lemma}
\newcommand{\be}{\begin{equation}}
\newcommand{\ee}{\end{equation}}
\newcommand{\bee}{\begin{eqnarray*}}
\newcommand{\eee}{\end{eqnarray*}}
\newcommand{\bel}{\begin{eqnarray}}
\newcommand{\eel}{\end{eqnarray}}
\newcommand{\bec}{\begin{cases}}
\newcommand{\eec}{\end{cases}}
\newcommand{\bem}{\begin{bmatrix}}
\newcommand{\eem}{\end{bmatrix}}

\newcommand{\la}{\label}
\newcommand{\li}{\left}
\newcommand{\ri}{\right}
\newcommand{\DEF}{\stackrel{\mathrm{def}}{=}}
\newcommand{\ovl}{\overline}
\newcommand{\udl}{\underline}
\newcommand{\lc}{\lceil}
\newcommand{\rc}{\rceil}

\newcommand{\ep}{\epsilon}
\newcommand{\vep}{\varepsilon}

\newcommand{\si}{\sigma}

\newcommand{\de}{\delta}

\newcommand{\ga}{\gamma}

\newcommand{\vse}{\vartheta}
\newcommand{\se}{\theta}
\newcommand{\Se}{\Theta}

\newcommand{\ze}{\zeta}
\newcommand{\al}{\alpha}

\newcommand{\ro}{\rho}

\newcommand{\om}{\omega}
\newcommand{\Om}{\Omega}

\newcommand{\f}{\frac}
\newcommand{\sq}{\sqrt}
\newcommand{\cd}{\cdots}

\newcommand{\qu}{\quad}
\newcommand{\qqu}{\qquad}
\newcommand{\fa}{\forall}

\newcommand{\mscr}{\mathscr}
\newcommand{\mcal}{\mathcal}
\newcommand{\mbf}{\mathbf}
\newcommand{\bb}{\mathbb}

\newcommand{\wh}{\widehat}
\newcommand{\wt}{\widetilde}
\newcommand{\mrm}{\mathrm}
\newcommand{\bs}{\boldsymbol}

\newcommand{\ap}{\approx}

\newcommand{\sh}{\slash}
\newcommand{\tx}{\text}
\newcommand{\iy}{\infty}

\newcommand{\pa}{\partial}

\newcommand{\bed}{\begin{description}}
\newcommand{\eed}{\end{description}}
\newcommand{\bei}{\begin{itemize}}
\newcommand{\eei}{\end{itemize}}
\newcommand{\ben}{\begin{enumerate}}
\newcommand{\een}{\end{enumerate}}
\newcommand{\bib}{\bibitem}
\newcommand{\beL}{\begin{lemma}}
\newcommand{\eeL}{\end{lemma}}
\newcommand{\beT}{\begin{theorem}}
\newcommand{\eeT}{\end{theorem}}
\newcommand{\sect}{\section}

\newcommand{\bpf}{\begin{pf}}
\newcommand{\epf}{\end{pf}}
\newcommand{\bsk}{\bigskip}
\newcommand{\bi}{\binom}

\newcommand{\pfbox}{\hfill\mbox{$\Box$}}

\newenvironment{pf}{\paragraph*{Proof{\rm.}}}{\pfbox\bigskip}

\begin{document}

\title{{\bf Exact Methods for Multistage Estimation of a Binomial  Proportion} }

\author{${\mrm{Zhengjia} \; \mrm{Chen} \;}$ and Xinjia Chen \thanks{Dr. Zhengjia Chen is working with Department of Biostatistics and Bioinformatics, Emory University, Atlanta, GA 30322; Email: zchen38@emory.edu.
Dr. Xinjia Chen is working with Department of Electrical Engineering, Southern University at Baton Rouge, LA 70813; Email:
xinjiachen@engr.subr.edu.} }

\date{February 2013}

\maketitle

\begin{abstract}

We first review existing sequential methods for estimating a binomial proportion. Afterward, we propose a new family of group sequential
sampling schemes for estimating a binomial proportion with prescribed margin of error and confidence level. In particular, we establish the
uniform controllability of coverage probability and the asymptotic optimality for such a family of sampling schemes. Our theoretical results
establish the possibility that the parameters of this family of sampling schemes can be determined so that the prescribed level of confidence is
guaranteed with little waste of samples.  Analytic bounds for the cumulative distribution functions and expectations of sample numbers are
derived.   Moreover, we discuss the inherent connection of various sampling schemes. Numerical issues are addressed for improving the accuracy
and efficiency of computation. Computational experiments are conducted for comparing sampling schemes. Illustrative examples are given for
applications in clinical trials.

\end{abstract}

\section{Introduction} \la{SecInct}

Estimating a binomial proportion is a problem of ubiquitous significance in many areas of engineering and sciences.  For economical reasons and
other concerns, it is important to use as fewer as possible samples to guarantee the required reliability of estimation.  To achieve this goal,
sequential sampling schemes can be very useful.  In a sequential sampling scheme, the total number of observations is not fixed in advance.  The
sampling process is continued stage by stage until a pre-specified stopping rule is satisfied.  The stopping rule is evaluated with accumulated
observations.  In many applications, for administrative feasibility, the sampling experiment is performed in a group fashion.  Similar to group
sequential tests \cite[Section 8]{CSW}, \cite{Jennsion}, an estimation method based on taking samples by groups and evaluating them sequentially
is referred to as a group sequential estimation method. It should be noted that group sequential estimation methods  are  general enough to
include fixed-sample-size and fully sequential procedures as special cases. Particularly, a fixed-sample-size method can be viewed as a group
sequential procedure of only one stage.  If the increment between the sample sizes of consecutive stages is equal to $1$, then the group
sequential method is actually a fully sequential method.

It is a common contention that statistical inference, as a unique science to quantify the uncertainties of inferential statements, should avoid
errors in the quantification of uncertainties, while minimizing the sampling cost.  That is, a statistical inferential method is expected to be
exact and efficient.  The conventional notion of exactness is that no approximation is involved, except the roundoff error due to finite word
length of computers. Existing sequential methods for estimating a binomial proportion are dominantly of asymptotic nature (see, e.g.,
\cite{Chow1965, BKGosh, MGosh, Lai, Siegmund} and the references therein). Undoubtedly, asymptotic techniques provide approximate solutions and
important insights for the relevant problems. However, any asymptotic method inevitably introduces unknown error in the resultant approximate
solution due to the necessary use of a finite number of samples.  In the direction of non-asymptotic sequential estimation, the primary goal is
to ensure that the true coverage probability is above the pre-specified confidence level for any value of the associated parameter, while the
required sample size is as low as possible. In this direction,  Mendo and Hernando \cite{Mendo2}  developed an inverse binomial sampling scheme
for estimating a binomial proportion with relative precision.  Tanaka \cite{Tanaka} developed a rigorous method for constructing fixed-width
sequential confidence intervals for a binomial proportion. Although no approximation is involved, Tanaka's method is very conservative due to
the bounding techniques employed in the derivation of sequential confidence intervals. Franz\'{e}n \cite{frazen} studied the construction of
fixed-width sequential confidence intervals for a binomial proportion. However, no effective method for defining stopping rules is proposed in
\cite{frazen}. In his later paper \cite{frazenB}, Franz\'{e}n proposed to construct fixed-width confidence intervals based on sequential
probability ratio tests (SPRTs) invented by Wald \cite{Wald}. His method can generate fixed-sample-size confidence intervals based on SPRTs.
Unfortunately, he made a fundamental flaw by mistaking that if the width of the fixed-sample-size confidence interval decreases to be smaller
than the pre-specified length as the number of samples is increasing, then the fixed-sample-size confidence interval at the termination of
sampling process is the desired fixed-width sequential confidence interval guaranteeing the prescribed confidence level.   More recently, Jesse
Frey published a paper \cite{Frey} in {\it The American Statistician} (TAS) on the classical problem of sequentially estimating a binomial
proportion with prescribed margin of error and confidence level.   Before Frey submitted his original manuscript to TAS in July 2009, a general
framework of multistage parameter estimation had been established by Chen \cite{Chen1241v1, Chen1241v12, Chen4679v1, Chen0430v1, Chen3458v2},
which provides exact methods for estimating parameters of common distributions with various error criterion.   This framework is also proposed
in \cite{Chen_SPIE}. The approach of Frey \cite{Frey} is similar to that of Chen \cite{Chen1241v1, Chen1241v12, Chen4679v1, Chen0430v1,
Chen3458v2} for the specific problem of estimating a binomial proportion with prescribed margin of error and confidence level.

In this paper, our primary interests are in the exact sequential methods for the estimation of a binomial proportion with prescribed margin of
error and confidence level.  We first introduce the exact approach established in \cite{Chen1241v1, Chen1241v12, Chen4679v1, Chen0430v1,
Chen3458v2}.  In particular, we introduce the inclusion principle proposed in \cite{Chen3458v2} and its applications to the construction of
concrete stopping rules.  We investigate the connection among various stopping rules.  Afterward, we propose a new family of stopping rules
which are extremely simple and accommodate some existing stopping rules as special cases. We provide rigorous justification for the feasibility
and asymptotic optimality of such stopping rules. We prove that the prescribed confidence level can be guaranteed uniformly for all values of a
binomial proportion by choosing appropriate parametric values for the stopping rule.  We show that as the margin of error tends to zero, the
sample size tends to the attainable minimum as if the binomial proportion were exactly known. We derive analytic bounds for distributions and
expectations of sample numbers. In addition, we address some critical computational issues and propose methods to improve the accuracy and
efficiency of numerical calculation. We conduct extensive numerical experiment to study the performance of various stopping rules. We determine
parametric values for the proposed stopping rules to achieve unprecedentedly efficiency while guaranteeing prescribed confidence levels.  We
attempt to make our proposed method as user-friendly as possible so that it can be immediately applicable even for layer persons.

The remainder of the paper is organized as follows.  In Section 2, we introduce the exact approach proposed in \cite{Chen1241v1, Chen1241v12,
Chen4679v1, Chen0430v1, Chen3458v2}. In Section 3, we discuss the general principle of constructing stopping rules.  In Section 4, we propose a
new family of sampling schemes and investigate their feasibility, optimality and analytic bounds of the distribution and expectation of sample
numbers.  In Section 5, we compare various computational methods.  In particular, we illustrate why the natural method of evaluating coverage
probability based on gridding parameter space is neither rigorous nor efficient. In Section 6, we present numerical results for various sampling
schemes. In Section 7, we illustrate the applications of our group sequential method in clinical trials. Section 8 is the conclusion. The proofs
of theorems are given in appendices. Throughout this paper, we shall use the following notations. The empty set is denoted by $\emptyset$.   The
set of positive integers is denoted by $\bb{N}$. The ceiling function is denoted by $\lc . \rc$. The notation $\Pr \{ E \mid \se \}$ denotes the
probability of the event $E$ associated with parameter $\se$.  The expectation of a random variable is denoted by $\bb{E}[.]$.  The standard
normal distribution is denoted by $\Phi(.)$. For $\al \in (0, 1)$, the notation $\mcal{Z}_\al$ denotes the critical value such that $\Phi (
\mcal{Z}_\al ) = 1 - \al$. For $n \in \bb{N}$, in the case that $X_1, \cd, X_n$ are i.i.d. samples of $X$, we denote the sample mean $\f{
\sum_{i=1}^n X_i }{n}$ by $\ovl{X}_n$, which is also called the relative frequency when $X$ is a Bernoulli random variable. The other notations
will be made clear as we proceed.

\sect{How Can It Be Exact?}  \la{SecHow}

In many areas of scientific investigation, the outcome of an experiment is of dichotomy nature and can be modeled as a Bernoulli random variable
$X$, defined in probability space $(\Om, \Pr, \mscr{F})$,  such that
\[
\Pr \{ X = 1 \} = 1 - \Pr \{ X = 0 \} = p \in (0, 1), \] where $p$ is referred to as a binomial proportion. In general, there is no analytical
method for evaluating the binomial proportion $p$.  A frequently-used approach is to estimate $p$ based on i.i.d. samples $X_1, X_2, \cd$ of
$X$.  To reduce the sampling cost, it is appropriate to estimate $p$ by a multistage sampling procedure. More formally, let $\vep \in (0, 1)$
and $1 - \de$, with $\de \in (0, 1)$,  be the pre-specified margin of error and confidence level respectively. The objective is to construct a
sequential estimator $\wh{\bs{p}}$ for $p$ based on a multistage sampling scheme such that \be \la{ineqcov} \Pr \{ | \wh{\bs{p}} - p | < \vep
\mid p \} \geq 1 - \de \ee for any $p \in (0, 1)$. Throughout this paper, the probability $\Pr \{ | \wh{\bs{p}} - p | < \vep \mid p \}$ is
referred to as the {\it coverage probability}. Accordingly, the probability $\Pr \{ | \wh{\bs{p}} - p | \geq \vep \mid p \}$ is referred to as
the {\it complementary coverage probability}.   Clearly, a complete construction of a multistage estimation scheme needs to determine the number
of stages, the sample sizes for all stages, the stopping rule, and the estimator for $p$.  Throughout this paper, we let $s$ denote the number
of stages and let $n_\ell$ denote the number of samples at the $\ell$-th stages.  That is, the sampling process consists of $s$ stages with
sample sizes $n_1 < n_2 < \cd < n_s$. For $\ell = 1, 2, \cd, s$, define $K_\ell = \sum_{i=1}^{n_\ell} X_i$ and $\wh{\bs{p}}_\ell = \f{ K_\ell
}{n_\ell}$.  The stopping rule is to be defined in terms of $\wh{\bs{p}}_\ell, \; \ell = 1, \cd, s$.  Of course, the index of stage at the
termination of the sampling process, denoted by $\bs{l}$, is a random number. Accordingly, the number of samples at the termination of the
experiment, denoted by $\mbf{n}$,  is a random number which equals $n_{\bs{l}}$.  Since for each $\ell$, $\wh{\bs{p}}_\ell$ is a
maximum-likelihood and minimum-variance unbiased estimator of $p$, the sequential estimator for $p$ is taken as \be \la{estdef} \wh{\bs{p}} =
\wh{\bs{p}}_{\bs{l}} = \f{ \sum_{i = 1}^{n_{\bs{l}}} X_i }{ n_{\bs{l}} } = \f{ \sum_{i = 1}^{\mbf{n}} X_i }{ \mbf{n} }. \ee In the above
discussion, we have outlined the general characteristics of a multistage sampling scheme for estimating a binomial proportion. It remains to
determine the number of stages, the sample sizes for all stages, and the stopping rule so that the resultant estimator $\wh{\bs{p}}$ satisfies
(\ref{ineqcov}) for any $p \in (0, 1)$.

Actually, the problem of sequential estimation of a binomial proportion has been treated by Chen \cite{Chen1241v1, Chen1241v12, Chen4679v1,
Chen0430v1, Chen3458v2} in a general framework of multistage parameter estimation.  The techniques of \cite{Chen1241v1, Chen1241v12, Chen4679v1,
Chen0430v1, Chen3458v2} are sufficient to offer exact solutions for a wide range of sequential estimation problems, including the estimation of
a binomial proportion as a special case. The central idea of the approach in \cite{Chen1241v1, Chen1241v12, Chen4679v1, Chen0430v1, Chen3458v2}
is the control of coverage probability by a single parameter $\ze$, referred to as the {\it coverage tuning parameter},  and the adaptive
rigorous checking of coverage guarantee by virtue of bounds of coverage probabilities. It is recognized in \cite{Chen1241v1, Chen1241v12,
Chen4679v1, Chen0430v1, Chen3458v2} that, due to the discontinuity of the coverage probability on parameter space, the conventional method of
evaluating the coverage probability for a finite number of parameter values is neither rigorous not computationally efficient for checking the
coverage probability guarantee.

As mentioned in the introduction, Frey published an article \cite{Frey} in TAS on the sequential estimation of a binomial proportion with
prescribed margin of error and confidence level.  For clarity of presentation, the comparison of the works of Chen and Frey is given in Section
5.4. In the remainder of this section, we shall only introduce the idea and techniques of \cite{Chen1241v1, Chen1241v12, Chen4679v1, Chen0430v1,
Chen3458v2}, which had been precedentially developed  by Chen before Frey submitted his original manuscript to TAS in July 2009. We will
introduce the approach of \cite{Chen1241v1, Chen1241v12, Chen4679v1, Chen0430v1, Chen3458v2} with a focus on the special problem of estimating a
binomial proportion with prescribed margin of error and confidence level.

\subsection{Four Components Suffice} \la{SubsecECS}

The exact methods of \cite{Chen1241v1, Chen1241v12, Chen4679v1, Chen0430v1, Chen3458v2} for multistage parameter estimation have four main
components as follows:

\noindent (I) Stopping rules parameterized by the coverage tuning parameter $\ze > 0$ such that the associated coverage probabilities can be
made arbitrarily close to $1$ by choosing $\ze > 0$ to be a sufficiently small number.

\noindent (II) Recursively computable lower and upper bounds for the complementary coverage probability for a given $\ze$ and an interval of
parameter values.

\noindent (III) Adapted Branch and Bound Algorithm.

\noindent (IV) Bisection coverage tuning.

Without looking at the technical details, one can see that these four components are sufficient for constructing a sequential estimator so that
the prescribed confidence level is guaranteed.  The reason is as follows:  As lower and upper bounds for the complementary coverage probability
are available, the global optimization technique, Branch and Bound (B\&B) Algorithm \cite{Land}, can be used to compute exactly the maximum of
complementary coverage probability on the whole parameter space.  Thus, it is possible to check rigorously whether the coverage probability
 associated with a given $\ze$ is no less than the pre-specified confidence level.
 Since the coverage probability can be controlled by $\ze$, it is possible
to determine $\ze$ as large as possible to guarantee the desired confidence level by a bisection search. This process is referred to as
bisection coverage tuning in \cite{Chen1241v1, Chen1241v12, Chen4679v1, Chen0430v1, Chen3458v2}. Since a critical subroutine needed for
bisection coverage tuning is to check whether the coverage probability is no less than the pre-specified confidence level, it is not necessary
to compute exactly the maximum of the complementary coverage probability. Therefore, Chen revised the standard B\&B algorithm to reduce the
computational complexity and called the improved algorithm  as the Adapted B\&B Algorithm. The idea is to adaptively partition the parameter
space as many subintervals. If for all subintervals, the upper bounds of the complementary coverage probability are no greater than $\de$,  then
declare that the coverage probability is guaranteed. If there exists a subinterval for which the lower bound of the complementary coverage
probability is greater than $\de$,  then declare that the coverage probability is not guaranteed. Continue partitioning the parameter space if
no decision can be made.  The four components are illustrated in the sequel under the headings of stopping rules, interval bounding, adapted
branch and bound, and bisection coverage tuning.

\subsection{Stopping Rules} \la{SubSecSR}

The first component for the exact sequential estimation of a binomial proportion is the stopping rule for constructing a sequential estimator
such that the coverage probability can be controlled by the coverage tuning parameter $\ze$.  For convenience of describing some concrete
stopping rules, define {\small
\[ \mscr{M} (z,\se) = \bec z \ln \f{\se}{z} + (1 - z) \ln \f{1 - \se}{1 - z} &
\tx{for} \; z \in (0,1) \; \tx{and} \; \se \in (0, 1),\\
\ln(1-\se) & \tx{for} \; z = 0 \; \tx{and} \; \se \in (0, 1),\\
\ln \se &  \tx{for} \; z = 1 \; \tx{and} \; \se \in (0, 1),\\
- \iy &  \tx{for} \; z \in [0, 1] \; \tx{and} \; \se \notin (0, 1) \eec
\]} and
\[
S (k, l, n, p) = \bec \sum_{i = k}^l \bi{n}{i} p^i (1 -
p)^{n - i} & \tx{for} \; p \in (0, 1),\\
0 & \tx{for} \; p \notin (0, 1)  \eec
\]
where $k$ and $l$ are integers such that $0 \leq k \leq l \leq n$. Assume that $0 < \ze \de < 1$.  For the purpose of controlling the coverage
probability $\Pr \{ | \wh{\bs{p}} - p | < \vep \mid p \}$ by the coverage tuning parameter, Chen has proposed four stopping rules as follows:

\noindent {\bf Stopping Rule A}: Continue sampling until  $\mscr{M} ( \f{1}{2} - |\f{1}{2} - \wh{\bs{p}}_\ell | , \f{1}{2} - |\f{1}{2} -
\wh{\bs{p}}_\ell | + \vep) \leq \f{ \ln ( \ze \de  ) } { n_\ell }$ for some $\ell \in \{1, \cd, s\}$.

\noindent {\bf Stopping Rule B}: Continue sampling until $( | \wh{\bs{p}}_\ell - \f{1}{2} | - \f{2 }{3} \vep )^2 \geq \f{1}{4} + \f{ \vep^2
n_\ell } {2 \ln (\ze \de) }$ for some $\ell \in \{1, \cd, s\}$.

\noindent {\bf Stopping Rule C}: Continue sampling until $S (K_\ell, n_\ell, n_\ell, \wh{\bs{p}}_\ell - \vep ) \leq \ze \de$ and $S (0, K_\ell,
n_\ell, \wh{\bs{p}}_\ell + \vep ) \leq \ze \de$ for some $\ell \in \{1, \cd, s\}$.

\noindent {\bf Stopping Rule D}: Continue sampling until {\small $n_\ell \geq \wh{\bs{p}}_\ell (1 - \wh{\bs{p}}_\ell) \f{2}{\vep^2} \ln
\f{1}{\ze \de}$} for some $\ell \in \{1, \cd, s\}$.

Stopping Rule A was first proposed in \cite[Theorem 7]{Chen1241v1} and restated in \cite[Theorem 16]{Chen1241v12}. Stopping Rule B was first
proposed in \cite[Theorem 1]{Chen4679v1} and represented as the third stopping rule in \cite[Section 4.1.1]{Chen1241v16}. Stopping Rule C
originated from \cite[Theorem 1]{Chen0430v1} and was restated as the first stopping rule in \cite[Section 4.1.1]{Chen1241v16}. Stopping Rule D
was described in the remarks following Theorem 7 of \cite{Chen1241v4}.  All these stopping rules can be derived from the general principles
proposed in \cite[Section 3]{Chen3458v2} and \cite[Section 2.4]{Chen_SPIE}.

Given that a stopping rule can be expressed in terms of $\wh{\bs{p}}_\ell$ and $n_\ell$ for $\ell = 1, \cd, s$, it is possible to find a
bivariate function $\mscr{D} (.,.)$ on $\{ (z, n): z \in [0,1], \; n \in \bb{N} \}$, taking values from $\{0, 1\}$, such that the stopping rule
can be stated as: Continue sampling until $\mscr{D}(\wh{\bs{p}}_\ell, n_\ell) = 1$ for some $\ell \in \{1, \cd, s \}$.   It can be checked that
such representation applies to Stopping Rules A, B, C, and D.  For example, Stopping Rule B can be expressed in this way by virtue of function
$\mscr{D} (.,.)$ such that
\[
\mscr{D} (z, n) = \bec 1 & \tx{if {\small $(|z - \f{1}{2}| -
\f{2}{3} \vep)^2 \geq \f{1}{4} + \f{\vep^2 n}{2 \ln (\ze \de)}$}},\\
0  & \tx{otherwise} \eec  \]
The motivation of introducing function $\mscr{D} (.,.)$ is to parameterize the stopping rule in terms of design
parameters.   The function $\mscr{D} (.,.)$ determines the form of the stopping rule and consequently, the sample sizes for all stages can be
chosen as functions of design parameters.   Specifically, let {\small \bel &  & N_{\mrm{min}} = \min \li \{ n \in \bb{N}:
\mscr{D} \li ( \f{k}{n}, n \ri ) = 1 \; \tx{for some nonnegative integer $k$ not exceeding $n$} \ri \}, \qqu \qqu \la{def8a}\\
&  & N_{\mrm{max}} = \min \li \{ n \in \bb{N}: \mscr{D} \li ( \f{k}{n}, n \ri ) = 1 \; \tx{for all nonnegative integer $k$ not exceeding $n$}
\ri \}.  \la{def8b} \eel} To avoid unnecessary checking of the stopping criterion and thus reduce administrative cost, there should be a
possibility that the sampling process is terminated at the first stage.  Hence,  the minimum sample size $n_1$ should be chosen to ensure that
$\{ \mbf{n} = n_1 \} \neq \emptyset$. This implies that the sample size $n_1$ for the first stage can be taken as $N_{\mrm{min}}$.  On the other
hand, since the sampling process must be terminated at or before the $s$-th stage, the maximum sample size $n_s$ should be chosen to guarantee
that $\{ \mbf{n}
> n_s \} = \emptyset$.  This implies that the sample size $n_s$ for the last stage can be taken as $N_{\mrm{max}}$.  If the number of stages $s$ is given, then the
sample sizes for stages in between $1$ and $s$ can be chosen as $s-2$ integers between $N_{\mrm{min}}$ and $N_{\mrm{max}}$.  Specially, if the
group sizes are expected to be approximately equal, then the sample sizes can be taken as \be \la{equalss} n_\ell = \li \lc N_{\mrm{min}} +
\f{\ell - 1}{s - 1} ( N_{\mrm{max}} -  N_{\mrm{min}}) \ri \rc, \qqu \ell = 1, \cd, s. \ee
 Since the stopping rule is associated with the coverage tuning parameter $\ze$, it follows that the number of
stages $s$ and the sample sizes $n_1, n_2, \cd, n_s$ can be expressed as functions of $\ze$. In this sense, it can be said that the stopping
rule is parameterized by the coverage tuning parameter $\ze$.  The above method of parameterizing stopping rules has been used in
\cite{Chen1241v1, Chen1241v12, Chen4679v1, Chen0430v1} and proposed in \cite[Section 2.1, page 9]{Chen1241v16}.

\subsection{Interval Bounding} \la{SubsecIB}

The second component for the exact sequential estimation of a binomial proportion is the method of bounding the complementary coverage
probability $\Pr \{ | \wh{\bs{p}} - p | \geq \vep \mid p \}$ for $p$ in an interval $[a, b]$ contained by interval $(0, 1)$. Applying Theorem 8
of \cite{Chen1241v12} to the special case of a Bernoulli distribution immediately yields {\small \be \la{bounbi} \Pr \{ \wh{\bs{p}} \leq a -
\vep \mid b \} + \Pr \{ \wh{\bs{p}} \geq b + \vep \mid a \} \leq \Pr \{ | \wh{\bs{p}} - p | \geq \vep \mid p \} \leq \Pr \{ \wh{\bs{p}} \leq b -
\vep \mid a \} + \Pr \{ \wh{\bs{p}} \geq a + \vep \mid b \} \qqu \ee} for all $p \in [a, b] \subseteq (0, 1)$.  The bounds of (\ref{bounbi}) can
be shown as follows:    Note that $\Pr \{ \wh{\bs{p}} \leq a - \vep \mid p \} + \Pr \{ \wh{\bs{p}} \geq b + \vep \mid p \} \leq \Pr \{ |
\wh{\bs{p}} - p | \geq \vep \mid p \} = \Pr \{ \wh{\bs{p}} \leq p - \vep \mid p \} + \Pr \{ \wh{\bs{p}} \geq p + \vep \mid p \} \leq \Pr \{
\wh{\bs{p}} \leq b - \vep \mid p \} + \Pr \{ \wh{\bs{p}} \geq a + \vep \mid p \}$ for $p \in [a, b] \subseteq (0, 1)$.  As a consequence of the
monotonicity of $\Pr \{ \wh{\bs{p}} \geq \vartheta \mid p \}$ and $\Pr \{ \wh{\bs{p}} \leq \vartheta \mid p \}$ with respect to $p$, where
$\vse$ is a real number independent of $p$, the lower and upper bounds of {\small $\Pr \{ | \wh{\bs{p}} - p | \geq \vep \mid p \}$} for $p \in
[a, b] \subseteq (0, 1)$ can be given as $\Pr \{ \wh{\bs{p}} \leq a - \vep \mid b \} + \Pr \{ \wh{\bs{p}} \geq b + \vep \mid a \}$ and $\Pr \{
\wh{\bs{p}} \leq b - \vep \mid a \} + \Pr \{ \wh{\bs{p}} \geq a + \vep \mid b \} $ respectively.

In page 15, equation (1) of \cite{Chen1241v12}, Chen proposed to apply the recursive method of Schultz \cite[Section 2]{sch} to compute the
lower and upper bounds of $\Pr \{ | \wh{\bs{p}} - p | \geq \vep \mid p \}$ given by (\ref{bounbi}).  It should be pointed out that such lower
and upper bounds of $\Pr \{ | \wh{\bs{p}} - p | \geq \vep \mid p \}$ can also be computed  by the recursive path-counting method of Franz\'{e}n
\cite[page 49]{frazen}.

\subsection{Adapted Branch and Bound} \la{SubsecBB}

The third component for the exact sequential estimation of a binomial proportion is the Adapted B\&B Algorithm, which was proposed in
\cite[Section 2.8]{Chen1241v12},  for quick determination of whether the coverage probability is no less than $1 - \de$ for any value of the
associated parameter. Such a task of checking the coverage probability is also referred to as checking the coverage probability guarantee. Given
that lower and upper bounds of the complementary coverage probability on an interval of parameter values can be obtained by the interval
bounding techniques, this task can be accomplished by applying the B\&B Algorithm \cite{Land} to compute exactly the maximum of the
complementary coverage probability on the parameter space.   However, in our applications, it suffices to determine whether the maximum of the
complementary coverage probability $\Pr \{ | \wh{\bs{p}} - p | \geq \vep \mid p \}$ with respect to $p \in (0, 1)$ is greater than the
confidence parameter $\de$. For fast checking whether the maximal complementary coverage probability exceeds $\de$, Chen proposed to reduce the
computational complexity by revising the standard B\&B Algorithm as the Adapted B\&B Algorithm in \cite[Section 2.8]{Chen1241v12}.  To describe
this algorithm, let $\mcal{I}_{\mrm{init}}$ denote the parameter space $(0, 1)$.  For an interval $\mcal{I} \subseteq \mcal{I}_{\mrm{init}}$,
let $\max \Psi (\mcal{I})$ denote the maximum of the complementary coverage probability $\Pr \{ | \wh{\bs{p}} - p | \geq \vep \mid p \}$ with
respect to $p \in \mcal{I}$.  Let $\Psi_{\mrm{lb}} (\mcal{I})$ and $\Psi_{\mrm{ub}} (\mcal{I})$ be respectively the lower and upper bounds of
$\Psi (\mcal{I})$, which can be obtained by the interval bounding techniques introduced in Section \ref{SubsecIB}.  Let $\eta > 0 $ be a
pre-specified tolerance, which is much smaller than $\de$.  The Adapted B\&B Algorithm of \cite{Chen1241v12} is represented with a slight
modification as follows.

 \bsk

\begin{tabular} {|l |}
\hline $ \nabla \; \tx{Let} \; k \leftarrow 0, \;  l_0 \leftarrow \Psi_{\mrm{lb}} ( \mcal{I}_{\mrm{init}} ) \;
\tx{and} \;  u_0 \leftarrow \Psi_{\mrm{ub}} ( \mcal{I}_{\mrm{init}} )$.\\
$ \nabla \; \tx{Let} \; \mscr{S}_0  \leftarrow  \{ \mcal{I}_{\mrm{init}} \} \; \tx{if $u_0 > \de$.  Otherwise, let $\mscr{S}_0$ be empty}$.\\
$ \nabla \; \tx{While $\mscr{S}_k$ is nonempty, $l_k < \de$ and $u_k$ is greater than $\max \{ l_k + \eta, \;\de \}$, do the following:}$\\
$ \indent \indent \diamond \; \tx{Split each interval in $\mscr{S}_k$} \; \tx{as two new intervals of equal length}.$\\
$ \indent \indent \;\; \; \; \tx{Let $S_k$ denote the set of all new intervals obtained from this splitting procedure}$.\\
$ \indent \indent \diamond \; \tx{Eliminate any interval $\mcal{I}$ from $S_k$
such that $\Psi_{\mrm{ub}} (\mcal{I}) \leq \de$}$.\\
$ \indent \indent \diamond \; \tx{Let} \; \mscr{S}_{k+1} \; \tx{be the set $S_k$ processed by the above elimination procedure}$.\\
$ \indent \indent  \diamond \; \tx{Let} \; l_{k+1} \leftarrow \max_{\mcal{I} \in \mscr{S}_{k+1}} \Psi_{\mrm{lb}} (\mcal{I}) \; \tx{and} \;
 u_{k+1} \leftarrow \max_{\mcal{I} \in \mscr{S}_{k+1}} \Psi_{\mrm{ub}} (\mcal{I})$.  $\tx{Let} \; k \leftarrow k + 1$.\\
$ \nabla \;  \tx{If $\mscr{S}_k$ is empty and $l_k < \de$, then declare $\max \Psi ( \mcal{I}_{\mrm{init}} ) \leq \de$}$.\\
$ \; \; \; \; \tx{Otherwise, declare $\max \Psi ( \mcal{I}_{\mrm{init}} )
> \de$}$.\\ \hline
\end{tabular}

\bsk

It should be noted that for a sampling scheme of symmetrical stopping boundary, the initial interval $\mcal{I}_{\mrm{init}}$ may be taken as
$(0, \f{1}{2})$ for the sake of efficiency.  In Section 5.1, we will illustrate why the Adapted B\&B Algorithm is superior than the direct
evaluation based on gridding parameter space.   As will be seen in Section 5.2, the objective of the Adapted B\&B Algorithm can also be
accomplished by the Adaptive Maximum Checking Algorithm due to Chen \cite[Section 3.3 ]{Chen1241v16} and rediscovered by Frey in the second
revision of his manuscript submitted to TAS in April 2010 \cite[Appendix]{Frey}.   An explanation is given in Section 5.3 for the advantage of
working with the complementary coverage probability.

\subsection{Bisection Coverage Tuning}  The fourth component for the exact sequential estimation of a binomial proportion is Bisection Coverage Tuning.
Based on the adaptive rigorous checking of coverage probability, Chen proposed in \cite[Section 2.7]{Chen1241v1} and \cite[Section
2.6]{Chen1241v12} to apply a bisection search method to determine maximal $\ze$ such that the coverage probability is no less than $1 - \de$ for
any value of the associated parameter. Moreover, Chen has developed asymptotic results in \cite[page 21, Theorem 18]{Chen1241v12} for
determining the initial interval of $\ze$ needed for the bisection search.  Specifically, if the complementary coverage probability $\Pr \{ |
\wh{\bs{p}} - p | \geq \vep \mid p \}$ associated with $\ze = \ze_0$ tends to $\de$ as $\vep \to 0$, then the initial interval of $\ze$ can be
taken as $[\ze_0 2^{i}, \; \ze_0 2^{i+1}]$, where $i$ is the largest integer such that the complementary coverage probability associated with
$\ze = \ze_0 2^{i}$ is no greater than $\de$ for all $p \in (0, 1)$.  By virtue of a bisection search, it is possible to obtain $\ze^* \in
[\ze_0 2^{i}, \; \ze_0 2^{i+1}]$ such that the complementary coverage probability associated with $\ze = \ze^*$ is guaranteed to be no greater
than $\de$ for all $p \in (0, 1)$.

\section{Principle of Constructing Stopping Rules} \la{SecPrinciple}

In this section, we shall illustrate the inherent connection between various stopping rules. It will be demonstrated that a lot of stopping
rules can be derived by virtue of the inclusion principle proposed by Chen \cite[Section 3]{Chen3458v2}.

\subsection{Inclusion Principle}  The problem of estimating a binomial proportion can be considered as a special case of
parameter estimation for a random variable $X$ parameterized by {\small $\se \in \Se$}, where the objective is to construct a sequential
estimator $\wh{\bs{\se}}$ for $\se$ such that $\Pr \{  |\wh{\bs{\se}} - \se | < \vep \mid \se  \} \geq 1 - \de$ for any $\se \in \Se$.  Assume
that the sampling process consists of $s$ stages with sample sizes {\small $n_1 < n_2 < \cd < n_s$}.  For $\ell = 1, \cd, s$, define an
estimator {\small $\wh{\bs{\se}}_\ell$} for $\se$ in terms of samples {\small $X_1, \cd, X_{n_\ell}$} of $X$.  Let $[L_\ell, U_\ell], \; \ell =
1, 2, \cd, s$ be a sequence of confidence intervals such that for any $\ell$, $[L_\ell, U_\ell]$ is defined in terms of {\small $X_1, \cd,
X_{n_\ell}$} and that the coverage probability $\Pr \{ L_\ell \leq \se \leq U_\ell \mid \se \}$ can be made arbitrarily close to $1$ by choosing
$\ze > 0$ to be a sufficiently small number.  In Theorem 2 of \cite{Chen3458v2}, Chen proposed the following general stopping rule: \be \la{GST}
 \tx{Continue sampling until {\small $U_\ell - \vep \leq
\wh{\bs{\se}}_\ell \leq L_\ell + \vep$} for some $\ell \in \{1, \cd, s \}$}. \ee At the termination of the sampling process, a sequential
estimator for $\se$  is taken as  {\small $\wh{\bs{\se}} = \wh{\bs{\se}}_{\bs{l}}$}, where $\bs{l}$ is the index of stage at the termination of
sampling process.

Clearly, the general stopping rule (\ref{GST}) can be restated as follows:

Continue sampling until the confidence interval $[L_\ell, U_\ell]$ is included by interval $[\wh{\bs{\se}}_\ell - \vep, \; \wh{\bs{\se}}_\ell +
\vep]$ for some $\ell \in \{1, \cd, s \}$.

The sequence of confidence intervals are parameterized by $\ze$ for purpose of controlling the coverage probability $\Pr \{ | \wh{\bs{\se}} -
\se | < \vep \mid \se \}$.   Due to the inclusion relationship $[L_\ell, U_\ell] \subseteq  [\wh{\bs{\se}}_\ell - \vep, \; \wh{\bs{\se}}_\ell +
\vep]$, such a general methodology of using a sequence of confidence intervals to construct a stopping rule for controlling the coverage
probability is referred to as the {\it inclusion principle}. It is asserted by Theorem 2 of \cite{Chen3458v2} that \be \la{analytical}
 \Pr \{ | \wh{\bs{\se}} - \se | < \vep \mid \se \} \geq 1 - s \ze \de \qu \fa \se \in \Se
\ee provided that $\Pr \{  L_\ell < \se < U_\ell \mid \se  \} \geq 1 - \ze \de$ for $\ell = 1, \cd, s$ and $\se \in \Se$.  This demonstrates
that if the number of stages $s$ is bounded with respective to $\ze$, then the coverage probability $\Pr \{ | \wh{\bs{\se}} - \se | < \vep \mid
\se \}$ associated with the stopping rule derived from the inclusion principle can be controlled by $\ze$.  Actually, before explicitly
proposing the inclusion principle in \cite{Chen3458v2}, Chen  had extensively applied the inclusion principle in \cite{Chen1241v1, Chen1241v12,
Chen4679v1, Chen0430v1} to construct stopping rules for estimating parameters of various distributions such as binomial, Poisson, geometric,
hypergeometric, normal distributions, etc.    A more general version of the inclusion principle is proposed in \cite[Section 2.4]{Chen_SPIE}.
For simplicity of the stopping rule, Chen had made effort to eliminate the computation of confidence limits.

In the context of estimating a binomial proportion $p$, the inclusion principle immediately leads to the following general stopping rule: \be
\la{ruleBI}
 \tx{Continue sampling until $\wh{\bs{p}}_\ell - \vep \leq L_\ell \leq U_\ell \leq \wh{\bs{p}}_\ell + \vep$ for some $\ell \in \{1,
\cd, s \}$.} \ee  Consequently, the sequential estimator for $p$ is taken as $\wh{\bs{p}}$ according to (\ref{estdef}).  It should be pointed
out that the stopping rule (\ref{ruleBI}) had been rediscovered by Frey in Section 2, the 1st paragraph of \cite{Frey}. The four stopping rules
considered in his paper follow immediately from applying various confidence intervals to the general stopping rule (\ref{ruleBI}).

In the sequel, we will illustrate how to apply (\ref{ruleBI}) to the derivation of Stopping Rules A, B, C, D introduced in Section 2.2 and other
specific stopping rules.

\subsection{Stopping Rule from Wald Intervals}  \la{SubSRWI}

By virtue of Wald's method of interval estimation for a binomial proportion $p$, a sequence of confidence intervals $[L_\ell, U_\ell], \; \ell =
1, \cd, s$ for $p$ can be constructed such that
\[
L_\ell = \wh{\bs{p}}_\ell - \mcal{Z}_{\ze \de}  \sq{  \f{ \wh{\bs{p}}_\ell (1 - \wh{\bs{p}}_\ell) }{n_\ell} }, \qqu U_\ell = \wh{\bs{p}}_\ell +
\mcal{Z}_{\ze \de}  \sq{  \f{ \wh{\bs{p}}_\ell (1 - \wh{\bs{p}}_\ell) }{n_\ell} }, \qqu \ell = 1, \cd, s
\]
and that $\Pr \{ L_\ell \leq p \leq U_\ell \mid p \} \ap 1 - 2 \ze \de$ for $\ell = 1, \cd, s$ and $p \in (0, 1)$.  Note that, for $\ell = 1,
\cd, s$,  the event $\{ \wh{\bs{p}}_\ell - \vep \leq L_\ell \leq U_\ell \leq \wh{\bs{p}}_\ell + \vep \}$ is the same as the event {\small $\li
\{ \li ( \wh{\bs{p}}_\ell - \f{1}{2} \ri )^2 \geq \f{1}{4} - n_\ell \li ( \f{ \vep } {\mcal{Z}_{\ze \de} } \ri )^2 \ri \}$}. So, applying this
sequence of confidence intervals to (\ref{ruleBI}) results in the stopping rule ``continue sampling until $\li ( \wh{\bs{p}}_\ell - \f{1}{2} \ri
)^2 \geq \f{1}{4} - n_\ell \li ( \f{ \vep } {\mcal{Z}_{\ze \de} } \ri )^2$ for some $\ell \in \{ 1, \cd, s \}$''.  Since for any $\ze \in (0,
\f{1}{\de})$, there exists a unique number $\ze^\prime \in(0, \f{1}{\de})$ such that $\mcal{Z}_{\ze \de} = \sqrt{2 \ln \f{1}{\ze^\prime \de}}$,
this stopping rule is equivalent to ``Continue sampling until $\li ( \wh{\bs{p}}_\ell - \f{1}{2} \ri )^2 \geq \f{1}{4} + \f{ \vep^2 n_\ell } { 2
\ln (\ze \de) }$ for some $\ell \in \{1, \cd, s \}$.'' This stopping rule is actually the same as Stopping Rule D, since {\small $\li \{ \li (
\wh{\bs{p}}_\ell - \f{1}{2} \ri )^2 \geq \f{1}{4} + \f{ \vep^2 n_\ell } { 2 \ln (\ze \de) } \ri \} = \li \{ n_\ell \geq \wh{\bs{p}}_\ell (1 -
\wh{\bs{p}}_\ell) \f{2}{\vep^2} \ln \f{1}{\ze \de} \ri \}$} for  $\ell \in \{ 1, \cd, s \}$.

\subsection{Stopping Rule from Revised Wald Intervals} \la{SubsecSRRW}

Define $\wt{\bs{p}}_\ell = \f{  n_\ell \; \wh{\bs{p}}_\ell + a }{n_\ell + 2 a }$ for $\ell = 1, \cd, s$, where $a$ is a positive number.
Inspired by Wald's method of interval estimation for $p$, a sequence of confidence intervals $[L_\ell, U_\ell], \; \ell = 1, \cd, s$ can be
constructed such that
\[
L_\ell = \wh{\bs{p}}_\ell - \mcal{Z}_{\ze \de}  \sq{  \f{ \wt{\bs{p}}_\ell (1 - \wt{\bs{p}}_\ell) }{n_\ell} }, \qqu U_\ell = \wh{\bs{p}}_\ell +
\mcal{Z}_{\ze \de}  \sq{  \f{ \wt{\bs{p}}_\ell (1 - \wt{\bs{p}}_\ell) }{n_\ell} }
\]
and that $\Pr \{ L_\ell \leq p \leq U_\ell \mid p \} \ap 1 - 2 \ze \de$ for $\ell = 1, \cd, s$ and $p \in (0, 1)$.  This sequence of confidence
intervals was applied by Frey \cite{Frey} to the general stopping rule (\ref{ruleBI}).  As a matter of fact,  such idea of revising Wald
interval {\small $\li [ \; \ovl{X}_n - \mcal{Z}_{\ze \de} \sq{ \f{ \ovl{X}_n (1 - \ovl{X}_n) }{n}   }, \; \; \ovl{X}_n + \mcal{Z}_{\ze \de} \sq{
\f{ \ovl{X}_n (1 - \ovl{X}_n) }{n}   } \; \ri ]$}  by replacing the relative frequency $\ovl{X}_n = \f{ \sum_{i=1}^n X_i }{n}$ involved in the
confidence limits with $\wt{p}_a = \f{n \ovl{X}_n + a}{n + 2a}$ had been proposed by H. Chen \cite[Section 4]{CHen_CI}.

As can be seen from Section 2, page 243, of Frey \cite{Frey}, applying (\ref{ruleBI}) with the sequence of revised Wald intervals yields the
stopping rule ``Continue sampling until $\li ( \wt{\bs{p}}_\ell - \f{1}{2} \ri )^2 \geq \f{1}{4} + \f{ \vep^2 n_\ell } { 2 \ln (\ze \de) }$ for
some $\ell \in \{1, \cd, s \}$.'' Clearly, replacing $\wh{\bs{p}}_\ell$ in Stopping Rule D with {\small $\wt{\bs{p}}_\ell = \f{ a + n_\ell
\wh{\bs{p}}_\ell} { n_\ell + 2 a}$} also leads to this stopping rule.

\subsection{Stopping Rule from Wilson's Confidence Intervals} \la{SubsecSRWI}

Making use of the interval estimation method of Wilson  \cite{Wilson}, one can obtain a sequence of confidence intervals $[L_\ell, U_\ell], \;
\ell = 1, \cd, s$ for $p$ such that {\small \[ L_\ell = \max \li \{ 0, \; \f{ \wh{\bs{p}}_\ell + \f{ \mcal{Z}_{\ze \de}^2 }{ 2 n_\ell}  -
\mcal{Z}_{\ze \de} \sq{ \f{ \wh{\bs{p}}_\ell ( 1 - \wh{\bs{p}}_\ell ) }{n_\ell} + \li ( \f{ \mcal{Z}_{\ze \de} }{ 2 n_\ell } \ri )^2 } }{ 1 +
\f{ \mcal{Z}_{\ze \de}^2 }{n_\ell} } \ri \}, \qu  U_\ell =  \min \li \{ 1, \;  \f{ \wh{\bs{p}}_\ell + \f{ \mcal{Z}_{\ze \de}^2 }{ 2 n_\ell}  +
\mcal{Z}_{\ze \de} \sq{ \f{ \wh{\bs{p}}_\ell ( 1 - \wh{\bs{p}}_\ell ) }{n_\ell} + \li ( \f{ \mcal{Z}_{\ze \de} }{ 2 n_\ell } \ri )^2 } }{ 1 +
\f{ \mcal{Z}_{\ze \de}^2 }{n_\ell} } \ri \} \;\; \]} and that $\Pr \{ L_\ell \leq p \leq U_\ell \mid p \} \ap 1 - 2 \ze \de$ for $\ell = 1, \cd,
s$ and $p \in (0, 1)$. It should be pointed out that the sequence of Wilson's confidence intervals  has been applied by Frey \cite[Section 2,
page 243]{Frey} to the general stopping rule (\ref{ruleBI}) for estimating a binomial proportion.

Since a stopping rule directly involves the sequence of  Wilson's confidence intervals is cumbersome, it is desirable to eliminate the
computation of Wilson's confidence intervals in the stopping rule. For this purpose, we need to use the following result.

\beT \la{Wisleq} Assume that $0 < \ze \de < 1$ and $0 < \vep < \f{1}{2}$. Then, Wilson's confidence intervals satisfy $\{ \wh{\bs{p}}_\ell -
\vep \leq L_\ell \leq U_\ell \leq \wh{\bs{p}}_\ell + \vep \} = \li \{ \li ( \li | \wh{\bs{p}}_\ell - \f{1}{2} \ri | - \vep \ri )^2 \geq \f{1}{4}
- n_\ell \li ( \f{ \vep } {\mcal{Z}_{\ze \de} } \ri )^2 \ri \}$ for $\ell = 1, \cd, s$. \eeT

See Appendix  \ref{Wisleq_app} for a proof.  As a consequence of Theorem \ref{Wisleq} and the fact that for any $\ze \in (0, \f{1}{\de})$, there
exists a unique number $\ze^\prime \in(0, \f{1}{\de})$ such that $\mcal{Z}_{\ze \de} = \sqrt{2 \ln \f{1}{\ze^\prime \de}}$, applying the
sequence of  Wilson's confidence intervals to (\ref{ruleBI}) leads to the following stopping rule: Continue sampling until \be
\la{normalmassart} \li ( \li | \wh{\bs{p}}_\ell - \f{1}{2} \ri | - \vep \ri )^2 \geq \f{1}{4} + \f{ \vep^2 n_\ell } { 2 \ln (\ze \de) } \ee for
some $\ell \in \{1, \cd, s \}$.

\subsection{Stopping Rule from Clopper-Pearson Confidence Intervals} \la{SubsecSRCPI}

Applying the interval estimation method of Clopper-Pearson \cite{Clopper}, a sequence of confidence intervals $[L_\ell, U_\ell], \; \ell = 1,
\cd, s$ for $p$ can be obtained such that $\Pr \{ L_\ell \leq p \leq U_\ell \mid p \} \geq 1 - 2 \ze \de$ for $\ell = 1, \cd, s$ and $p \in (0,
1)$, where the upper confidence limit $U_\ell$ satisfies the equation $S (0, K_\ell, n_\ell, U_\ell) = \ze \de$ if $K_\ell < n_\ell$; and the
lower confidence limit $L_\ell$ satisfies the equation $S (K_\ell, n_\ell, n_\ell,  L_\ell) = \ze \de$ if $K_\ell >  0$.  The well known
equation (10.8) in \cite[page 173]{Feller} implies that $S(0, k, n, p)$, with $0 \leq k < n$,  is decreasing with respect to $p \in (0, 1)$ and
that $S(k, n, n, p)$, with $0 < k \leq n$, is increasing with respect to $p \in (0, 1)$. It follows that {\small  \bee  \li \{ \wh{\bs{p}}_\ell
- \vep \leq L_\ell  \ri \} = \li \{ 0 < \wh{\bs{p}}_\ell - \vep \leq L_\ell \ri \} \cup \{ \wh{\bs{p}}_\ell \leq \vep \} = \li \{
\wh{\bs{p}}_\ell
> \vep, \; S (K_\ell, n_\ell, n_\ell, \wh{\bs{p}}_\ell - \vep) \leq \ze \de
 \ri \} \cup \{ \wh{\bs{p}}_\ell \leq  \vep \} \qqu \qqu \; \; &  & \\
= \li \{ \wh{\bs{p}}_\ell > \vep, \; S (K_\ell, n_\ell, n_\ell, \wh{\bs{p}}_\ell - \vep) \leq \ze \de
 \ri \} \cup \{ \wh{\bs{p}}_\ell \leq  \vep, \; S (K_\ell, n_\ell, n_\ell, \wh{\bs{p}}_\ell - \vep) \leq \ze \de \}
  = \li \{ S (K_\ell, n_\ell, n_\ell, \wh{\bs{p}}_\ell - \vep) \leq \ze \de
 \ri \}  &  &
 \eee}
and {\small \bee  \li \{ \wh{\bs{p}}_\ell + \vep \geq U_\ell  \ri \} = \li \{ 1 > \wh{\bs{p}}_\ell + \vep \geq U_\ell \ri \} \cup \{
\wh{\bs{p}}_\ell \geq 1 - \vep \} = \li \{ \wh{\bs{p}}_\ell < 1 - \vep, \; S (0, K_\ell, n_\ell, \wh{\bs{p}}_\ell + \vep) \leq \ze \de
 \ri \} \cup \{ \wh{\bs{p}}_\ell \geq  1 - \vep \} &  & \\
 = \li \{ \wh{\bs{p}}_\ell < 1 - \vep, \; S (0, K_\ell, n_\ell, \wh{\bs{p}}_\ell + \vep) \leq \ze \de
 \ri \} \cup \{ \wh{\bs{p}}_\ell \geq  1 - \vep, \; S (0, K_\ell, n_\ell, \wh{\bs{p}}_\ell + \vep) \leq \ze \de \} \qqu \qqu \; \; &  & \\
  = \li \{ S (0, K_\ell, n_\ell, \wh{\bs{p}}_\ell + \vep) \leq \ze \de
 \ri \}  \qqu \qqu \qqu \qqu  \qqu \qqu \qqu \qqu \qqu \qqu \qqu \qqu \qqu \qu \; \; &  &
 \eee}
for $\ell = 1, \cd, s$.  Consequently,
\[
\{ \wh{\bs{p}}_\ell - \vep \leq L_\ell \leq U_\ell  \leq \wh{\bs{p}}_\ell + \vep  \} = \{ S (K_\ell , n_\ell, n_\ell, \wh{\bs{p}}_\ell - \vep)
\leq \ze \de, \; S (0, K_\ell, n_\ell, \wh{\bs{p}}_\ell + \vep) \leq \ze \de \} \]
 for $\ell = 1, \cd, s$.  This demonstrates that applying the
sequence of Clopper-Pearson confidence intervals to the general stopping rule (\ref{ruleBI}) gives Stopping Rule C.

It should be pointed out that Stopping Rule C was rediscovered by J. Frey  as the third stopping rule in Section 2, page 243 of his paper
\cite{Frey}.

\subsection{Stopping Rule from Fishman's Confidence Intervals}  \la{SubsecSRFI}

By the interval estimation method of Fishman \cite{Fishman}, a sequence of confidence intervals $[L_\ell, U_\ell], \; \ell = 1, \cd, s$  for $p$
can be obtained such that
{\small \bee L_\ell = \bec 0 & \tx{if} \; \wh{\bs{p}}_\ell = 0,\\
\{ \se_\ell \in (0, \wh{\bs{p}}_\ell): \mscr{M} ( \wh{\bs{p}}_\ell, \se_\ell) = \f{ \ln(\ze \de)  }{n_\ell} \}  & \tx{if} \; \wh{\bs{p}}_\ell >
0 \eec  \qu \;  U_\ell = \bec 1 & \tx{if} \; \wh{\bs{p}}_\ell = 1,\\
\{ \se_\ell \in (\wh{\bs{p}}_\ell, 1): \mscr{M} ( \wh{\bs{p}}_\ell, \se_\ell) = \f{ \ln(\ze \de)  }{n_\ell} \}  & \tx{if} \; \wh{\bs{p}}_\ell <
1 \eec \eee}   Under the assumption that $0 < \ze \de < 1$ and $0 < \vep < \f{1}{2}$, by similar techniques as the proof of Theorem 7 of
\cite{Chen1241v4}, it can be shown that  $\{ \wh{\bs{p}}_\ell - \vep \leq L_\ell \leq U_\ell \leq \wh{\bs{p}}_\ell + \vep \} = \{  \mscr{M} (
\f{1}{2} - |\f{1}{2} - \wh{\bs{p}}_\ell | , \f{1}{2} - |\f{1}{2} - \wh{\bs{p}}_\ell | + \vep) \leq \f{ \ln ( \ze \de ) } { n_\ell } \}$ for
$\ell = 1, \cd, s$.   Therefore, applying the sequence of confidence intervals of Fishman to the general stopping rule (\ref{ruleBI}) gives
Stopping Rule A.

It should be noted that Fishman's confidence intervals are actually derived from the Chernoff bounds of the tailed probabilities of the sample
mean of Bernoulli random variable. Hence, Stopping Rule A is also referred to as the stopping rule from Chernoff bounds in this paper.

\subsection{Stopping Rule from Confidence Intervals of Chen et. al.}  \la{SubsecSRIC}

Using the interval estimation method of Chen et. al. \cite{CZA_CI}, a sequence of confidence intervals $[L_\ell, U_\ell], \; \ell = 1, \cd, s$
for $p$ can be obtained such that \bee  &  & L_\ell = \max \li \{ 0, \; \wh{\bs{p}}_\ell + \frac{3}{4} \; \frac{ 1 - 2 \wh{\bs{p}}_\ell - \sqrt{
1 + \f{9 n_\ell} { 2 \ln \f{1}{\ze \de} } \; \wh{\bs{p}}_\ell ( 1 - \wh{\bs{p}}_\ell ) } } {1 + \f{9 n_\ell} { 8 \ln \f{1}{\ze \de} } } \ri \}, \\
&  &  U_\ell = \min \li \{1, \; \wh{\bs{p}}_\ell + \frac{3}{4} \; \frac{ 1 - 2 \wh{\bs{p}}_\ell + \sqrt{ 1 + \f{9 n_\ell} { 2 \ln \f{1}{\ze \de}
} \; \wh{\bs{p}}_\ell  ( 1 - \wh{\bs{p}}_\ell ) } } {1 + \f{9 n_\ell} { 8 \ln \f{1}{\ze \de} } } \ri \} \eee and that $\Pr \{ L_\ell \leq p \leq
U_\ell \mid p \} \geq 1 - 2 \ze \de$ for $\ell = 1, \cd, s$ and $p \in (0, 1)$.  Under the assumption that $0 < \ze \de < 1$ and $0 < \vep <
\f{1}{2}$, by similar techniques as the proof of Theorem 1 of \cite{Chen4679v2}, it can be shown that $\{ \wh{\bs{p}}_\ell - \vep \leq L_\ell
\leq U_\ell \leq \wh{\bs{p}}_\ell + \vep \} = \{ ( | \wh{\bs{p}}_\ell - \f{1}{2} | - \f{2 }{3} \vep )^2 \geq \f{1}{4} + \f{ \vep^2 n_\ell } {2
\ln (\ze \de) } \}$ for $\ell = 1, \cd, s$. This implies that applying the sequence of confidence intervals of Chen et. al. to the general
stopping rule (\ref{ruleBI}) leads to Stopping Rule B.

Actually, the confidence intervals of Chen et. al. \cite{CZA_CI} are derived from Massart's inequality \cite{Massart} on the tailed
probabilities of the sample mean of Bernoulli random variable. For this reason, Stopping Rule B is also referred to as the stopping rule from
Massart's inequality  in \cite[Section 4.1.1]{Chen1241v16}.

\section{Double-Parabolic Sequential Estimation} \la{SecDPSE}

From Sections \ref{SubSecSR}, \ref{SubSRWI} and \ref{SubsecSRIC}, it can be seen that, by introducing a new parameter $\ro \in [0, 1]$ and
letting $\ro$ take values $\f{2}{3}$ and $0$ respectively, Stopping Rules B and D can be accommodated as special cases of the following general
stopping rule:

Continue the sampling process until \be \la{simplegeneral}
 \li (  \li | \wh{\bs{p}}_\ell - \f{1}{2}
\ri | - \ro  \vep \ri )^2 \geq \f{1}{4} + \f{ \vep^2 n_\ell}{2 \ln (\ze \de)} \ee for some $\ell \in \{1, 2, \cd, s \}$, where $\ze \in (0,
\f{1}{\de})$.

Moreover, as can be seen from (\ref{normalmassart}), the stopping rule derived from applying Wilson's confidence intervals to (\ref{ruleBI}) can
also be viewed as a special case of such general stopping rule with $\ro = 1$.

\begin{figure}[here]
\centerline{\psfig{figure=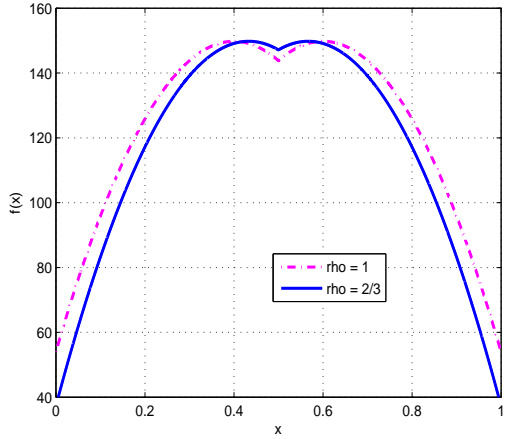, height=2.5in, width=3.0in }} \caption{Double-parabolic sampling}
\label{fig_CDMA}       
\end{figure}

From the stopping condition (\ref{simplegeneral}), it can be seen that the stopping boundary is associated with the double-parabolic function
$f(x) = \f{2}{\vep^2} \ln (\ze \de) \li [  \f{1}{4} -  \li (  \li | x - \f{1}{2} \ri | - \ro  \vep \ri )^2 \ri ]$ such that $x$ and $f(x)$
correspond to the sample mean and sample size respectively.  For $\vep = 0.1, \; \de = 0.05$ and $\ze = 1$, stopping boundaries with various
$\ro$  are shown by Figure \ref{fig_CDMA}.

For fixed $\vep$ and $\de$, the parameters $\ro$ and $\ze$ affect the shape of the stoping boundary in a way as follows.   As $\ro$ increases,
the span of stopping boundary is increasing in the axis of sample mean. By decreasing $\ze$, the stopping boundary can be dragged toward the
direction of increasing sample size. Hence, the parameter $\ro$ is referred to as the {\it dilation coefficient}.  The parameter $\ze$ is
referred to as the {\it coverage tuning parameter}. Since the stopping boundary consists of two parabolas, this approach of estimating a
binomial proportion is refereed to as the {\it double-parabolic sequential estimation} method.

\subsection{Parametrization of the Sampling Scheme} \la{SubsecSSS}

In this section, we shall parameterize the double-parabolic sequential sampling scheme by the method described in Section 2.2.    From the
stopping condition (\ref{simplegeneral}), the stopping rule can be restated as: Continue sampling until $\mscr{D}(\wh{\bs{p}}_\ell, n_\ell) = 1$
for some $\ell \in \{1, \cd, s \}$, where the function $\mscr{D} (.,.)$ is defined by \be \la{new D} \mscr{D} (z, n) = \bec 1 & \tx{if {\small
$(|z - \f{1}{2}| - \ro \vep)^2 \geq \f{1}{4} + \f{\vep^2 n}{2 \ln (\ze \de)}$}},\\
0  & \tx{otherwise} \eec  \ee  Clearly, the function $\mscr{D} (.,.)$ associated with the double-parabolic sequential sampling scheme depends on
the design parameters $\ro, \ze, \vep$ and $\de$.   Applying the function $\mscr{D} (.,.)$ defined by (\ref{new D}) to (\ref{def8a}) yields
{\small \be N_{\mrm{min}} = \min \li \{ n \in \bb{N}: \li ( \li | \f{k}{n} - \f{1}{2} \ri | - \ro \vep \ri )^2 \geq \f{1}{4} + \f{\vep^2 n}{2
\ln (\ze \de)} \; \tx{for some nonnegative integer $k$ not exceeding $n$} \ri \}. \la{call8} \ee} Since $\vep$ is usually small in practical
applications, we restrict $\vep$ to satisfy $0 < \ro \vep \leq \f{1}{4}$.  As a consequence of $0 \leq \ro \vep \leq \f{1}{4}$ and the fact that
$\li | z  - \f{1}{2} \ri | \leq \f{1}{2}$ for any $z \in [0, 1]$, it must be true that $\li ( \li | z - \f{1}{2} \ri | - \ro \vep \ri )^2 \leq
\li ( \f{1}{2} - \ro \vep \ri )^2$ for any $z \in [0, 1]$. It follows from (\ref{call8}) that $\li ( \f{1}{2} - \ro \vep \ri )^2 \geq \f{1}{4} +
\f{ \vep^2 N_{\mrm{min}} } { 2 \ln (\ze \de) }$, which implies that the minimum sample size can be taken as \be \la{keya}
 N_{\mrm{min}} = \li \lc 2 \ro \li ( \f{1}{ \vep
} - \ro \ri )  \ln \f{1}{\ze \de}  \ri \rc. \ee On the other hand, applying the function $\mscr{D} (.,.)$ defined by (\ref{new D}) to
(\ref{def8b}) gives {\small \be N_{\mrm{max}} = \min \li \{ n \in \bb{N}: \li ( \li | \f{k}{n} - \f{1}{2} \ri | - \ro \vep \ri )^2 \geq \f{1}{4}
+ \f{\vep^2 n}{2 \ln (\ze \de)} \; \tx{for all nonnegative integer $k$ not exceeding $n$} \ri \}. \la{call9} \ee} Since $\li ( \li | z -
\f{1}{2} \ri | - \ro \vep \ri )^2 \geq 0$ for any $z \in [0, 1]$, it follows from (\ref{call9}) that $\f{1}{4} + \f{ \vep^2 N_{\mrm{max}} } { 2
\ln (\ze \de) } \leq 0$, which implies that maximum sample size can be taken as \be \la{keyb} N_{\mrm{max}} = \li \lc \f{ 1 }{ 2 \vep^2} \ln
\f{1}{\ze \de} \ri \rc. \ee  Therefore, the sample sizes $n_1, \cd, n_s$ can be chosen as functions of $\ro, \ze, \vep$ and $\de$ which satisfy
the following constraint: \be \la{constraint}  N_{\mrm{min}} \leq n_1 < \cd < n_{s-1} <  N_{\mrm{max}} \leq n_s. \ee In particular, if the
number of stages $s$ is given and the group sizes are expected to be approximately equal, then the sample sizes, $n_1, \cd, n_s$,  for all
stages can be obtained by substituting $N_{\mrm{min}}$ defined by (\ref{keya}) and $N_{\mrm{max}}$ defined by (\ref{keyb}) into (\ref{equalss}).
For example, if the values of design parameters are $\vep = 0.05, \; \de = 0.05, \; \ro = \f{3}{4}, \; \ze = 2.6759$ and $s = 7$, then the
sample sizes of this sampling scheme are calculated as
\[ n_1 = 59, \; \; n_2 = 116, \; \; n_3 = 173, \; \; n_4 = 231, \; \; n_5 = 288, \; \; n_6 = 345, \; \; n_7 = 403.
\]
The stopping rule is completely determined by substituting the values of design parameters into (\ref{simplegeneral}).

\subsection{Uniform Controllability of Coverage Probability} \la{SubsecUCCP}

Clearly, for pre-specified  $\vep, \; \de$ and $\ro$, the coverage probability $\Pr \{  | \wh{\bs{p}} - p | < \vep \mid p \}$ depends on the
parameter $\ze$, the number of stages $s$, and the sample sizes $n_1, \cd,  n_s$.  As illustrated in Section \ref{SubsecSSS}, the number of
stages $s$ and the sample sizes $n_1, \cd, n_s$ can be defined as functions of $\ze \in (0, \f{1}{\de})$.  That is, the stopping rule can be
parameterized by $\ze$. Accordingly, for any $p \in (0, 1)$, the coverage probability $\Pr \{  | \wh{\bs{p}} - p | < \vep \mid p \}$ becomes a
function of $\ze$. The following theorem shows that it suffices to choose $\ze \in (0, \f{1}{\de})$ small enough to guarantee the pre-specified
confidence level.

\beT \la{UniformControl}
Let $\vep, \; \de \in (0, 1)$ and $\ro \in (0, 1]$ be fixed. Assume that the number of stages $s$ and the sample sizes
$n_1, \cd, n_s$ are functions of $\ze \in (0, \f{1}{\de})$ such that the constraint (\ref{constraint}) is satisfied. Then, $\Pr \{  |
\wh{\bs{p}} - p | < \vep \mid p \}$ is no less than $1 - \de$ for any $p \in (0, 1)$ provided that \[ 0 < \ze \leq \f{1}{\de} \exp \li (  \f{
\ln \f{\de}{2} + \ln \li [ 1 - \exp ( - 2  \vep^2 ) \ri ]}{ 4 \vep \ro ( 1 - \ro \vep) } \ri ).
\]
\eeT

See Appendix \ref{UniformControl_app} for a proof.  For Theorem \ref{UniformControl} to be valid, the choice of sample sizes is very flexible.
Specially, the sample sizes can be arithmetic or geometric progressions or any others, as long as the constraint (\ref{constraint}) is
satisfied. It can be seen that for the coverage probability to be uniformly controllable, the dilation coefficient $\ro$ must be greater than
$0$.   Theorem \ref{UniformControl} asserts that there exists $\ze > 0$ such that the coverage probability is no less than $1 - \de$, regardless
of the associated binomial proportion $p$. For the purpose of reducing sampling cost, we want to have a value of $\ze$ as large as possible such
that the pre-specified confidence level is guaranteed for any $p \in (0, 1)$. This can be accomplished by the technical components introduced in
Sections 2.1, 2.3, 2.4 and Section 2.5.  Clearly, for every value of $\ro$, we can obtain a corresponding value of $\ze$ (as large as possible)
to ensure the desired confidence level. However, the performance of resultant stopping rules are different.  Therefore, we can try a number of
values of $\ro$ and pick the best resultant stopping rule for practical use.

\subsection{Asymptotic Optimality of Sampling Schemes} \la{SubsecAOSS}

Now we shall provide an important reason why we propose the sampling scheme of that structure by showing its asymptotic optimality. Since the
performance of a group sampling scheme will be close to its fully sequential counterpart, we investigate the optimality of the fully sequential
sampling scheme. In this scenario, the sample sizes $n_1, n_2, \cd, n_s$ are consecutive integers such that \be \la{add338}
 \li \lc 2 \ro \li ( \f{1}{\vep}  -
\ro \ri ) \ln \f{1}{\ze \de}  \ri \rc = n_1 < n_2 < \cd < n_{s-1} < n_s = \li \lc \f{ 1 }{ 2 \vep^2} \ln \f{1}{\ze \de}  \ri \rc. \ee The fully
sequential sampling scheme can be viewed as a special case of a group sampling scheme of $s = n_s - n_1 + 1$ stages and group size $1$. Clearly,
if $\de, \; \ze$ and $\ro$ are fixed, the sampling scheme is dependent only on $\vep$.  Hence, for any $p \in (0, 1)$, if we allow $\vep$ to
vary in $(0, 1)$, then the coverage probability $\Pr \{ | \wh{\bs{p}} - p | < \vep \mid p \}$  and the average sample number $\bb{E} [ \mbf{n}
]$ are functions of $\vep$.  We are interested in knowing the asymptotic behavior of these functions as $\vep \to 0$, since $\vep$ is usually
small in practical situations. The following theorem provides us the desired insights.

\beT \la{AspOptimality} Assume that $\de \in (0, 1), \; \ze \in (0, \f{1}{\de})$ and $\ro \in (0, 1]$ are fixed.   Define $N (p, \vep, \de, \ze)
= \f{ 2 p(1 - p) \ln \f{1}{\ze \de} }{ \vep^2 }$ for $p \in (0, 1)$ and $\vep \in (0, 1)$.  Then, {\small \bel &  & \Pr \li \{ \lim_{\vep \to 0}
\f{ \mbf{n} } { N (p, \vep, \de, \ze) }
= 1 \mid p \ri \} = 1,  \nonumber \\
&  &  \lim_{\vep \to 0} \Pr \{ | \wh{\bs{p}} - p | < \vep \mid p \} = 2 \Phi \li ( \sq{ 2 \ln \f{1}{\ze \de} } \ri )  - 1, \la{asp8a}\\
&  &  \lim_{\vep \to 0} \f{ \bb{E} [ \mbf{n} ]   } { N (p, \vep, \de, \ze) } = 1 \la{asp8b} \eel} for any $p \in (0, 1)$.

\eeT

See Appendix \ref{AspOptimality_app} for a proof.  From (\ref{asp8a}), it can be seen that $\lim_{\vep \to 0} \Pr \{ | \wh{\bs{p}} - p | < \vep
\mid p \} = 1 - \de$ for any $p \in (0, 1)$ if $\ze  = \f{1}{\de} \exp  ( - \f{1}{2} \mcal{Z}_{\de \sh 2}^2)$.  Such value can be taken as an
initial value for the coverage tuning parameter $\ze$.   In addition to provide guidance on the coverage tuning techniques, Theorem
\ref{AspOptimality} also establishes the optimality of the sampling scheme.  To see this, let $\mscr{N} (p, \vep, \de)$ denote the minimum
sample size $n$ required for a fixed-sample-size procedure  to guarantee that $\Pr \{ | \ovl{X}_n - p | < \vep \mid p \} \geq 1 - \de$ for any
$p \in (0, 1)$, where $\ovl{X}_n = \f{\sum_{i = 1}^n X_i  }{n}$.  It is well known that from the central limit theorem, \be \la{appss}
 \lim_{\vep \to 0} \f{ \mscr{N} (p, \vep, \de) }{  p(1 - p)  \li ( \f{  \mcal{Z}_{\de
\sh 2} } { \vep } \ri )^2 } = 1. \ee  Making use of  (\ref {asp8b}), (\ref{appss}) and letting $\ze  = \f{1}{\de} \exp  ( - \f{1}{2}
\mcal{Z}_{\de \sh 2}^2)$, we have $\lim_{\vep \to 0} \f{ \mscr{N} (p, \vep, \de) } { N (p, \vep, \de, \ze) } = 1$ for $p \in (0, 1)$ and $\de
\in (0, 1)$, which implies the asymptotic optimality of the double-parabolic sampling scheme.  By virtue of (\ref {asp8b}), an approximate
formula for computing the average sample number is given as \be \la{ASNAP}
 \bb{E} [ \mbf{n} ] \ap N (p, \vep, \de, \ze) = \f{ 2 p(1 - p) \ln \f{1}{\ze \de} }{ \vep^2 }
\ee for $p \in (0, 1)$ and $\vep \in (0, 1)$.  From (\ref{appss}), one obtains $\mscr{N} (p, \vep, \de) \ap p(1 - p)  \li ( \f{  \mcal{Z}_{\de
\sh 2} } { \vep } \ri )^2$,  which is a well-known result in statistics.  In situations that no information of $p$ is available,  one usually
uses \be \la{SSNAP}
 N_{\mrm{normal}}  \DEF \li \lc \f{1}{4} \li ( \f{  \mcal{Z}_{\de \sh 2} } { \vep } \ri )^2 \ri \rc
\ee as the sample size for estimating the binomial proportion $p$ with prescribed margin of error $\vep$ and confidence level $1- \de$.  Since
the sample size formula (\ref{SSNAP}) can lead to under-coverage, researchers in many areas are willing to use a more conservative but rigorous
sample size formula \be \la{SSNCH}
 N_{\mrm{ch}}  \DEF \li \lc \f{ \ln \f{2}{\de} }{ 2 \vep^2} \ri \rc,
\ee which is derived from the Chernoff-Hoeffding bound \cite{Chernoff, Hoeffding}.  Comparing (\ref{ASNAP}) and (\ref{SSNCH}), one can see that
under the premise of guaranteeing the prescribed confidence level $1 - \de$, the double-parabolic sampling scheme can lead to a substantial
reduction of sample number when the unknown binomial proportion $p$ is close to $0$ or $1$.

\subsection{Bounds on Distribution and Expectation of Sample Number}  \la{SubsecBDESN}

We shall derive analytic bounds for the cumulative distribution function and expectation of the sample number $\mbf{n}$ associated with the
double-parabolic sampling scheme. In this direction, we have obtained the following results.  \beT \la{bounds} Let $p \in (0, \f{1}{2}]$.
Define $a_\ell = \f{1}{2} - \ro \vep - \sq{\f{1}{4} + \f{ \vep^2 n_\ell } { 2 \ln (\ze \de) }}$ for $\ell = 1, \cd, s$. Let $\tau$ denote the
index of stage such that $a_{\tau - 1} \leq p < a_\tau$. Then, $\Pr \{ \mbf{n}
> n_\ell \mid p \} \leq \exp (n_\ell \mscr{M} (a_\ell,p))$ for $\tau \leq \ell < s$. Moreover,  $\bb{E} [ \mbf{n} ] \leq n_\tau + \sum_{\ell =
\tau}^{s-1} ( n_{\ell + 1} - n_\ell ) \exp (n_\ell \mscr{M} (a_\ell,p))$.   \eeT

See Appendix \ref{App_bounds} for a proof.  By the symmetry of the double-parabolic sampling scheme, similar analytic bounds for the
distribution and expectation of the sample number can be derived for the case that $p \in [\f{1}{2}, 1)$.

\section{Comparison of Computational Methods} \la{SecCI}

In this section, we shall compare various computational methods.   First, we will illustrate why a frequently-used method of evaluating the
coverage probability based on gridding the parameter space is not rigorous and is less efficient as compared to the Adapted B\&B Algorithm.
Second, we will introduce the Adaptive Maximum Checking Algorithm of \cite{Chen1241v16} which has better computational efficiency as compared to
the Adapted B\&B Algorithm. Third, we will explain that it is more advantageous in terms of numerical accuracy to work with the complementary
coverage probability as compared to direct evaluation of the coverage probability. Finally, we will compare the computational methods of Chen
\cite{Chen1241v1, Chen1241v12, Chen4679v1, Chen0430v1, Chen3458v2} and Frey \cite{Frey} for the design of sequential procedures for estimating a
binomial proportion.

\subsection{Verifying Coverage Guarantee without Gridding Parameter Space} \la{SubsecVCGGPS}

For purpose of constructing a sampling scheme so that the prescribed confidence level $1 - \de$ is guaranteed, an essential task is to determine
whether the coverage probability $\Pr \{   | \wh{\bs{p}} - p |  < \vep \mid p \}$ associated with a given stopping rule is no less than $1 -
\de$. In other words, it is necessary to compare the infimum of coverage probability with $1 - \de$.  To accomplish such a task of checking
coverage guarantee, a natural method is to evaluate the infimum of coverage probability as follows:

{\bf (i)}: Choose $m$ grid points $p_1, \cd, p_m$ from parameter space $(0, 1)$.

{\bf (ii)}: Compute $c_j = \Pr \{   | \wh{\bs{p}} - p |  < \vep \mid p_j \}$ for $j = 1, \cd, m$.

{\bf (iii)}: Take $\min \{ c_1, \cd, c_m \} $ as $\inf_{p \in (0, 1)} \Pr \{   | \wh{\bs{p}} - p |  < \vep \mid p \}$.

This method  can be easily mistaken as an exact approach and has been frequently used for evaluating coverage probabilities in many problem
areas.

It is not hard to show that if the sample size $\mbf{n}$ of a sequential procedure has a support $\mscr{S}$, then the coverage probability $\Pr
\{ | \wh{\bs{p}} - p | < \vep  \mid p \}$ is discontinuous at $p \in \mscr{P} \cap (0, 1)$, where $\mscr{P} = \{ \f{k}{n} \pm \vep : \tx{$k$ is
a nonnegative integer}$ $\tx{no greater than $n \in \mscr{S}$} \}$.  The set $\mscr{P}$ typically has a large number of parameter values. Due to
the discontinuity of the coverage probability as a function of $p$, the coverage probabilities can differ significantly for two parameter values
which are extremely close. This implies that an intolerable error can be introduced by taking the minimum of coverage probabilities of a finite
number of parameter values as the infimum of coverage probability on the whole parameter space. So, if one simply uses the minimum of the
coverage probabilities of a finite number of parameter values as the infimum of coverage probability to check the coverage guarantee, the
sequential estimator $\wh{\bs{p}}$ of the resultant stopping rule will fail to guarantee the prescribed confidence level.

In addition to the lack of rigorousness, another drawback of checking coverage guarantee based on the method of gridding parameter space is its
low efficiency.  A critical issue is on the choice of the number, $m$, of grid points.  If the number $m$ is too small, the induced error can be
substantial. On the other hand, choosing a large number for $m$  results in high computational complexity.

In contrast to the method based on gridding parameter space, the Adapted B\&B Algorithm is a rigorous approach for checking coverage guarantee
as a consequence of the mechanism for comparing the bounds of coverage probability with the prescribed confidence level.  The algorithm is also
efficient due to the mechanism of pruning branches.

\subsection{Adaptive Maximum Checking Algorithm} \la{SubsecAMCA}

As illustrated in Section \ref{SecHow}, the techniques developed in \cite{Chen1241v1, Chen1241v12,  Chen4679v1, Chen0430v1, Chen3458v2} are
sufficient to provide exact solutions for a wide range of sequential estimation problems.  However, one of the four components, the Adapted B\&B
Algorithm,  requires computing both the lower and upper bounds of the complementary coverage probability.  To further reduce the computational
complexity, it is desirable to have a checking algorithm which needs only one of the lower and upper bounds. For this purpose, Chen had
developed the Adaptive Maximum Checking Algorithm (AMCA) in  \cite[Section 3.3]{Chen1241v16} and \cite[Section 2.7]{Chen_SPIE}.   In the
following introduction of the AMCA, we shall follow the description of \cite{Chen1241v16}.  The AMCA can be applied to a wide class of
computational problems dependent on the following critical subroutine:

Determine whether a function $C(\se)$ is smaller than a prescribed number $\de$ for every value of $\se$ contained in interval $[ \udl{\se},
\ovl{\se}]$.

Specially, for checking the coverage guarantee in the context of estimating a binomial proportion, the parameter $\se$ is the binomial
proportion $p$ and the function $C(\se)$ is actually the complementary coverage probability. In many situations, it is impossible or very
difficult to evaluate $C(\se)$ for every value of $\se$ in interval $[ \udl{\se}, \ovl{\se}]$, since the interval may contain infinitely many or
an extremely large number of values. Similar to the Adapted B\&B Algorithm, the purpose of AMCA is to reduce the computational complexity
associated with the problem of determining whether the maximum of $C(\se)$ over $[ \udl{\se}, \ovl{\se}]$ is less than $\de$. The only
assumption required for AMCA is that, for any interval $[a, b] \subseteq [ \udl{\se}, \ovl{\se}]$, it is possible to compute an upper bound
$\ovl{C} (a, b)$ such that $C(\se) \leq \ovl{C} (a, b)$ for any $\se \in [a, b]$ and that the upper bound converges to $C(\se)$ as the interval
width $b - a$ tends to $0$. The backward AMCA proceeds as follows:

\bsk

\begin{tabular} {|l |}
\hline
$ \nabla  \; \tx{Choose initial step size $d > \eta$}$.\\
$ \nabla \; \tx{Let $F \leftarrow 0, \; T \leftarrow 0$ and $b \leftarrow \ovl{\se}$}$.\\
$ \nabla \; \tx{While $F = T = 0$, do the following}$:\\
$ \indent \indent \; \diamond  \; \tx{Let $\tx{st} \leftarrow 0$ and $\ell \leftarrow 2$}$;\\
$ \indent \indent \; \diamond \; \tx{While $\tx{st} = 0$, do the following}$:\\
$\indent \indent \indent \indent \; \star \; \tx{Let $\ell \leftarrow \ell - 1$ and $d \leftarrow d  2^\ell$}$.\\
$ \indent \indent \indent \indent \; \star \; \tx{If $b - d > \udl{\se}$, then let $a \leftarrow b - d$ and $T \leftarrow 0$}$.\\
$ \indent \indent \indent \indent \; \qu \; \; \tx{Otherwise, let $a \leftarrow \udl{\se}$ and $T \leftarrow 1$}$.\\
$\indent \indent \indent \indent \; \star \; \tx{If $\ovl{C} (a, b) < \de$, then let $\tx{st} \leftarrow 1$ and $b \leftarrow a$}$.\\
$ \indent \indent \indent \indent \; \star \; \tx{If $d < \eta$, then let $\tx{st} \leftarrow 1$ and $F \leftarrow 1$}$.\\
$ \nabla \; \tx{Return $F$}$.
\\ \hline
\end{tabular}

\bsk

The output of the backward AMCA is a binary variable $F$ such that ``$F = 0$'' means ``$C(\se) < \de$'' and ``$F = 1$'' means ``$C(\se) \geq
\de$''. An intermediate variable $T$ is introduced in the description of AMCA such that ``$T = 1$'' means that the left endpoint of the interval
is reached. The backward AMCA starts from the right endpoint of the interval (i.e., $b = \ovl{\se}$) and attempts to find an interval $[a, b]$
such that $\ovl{C} (a, b) < \de$. If such an interval is available, then, attempt to go backward to find the next consecutive interval with
twice width. If doubling the interval width fails to guarantee $\ovl{C} (a, b) < \de$, then try to repeatedly cut the interval width in half to
ensure that $\ovl{C} (a, b) < \de$. If the interval width becomes smaller than a prescribed tolerance $\eta$, then AMCA declares that ``$F =
1$''. For our relevant statistical problems, if $C(\se) \geq \de$ for some $\se \in [\udl{\se}, \ovl{\se}]$, it is sure that ``$F = 1$'' will be
declared.  On the other hand, it is possible that ``$F = 1$'' is declared even though $C(\se) < \de$ for any $\se \in [\udl{\se}, \ovl{\se}]$.
However, such situation can be made extremely rare and immaterial if we choose $\eta$ to be a very small number.  Moreover, this will only
introduce negligible conservativeness in the evaluation of $C(\se)$ if $\eta$ is chosen to be sufficiently small (e.g., $\eta = 10^{-15}$).
Clearly, the backward AMCA can be easily modified as forward AMCA. Moreover, the AMCA can also be easily modified as  Adaptive Minimum Checking
Algorithm (forward and backward).  For checking the maximum of complementary coverage probability $\Pr \{ | \wh{\bs{p}} - p | \geq \vep \mid p
\}$, one can use the AMCA with $C(p) = \Pr \{ | \wh{\bs{p}} - p | \geq \vep \mid p \}$ over interval $[0, \f{1}{2}]$.   We would like to point
out that, in contrast to the Adapted B\&B Algorithm,  it seems difficult to generalize the AMCA to problems involving multidimensional parameter
spaces.

\subsection{Working with Complementary Coverage Probability} \la{SubsecWCCP}

We would like to point out that, instead of evaluating the coverage probability as in \cite{Frey}, it is better to evaluate the complementary
coverage probability for purpose of reducing numerical error.  The advantage of working on the complementary coverage probability can be
explained as follows: Note that, in many cases, the coverage probability is very close to $1$ and the complementary coverage probability is very
close to $0$. Since the absolute precision for computing a number close to $1$ is much lower than the absolute precision for computing a number
close to $0$, the method of directly evaluating the coverage probability will lead to intolerable numerical error for problems involving small
$\de$.  As an example, consider a situation that the complementary coverage probability is in the order of $10^{-5}$. Direct computation of the
coverage probability can easily lead to an absolute error of the order of $10^{-5}$. However, the absolute error of computing the complementary
coverage probability can be readily controlled at the order of $10^{-9}$.

\subsection{Comparison of Approaches of Chen and J. Frey}

As mentioned in the introduction, J. Frey published a paper \cite{Frey} in {\it The American Statistician} (TAS) on the sequential estimation of
a binomial proportion with prescribed margin of error and confidence level.  The approaches of Chen and Frey are based on the same strategy as
follows:  First, construct a family of stopping rules parameterized by $\ga$ (and possibly other design parameters) so that the associated
coverage probability $\Pr \{ |\wh{\bs{p}} - p | < \vep \mid p \}$ can be controlled by parameter $\ga$ in the sense that the coverage
probability  can be made arbitrarily close to $1$ by increasing $\ga$.  Second, adaptively and rigorously check the coverage guarantee by virtue
of bounds of coverage probabilities. Third, apply a bisection search method to determine the parameter $\ga$ so that the coverage probability is
no less than the prescribed confidence level $1 - \de$ for any $p \in (0, 1)$.

For the purpose of controlling the coverage probability,  Frey \cite{Frey} applied the inclusion principle previously proposed in \cite[Section
3]{Chen3458v2} and used in \cite{Chen1241v1, Chen1241v12,  Chen4679v1, Chen0430v1}.  As illustrated in Section 3, the central idea of inclusion
principle is to use a sequence of confidence intervals to construct stopping rules so that the sampling process is continued until a confidence
interval is included by an interval defined in terms of the estimator and margin of error. Due to the inclusion relationship, the associated
coverage probability can be controlled by the confidence coefficients of the sequence of confidence intervals. The critical value $\ga$ used by
Frey plays the same role for controlling coverage probabilities as that of the coverage tuning parameter $\ze$ used by Chen. Frey \cite{Frey}
stated stopping rules in terms of confidence limits. This way of expressing stopping rules is straightforward and insightful, since one can
readily seen the principle behind the construction.  For convenience of practical use, Chen proposed to eliminate the necessity of computing
confidence limits.

Frey's method for checking coverage guarantee differs from the Adapted B\&B Algorithm, but coincides with other techniques of Chen
\cite{Chen1241v16}. On September 18, 2011, in response to an inquiry on the coincidence of the research results, Frey simultaneously emailed
Xinjia Chen (the coauthor of the present paper) and TAS Editor John Stufken all pre-final revisions of his manuscript for the paper \cite{Frey}.
In his original manuscript submitted to TAS in July 2009, Frey's method was to ``simply approximate $CP(\ga)$ by taking the minimum over the
grid of values $p = 1/2001, . . . , 2000/2001$.'' In the first revision of his manuscript submitted to TAS in November 2009, Frey's method was
to ``approximate $CP(\ga)$ by taking the minimum of $T(p; \ga)$ over the grid of values $p = 1/2001, . . . , 2000/2001$ and the set of values of
the form $p = c \pm \ep$, where $c \in C$ and $\ep = 10^{-10}$.''    In Frey's notational system, $\ga$ is the critical value which plays the
same role as that of the coverage tuning parameter $\ze$ in the present paper, $T(p; \ga)$ is the coverage probability, $CP(\ga)$ is the infimum
of coverage probability for $p \in (0,1)$, and $C = \{ \wh{p} \pm \vep: \tx{ $\wh{p}$ is a possible value of $\wh{\bs{p}}$} \} \cap (0, 1)$.
From the original and the first revision of his manuscript submitted to TAS before April 2010, it can be seen that Frey's method of checking
coverage guarantee was dependent on taking the minimum of coverage probabilities for a finite number of gridding points of $p \in (0, 1)$ as the
infimum coverage probability for $p \in (0, 1)$.   As can be seen from Section 5.1 of the present paper,  such method lacks rigorousness and
efficiency. In the second revision of his manuscript submitted to TAS in April 2010, for the purpose of checking coverage guarantee, Frey
replaced the method of gridding parameter space with an interval bounding technique and proposed a checking algorithm which is essentially the
same as the AMCA precedentially established by Chen \cite[Section 3.3]{Chen1241v16} in November 2009.

Similar to the AMCA of \cite[Section 3.3]{Chen1241v16}, the algorithm of Frey \cite[Appendix]{Frey} for checking coverage guarantee adaptively
scans the parameter space based on interval bounding. The adaptive method used by Frey for updating step size is essentially the same as that of
the AMCA. Ignoring the number $0.01$ in Frey's expression ``$\ep_i = \min \{0.01, 2(p_{i-1} - p_{i-2})\}$'', which has very little impact on the
computational efficiency, Frey's step size $\ep_i$ can be identified as the adaptive step size $d$ in the AMCA. The operation associated with
``$\ep_i = \min \{0.01, 2(p_{i-1} - p_{i-2})\}$'' has a similar function as that of the command ``Let $\tx{st} \leftarrow 0$ and $\ell
\leftarrow 2$'' in the outer loop of the AMCA.  The operation associated with Frey's expression ``$p_{i-1} + \ep_i/2^j, \;  j \geq 0$'' is
equivalent to that of the command ``Let $\ell \leftarrow \ell - 1$ and $d \leftarrow d 2^\ell$'' in the inner loop of the AMCA. Frey proposed to
declare a failure of coverage guarantee if ``the distance from $p_{i-1}$ to the candidate value for $p_i$ falls below $10^{-14}$''. The number
``$10^{-14}$'' actually plays the same role as ``$\eta$'' in the AMCA, where ``$\eta = 10^{-15}$'' is recommended by \cite{Chen1241v16}.

\sect{Numerical Results} \la{SubsecNR}

In this section, we shall illustrate the proposed double-parabolic sampling scheme through examples.  As demonstrated in Section 2.2 and Section
4, the double-parabolic sampling scheme can be parameterized by the dilation coefficient $\ro$ and the coverage tuning parameter $\ze$. Hence,
the performance of the resultant stopping rule can be optimized with respect to $\ro \in (0, 1]$ and $\ze$ by choosing various values of $\ro$
from interval $(0, 1]$ and determining the corresponding values of $\ze$ by the computational techniques introduced in Section 2 to guarantee
the desired confidence interval.

\subsection{Asymptotic Analysis May Be Inadequate} \la{SubsecAAM}

\begin{figure}[here]
\includegraphics{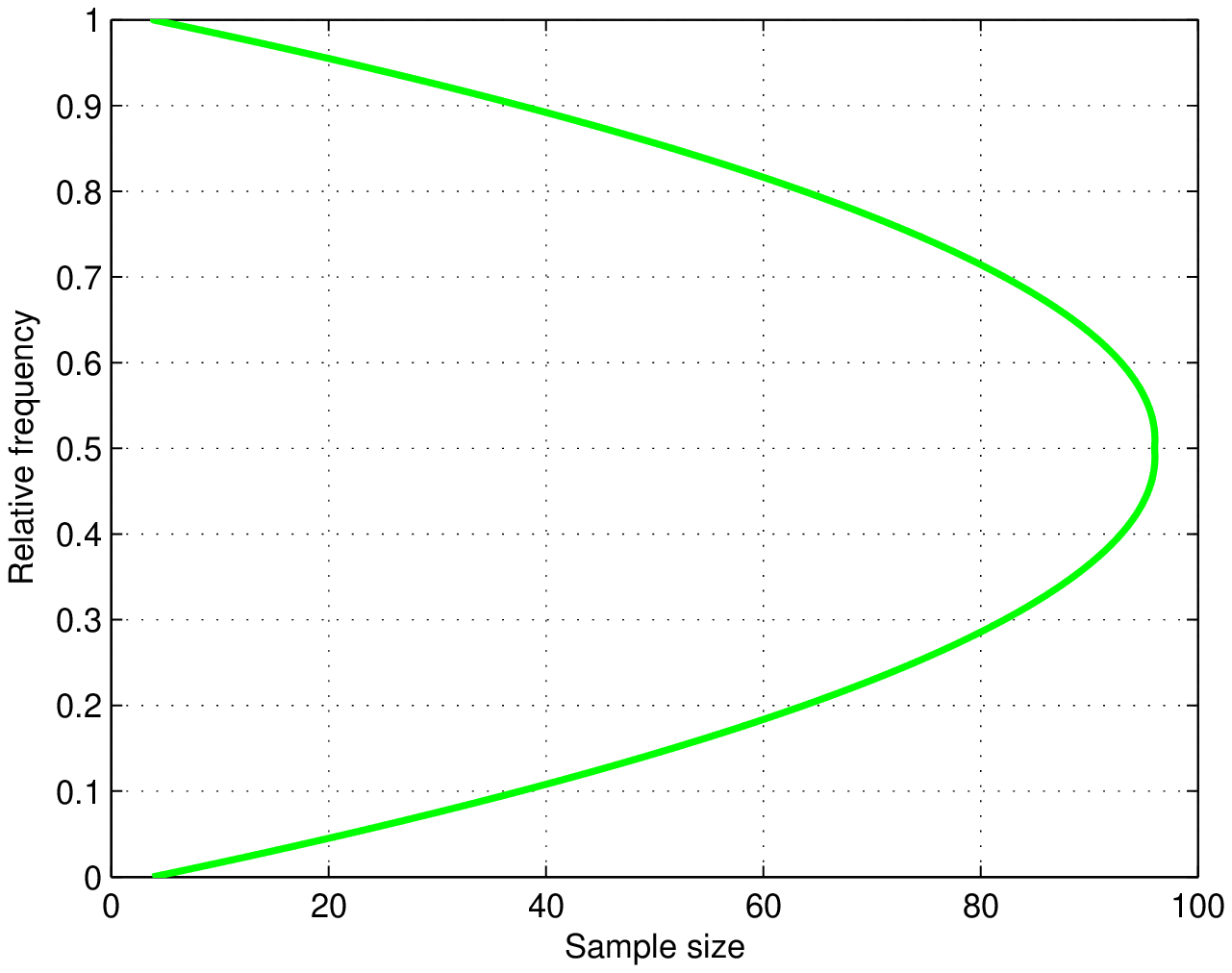} \includegraphics{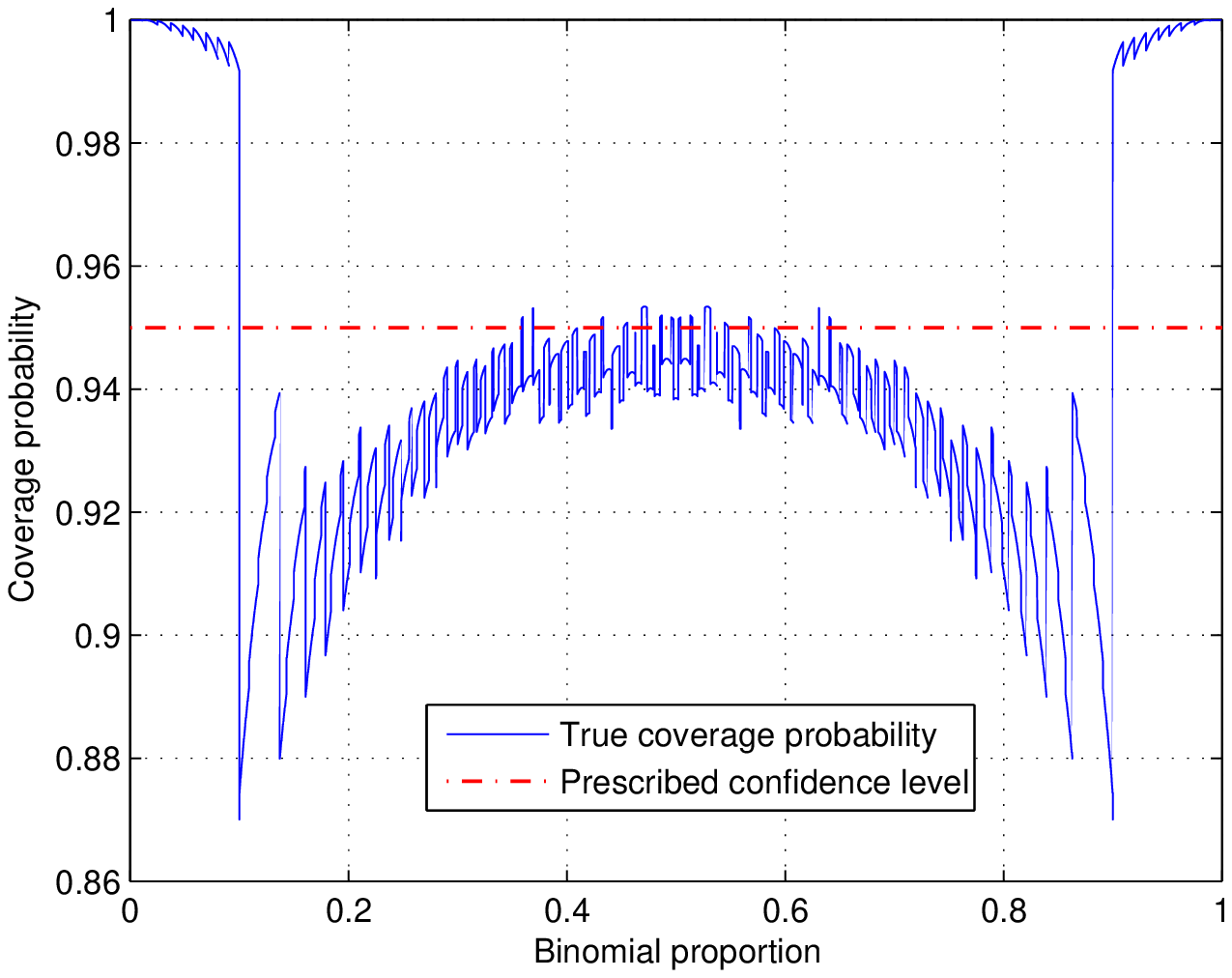} \caption{Double-parabolic sampling with $\vep = 0.1, \; \de = 0.05, \; \ro = \f{1}{10}$ and $\ze = 2.93$ } \label{fig_ST_Asp}
\end{figure}

For fully sequential cases, we have evaluated the double-parabolic sampling scheme with $\vep = 0.1, \; \de = 0.05, \; \ro = 0.1$ and {\small
$\ze = \f{1}{\de} \exp \li ( - \f{1}{2} \mcal{Z}_{\de\sh2}^2 \ri ) \ap 2.93$}.   The stopping boundary is displayed in the left side of Figure
\ref{fig_ST_Asp}. The function of coverage probability with respect to the binomial proportion is shown in the right side of Figure
\ref{fig_ST_Asp}, which indicates that the coverage probabilities are generally substantially lower than the prescribed confidence level $1 -
\de = 0.05$.  By considering $\vep = 0.1$ as a small number and applying the asymptotic theory, the coverage probability associated with the
sampling scheme is expected to be close to $0.95$. This numerical example demonstrates that although the asymptotic method is insightful and
involves virtually no computation, it may not be adequate.

In general, the main drawback of an asymptotic method is that there is no guarantee of coverage probability.  Although an asymptotical method
asserts that if the margin of error $\vep$ tends to $0$, the coverage probability will tend to the pre-specified confidence level $1 - \de$, it
is difficult to determine how small the margin of error $\vep$ is sufficient for the asymptotic method to be applicable. Note
 that $\vep \to 0$ implies the average sample size tends to $\iy$.  However, in reality, the sample sizes must be finite.
 Consequently, an asymptotic method inevitably introduces unknown statistical error.
 Since an asymptotic method does not necessarily guarantee the prescribed confidence level, it is not fair to compare its associated sample size with that
of  an exact method, which guarantees the pre-specified confidence level.

This example also indicates that, due to the discrete nature of the problem,  the coverage probability is a discontinuous and erratic function
of $p$, which implies that Monte Carlo simulation is not suitable for evaluating the coverage performance.

\bsk  \bsk \bsk \bsk \bsk

\vspace*{2.5in}

\subsection{Parametric Values of Fully Sequential Schemes} \la{SubsecPVFSS}

For fully sequential cases, to allow direct application of our double-parabolic sequential method, we have obtained values of coverage tuning
parameter $\ze$, which guarantee the prescribed confidence levels,  for  double-parabolic  sampling schemes with $\rho = \f{3}{4}$ and various
combinations of $(\vep, \de)$ as shown in Table \ref{table_abs}.  We used the computational techniques introduced in Section 2 to obtain this
table.

\begin{table}[here]
\caption{Coverage Tuning Parameter} \label{table_abs}
\begin{center}
\begin{tabular}{|c||c||c|||c||c||c|||c||c||c|}
\hline $\vep$ & $\de$ & $\zeta$ & $\vep$ & $\de$ & $\zeta$ & $\vep$ & $\de$
& $\zeta$\\
\hline
\hline $0.1$ & $0.1$ & $2.0427$ & $0.1$ & $0.05$ & $ 2.4174$ & $0.1$ & $0.01$ & $3.0608$\\
\hline $0.05$ & $0.1$ & $2.0503$ & $0.05$ & $0.05$ & $2.5862$ & $0.05$ & $0.01$ & $3.3125$\\
\hline $0.02$ & $0.1$ & $2.1725$ & $0.02$ & $0.05$ & $2.5592$ & $0.02$ & $0.01$ & $3.4461$\\
\hline $0.01$ & $0.1$ & $2.1725$ & $0.01$ & $0.05$ & $2.5592$ & $0.01$ & $0.01$ & $3.4461$\\
 \hline
\end{tabular}
\end{center}
\end{table}

\begin{figure}[here]
\includegraphics{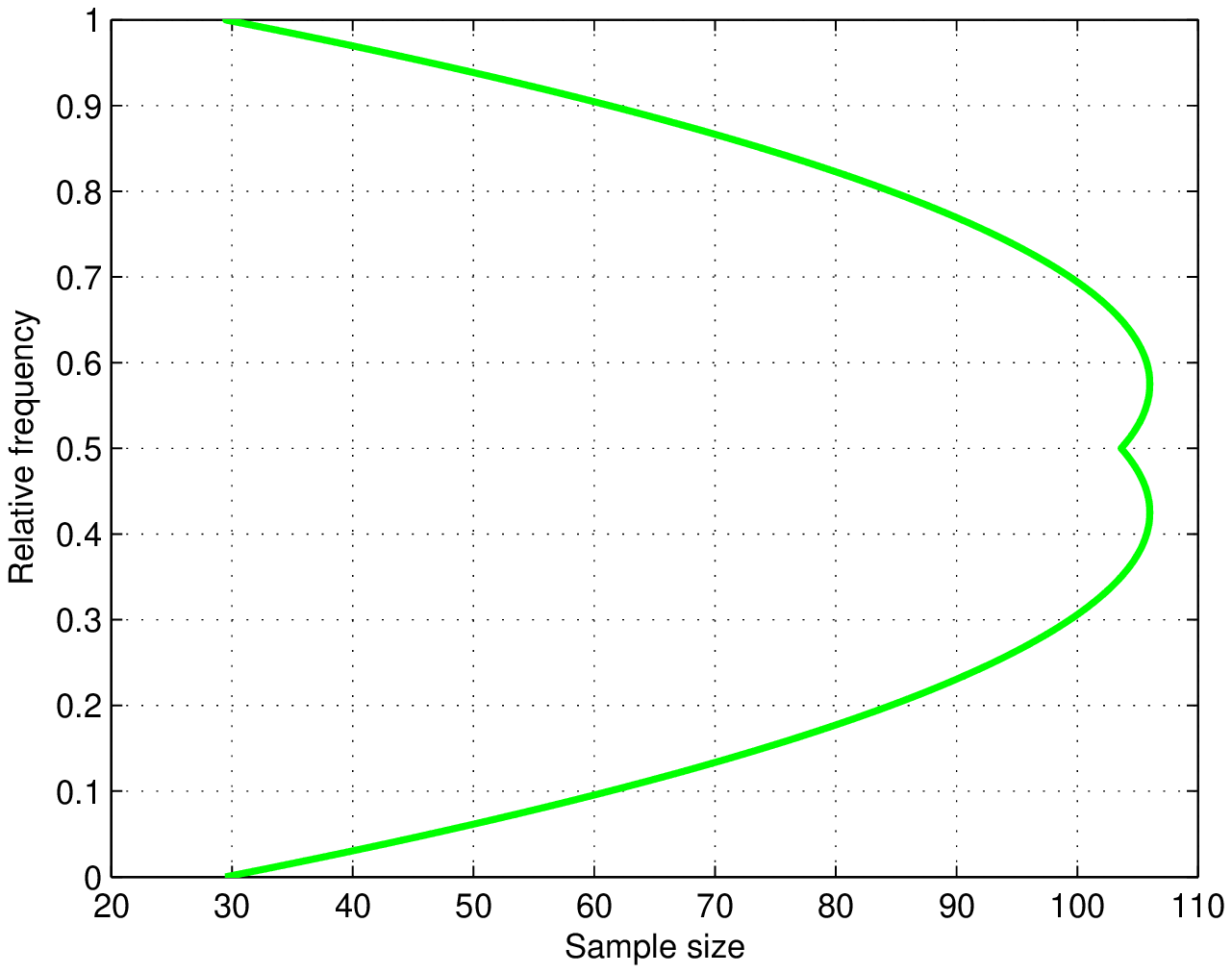} \includegraphics{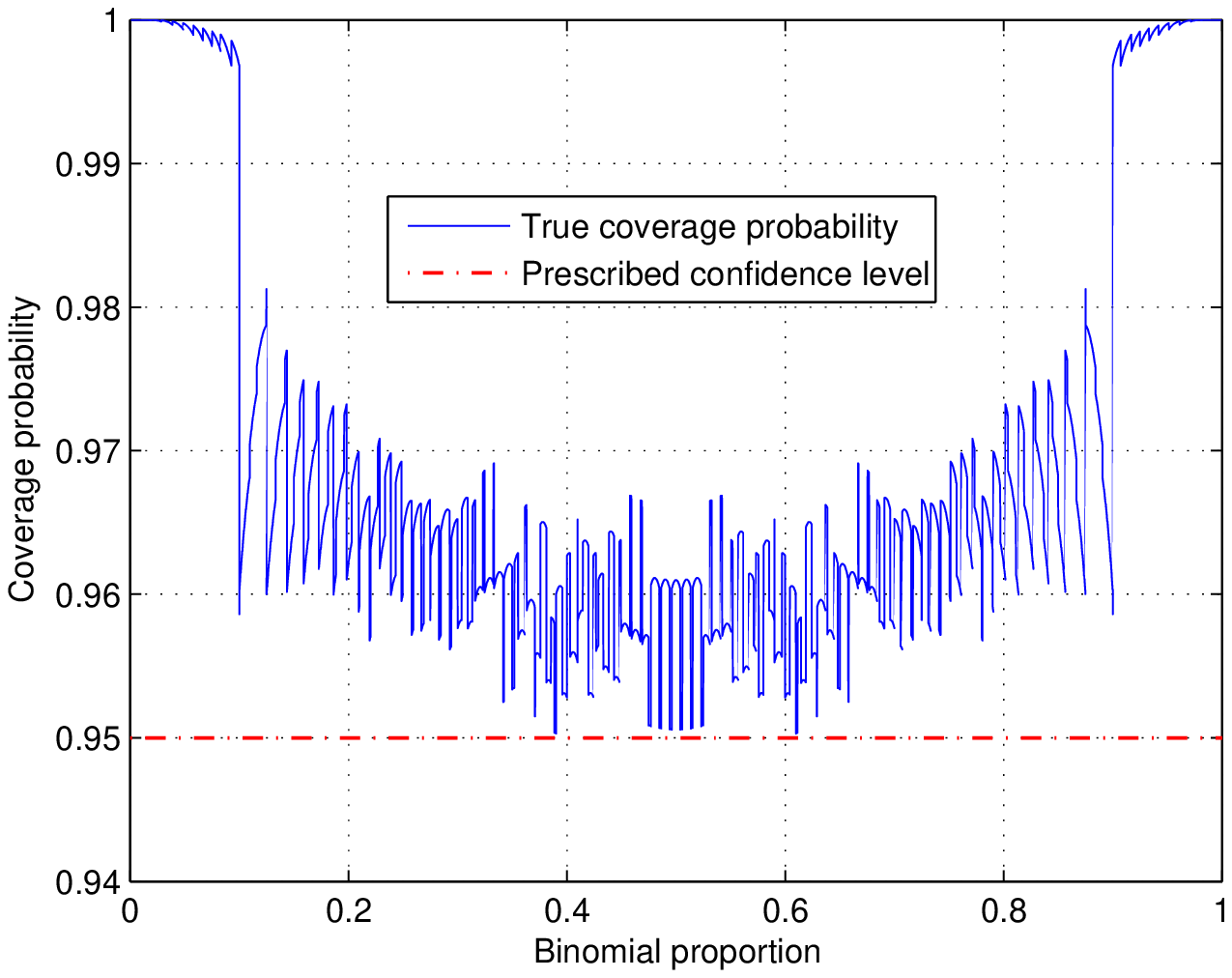} \caption{Double-parabolic sampling with $\vep = 0.1, \; \de = 0.05, \; \ro = \f{3}{4}$ and $\ze = 2.4174$}
\label{fig_ST_DP}
\end{figure}

To illustrate the use of Table \ref{table_abs}, suppose that one wants a fully sequential sampling procedure to ensure that $\Pr \{ |
\wh{\bs{p}} - p | < 0.1 \mid p \} > 0.95$ for any $p \in (0, 1)$.  This means that one can choose $\vep = 0.1, \; \de = 0.05$ and the range of
sample size is given by (\ref{add338}).  From Table \ref{table_abs}, it can be seen that the value of $\ze$ corresponding to $\vep = 0.1, \; \de
= 0.05$ is $2.4174$. Consequently, the stopping rule is completely determined by substituting the values of design parameters $\vep = 0.1, \;
\de = 0.05, \; \ro = \f{3}{4}, \; \ze = 2.4174$ into its definition. The stopping boundary of this sampling scheme is displayed in the left side
of Figure \ref{fig_ST_DP}. The function of coverage probability with respect to the binomial proportion is shown in the right side of Figure
\ref{fig_ST_DP}.

\bsk \bsk \bsk

\vspace*{2.7in}

\subsection{Parametric Values of Group Sequential Schemes} \la{SubsecPVGSS}

In many situations, especially in clinical trials, it is desirable to use group sequential sampling schemes.  In Tables \ref{table_abs_group05}
and \ref{table_abs_group01}, assuming that sample sizes satisfy (\ref{equalss}) for the purpose of having approximately equal group sizes, we
have obtained parameters for concrete schemes by the computational techniques introduced in Section 2.

For dilation coefficient $\rho = \f{3}{4}$ and confidence parameter $\de = 0.05$, we have obtained values of coverage tuning parameter $\ze$,
which guarantee the prescribed confidence level $0.95$,  for double-parabolic sampling schemes, with the number of stages $s$ ranging from $3$
to $10$, as shown in Table \ref{table_abs_group05}.

For dilation coefficient $\rho = \f{3}{4}$ and confidence parameter $\de = 0.01$, we have obtained values of coverage tuning parameter $\ze$,
which guarantee the prescribed confidence level $0.99$,  for double-parabolic sampling schemes, with the number of stages $s$ ranging from $3$
to $10$, as shown in Table \ref{table_abs_group01}.

\begin{table}[here]
\caption{Coverage Tuning Parameter} \label{table_abs_group05}
\begin{center}
\begin{tabular}{|c||c|c|c|c|c|c|c|c|}
\hline $ $ & $s = 3$ & $s = 4$ & $s = 5$ & $s = 6$ & $s = 7$ & $s = 8$ & $s = 9$ & $s = 10$\\
\hline
\hline $\vep = 0.1$ &  $2.6583$  &  $2.6583$  &   $2.5096$  &  $2.5946$  &   $2.4459$ &  $2.6512$  &   $2.5096$  &  $2.4459$\\
\hline $\vep = 0.05$ & $2.6759$  &  $2.6759$  &  $2.6759$  &  $2.6759$  &  $2.6759$  &  $2.6759$  &  $2.6759$   &  $2.6759$\\
\hline $\vep = 0.02$ & $2.6725$  &  $2.6725$  &  $2.6725$  &  $2.6725$  &  $2.6725$  &  $2.6725$  &  $2.6725$  &  $2.6725$\\
\hline $\vep = 0.01$ & $2.6796$  &  $2.6796$  &  $2.6796$  &  $2.6796$  &  $2.6796$  &  $2.5875$  &  $2.6796$  &  $2.6796$\\
 \hline
\end{tabular}
\end{center}
\end{table}

\begin{table}[here]
\caption{Coverage Tuning Parameter} \label{table_abs_group01}
\begin{center}
\begin{tabular}{|c||c|c|c|c|c|c|c|c|}
\hline $ $ & $s = 3$ & $s = 4$ & $s = 5$ & $s = 6$ & $s = 7$ & $s = 8$ & $s = 9$ & $s = 10$\\
\hline
\hline $\vep = 0.1$ &  $3.3322$  &  $3.3322$  &  $3.3322$  &   $3.3322$  &  $3.3322$  &  $3.2709$  &  $3.0782$  &   $3.3322$\\
\hline $\vep = 0.05$  &  $3.5074$  &  $3.5074$  &  $3.5074$  &  $3.5074$  &  $3.5074$  &  $3.5074$  &  $3.5074$  & $3.5074$\\
\hline $\vep = 0.02$ &  $3.5430$  &  $3.5430$  &  $3.5430$  &  $3.5430$  &  $3.5430$  &  $3.5430$  &  $3.5430$  &  $3.5430$\\
\hline $\vep = 0.01$ &  $3.5753$  &  $3.5753$  &  $3.5753$  &  $3.5753$  &  $3.5753$  &  $3.5753$  &  $3.5753$  &  $3.5753$\\
 \hline
\end{tabular}
\end{center}
\end{table}

\bsk \bsk

\vspace*{2.8in}
\begin{figure}[here]
\includegraphics{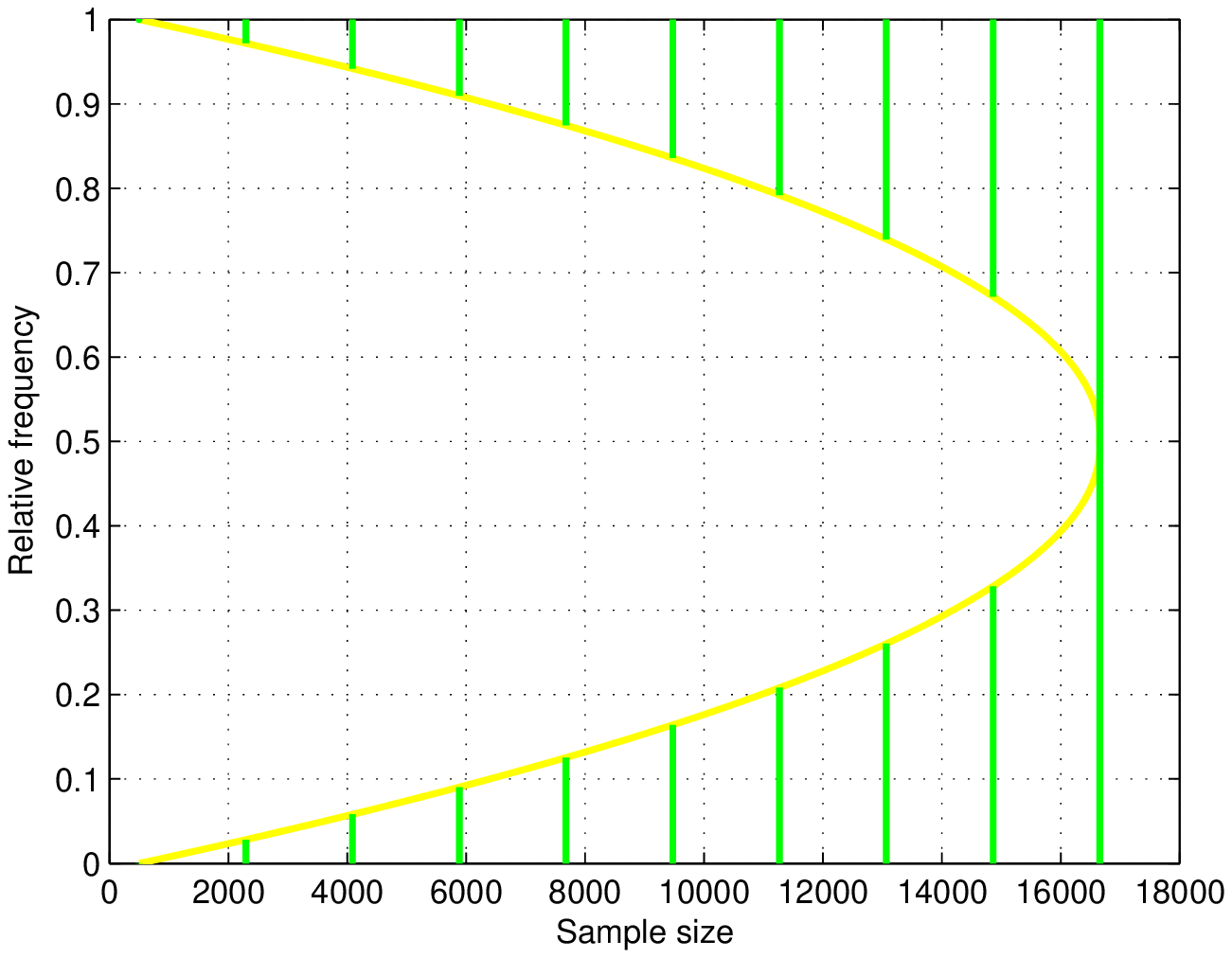} \includegraphics{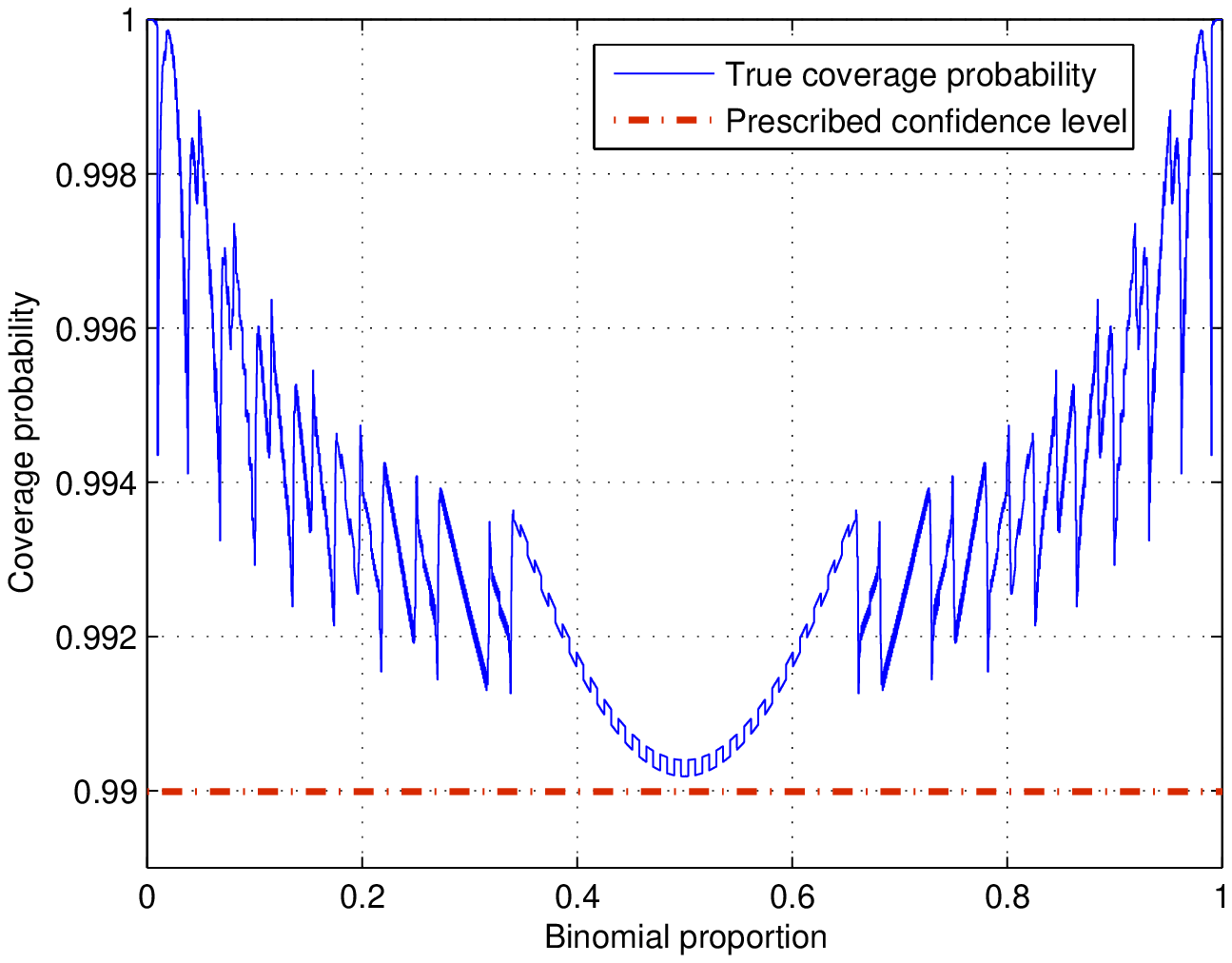} \caption{Double-parabolic sampling with $\vep = \de = 0.01, \; s = 10, \; \rho = \f{3}{4}$ and $\ze = 3.5753$}
\label{fig_ST_DP_0101_m10}
\end{figure}

To illustrate the use of these tables, suppose that one wants a ten-stage sampling procedure of approximately equal group sizes to ensure that
$\Pr \{ | \wh{\bs{p}} - p | < 0.01 \mid p \} > 0.99$ for any $p \in (0, 1)$.  This means that one can choose $\vep = \de = 0.01, \; s = 10$ and
sample sizes satisfying (\ref{equalss}).  To obtain appropriate parameter values for the sampling procedure, one can look at Table
\ref{table_abs_group01} to find the coverage tuning parameter $\ze$ corresponding to $\vep = 0.01$ and $s = 10$. From Table 3, it can be seen
that $\ze$ can be taken as $3.5753$. Consequently, the stopping rule is completely determined by substituting the values of design parameters
$\vep = 0.01, \; \de = 0.01, \; \ro = \f{3}{4}, \; \ze = 3.5753, \; s = 10$ into its definition and equation (\ref{equalss}).  The stopping
boundary of this sampling scheme and  the function of coverage probability with respect to the binomial proportion are displayed, respectively,
in the left and right sides of Figure \ref{fig_ST_DP_0101_m10}.

\subsection{Comparison of Sampling Schemes} \la{SecsubCSS}

\begin{figure}[here]
\includegraphics{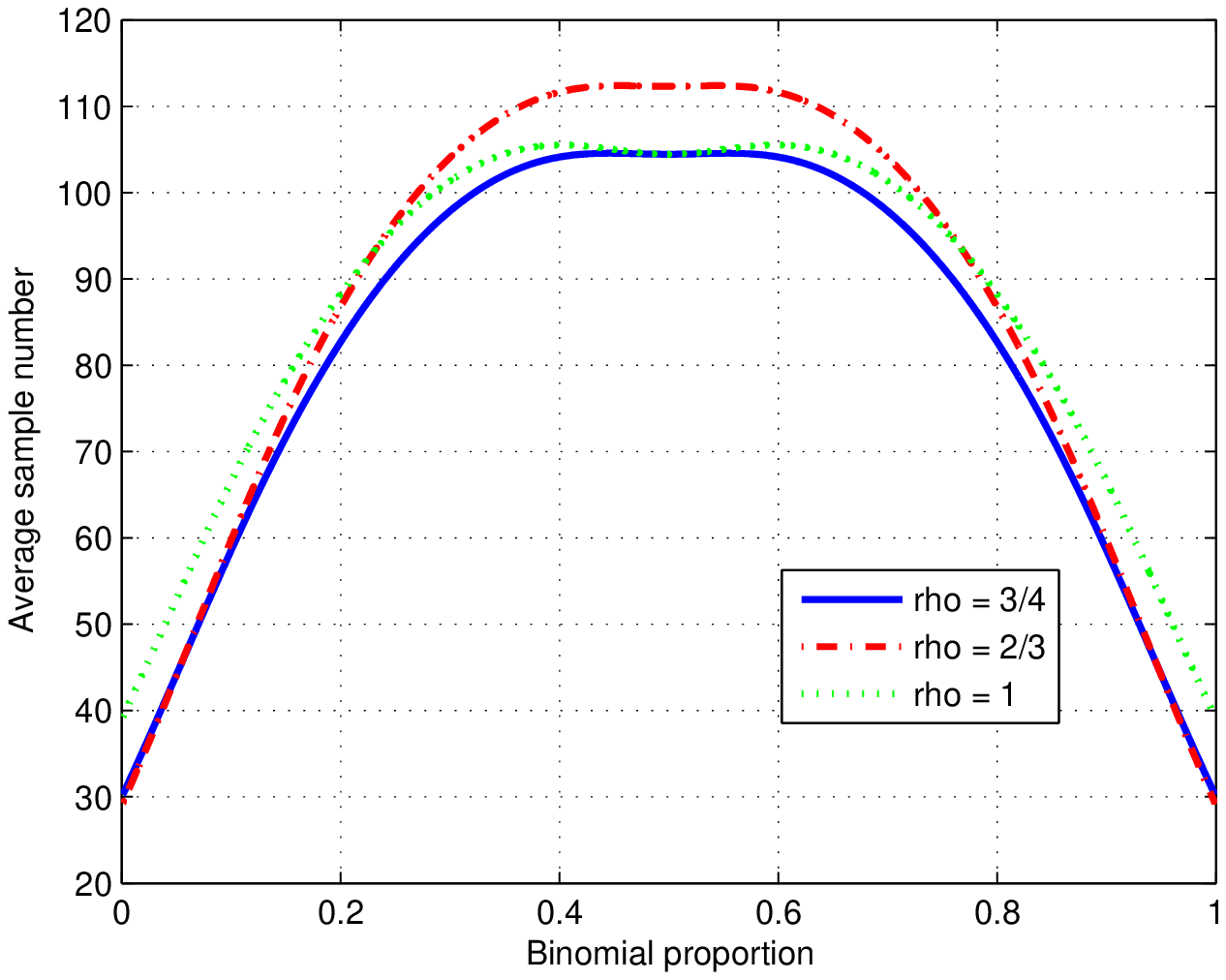} \includegraphics{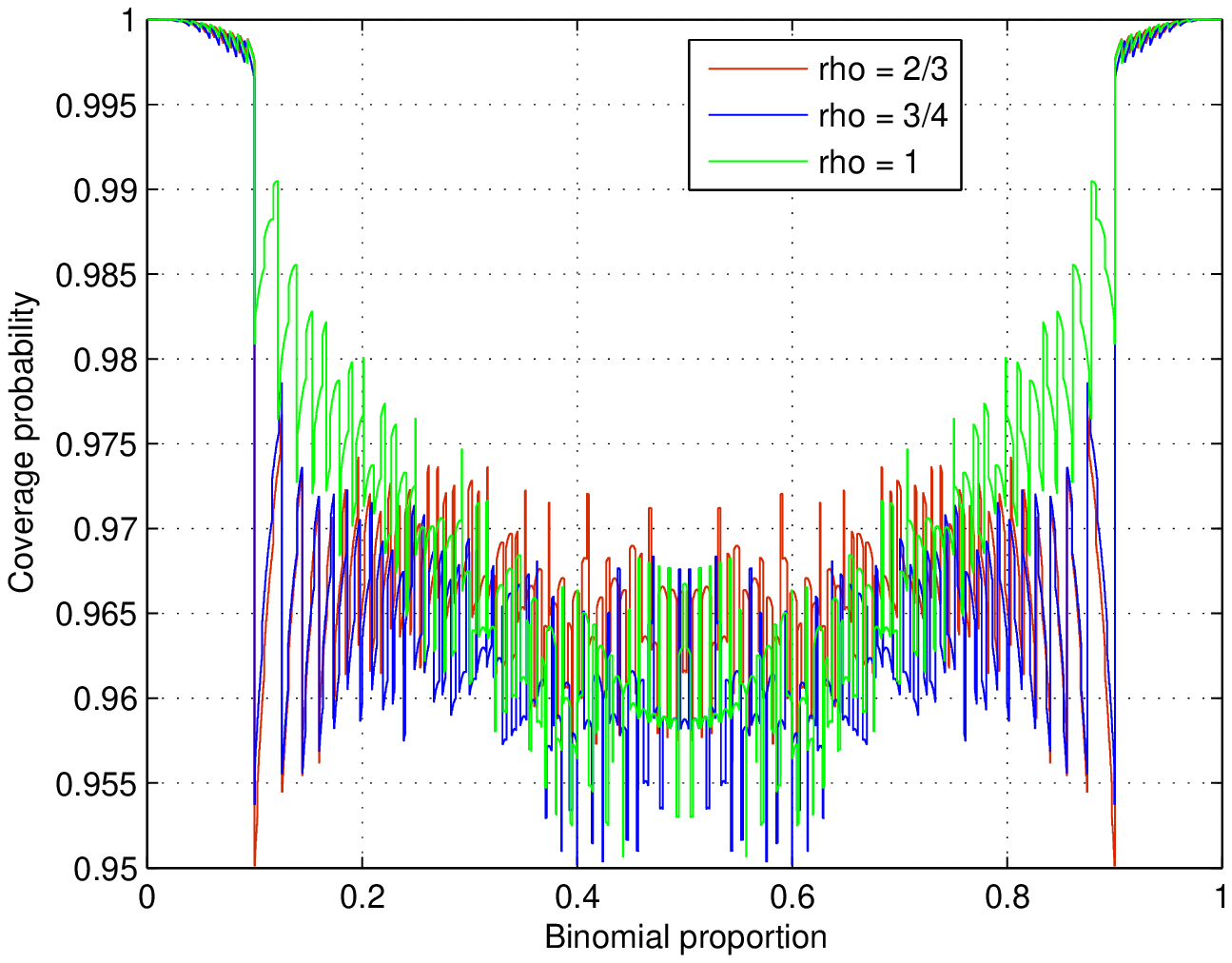} \caption{Double-parabolic sampling with various dilation coefficients} \label{fig_ASN_DP}
\end{figure}

We have conducted numerical experiments to investigate the impact of dilation coefficient $\ro$ on the performance of our double-parabolic
sampling schemes. Our computational experiences indicate that the dilation coefficient $\rho = \f{3}{4}$ is frequently  a good choice in terms
of average sample number and coverage probability.  For example, consider the case that the margin of error is given as $\vep = 0.1$ and the
prescribed confidence level is $1 - \de$ with $\de = 0.05$. For the double-parabolic sampling scheme with the dilation coefficient $\rho$ chosen
as $\f{2}{3}, \; \f{3}{4}$ and $1$, we have determined that,  to ensure the prescribed confidence level $1 - \de = 0.95$, it suffices to set the
coverage tuning parameter $\ze$ as $2.1, \; 2.4$ and $2.4$, respectively.  The average sample numbers of these sampling schemes and the coverage
probabilities as functions of the binomial proportion are shown, respectively, in the left and right sides of Figure \ref{fig_ASN_DP}.  From
Figure \ref{fig_ASN_DP}, it can be seen that a double-parabolic sampling scheme with dilation coefficient $\rho = \f{3}{4}$ has better
performance in terms of average sample number and coverage probability as compared to that of the double-parabolic sampling scheme with smaller
or larger values of dilation coefficient.

\bsk \bsk \bsk \bsk

\vspace*{2.5in}

\begin{figure}[here]
\includegraphics{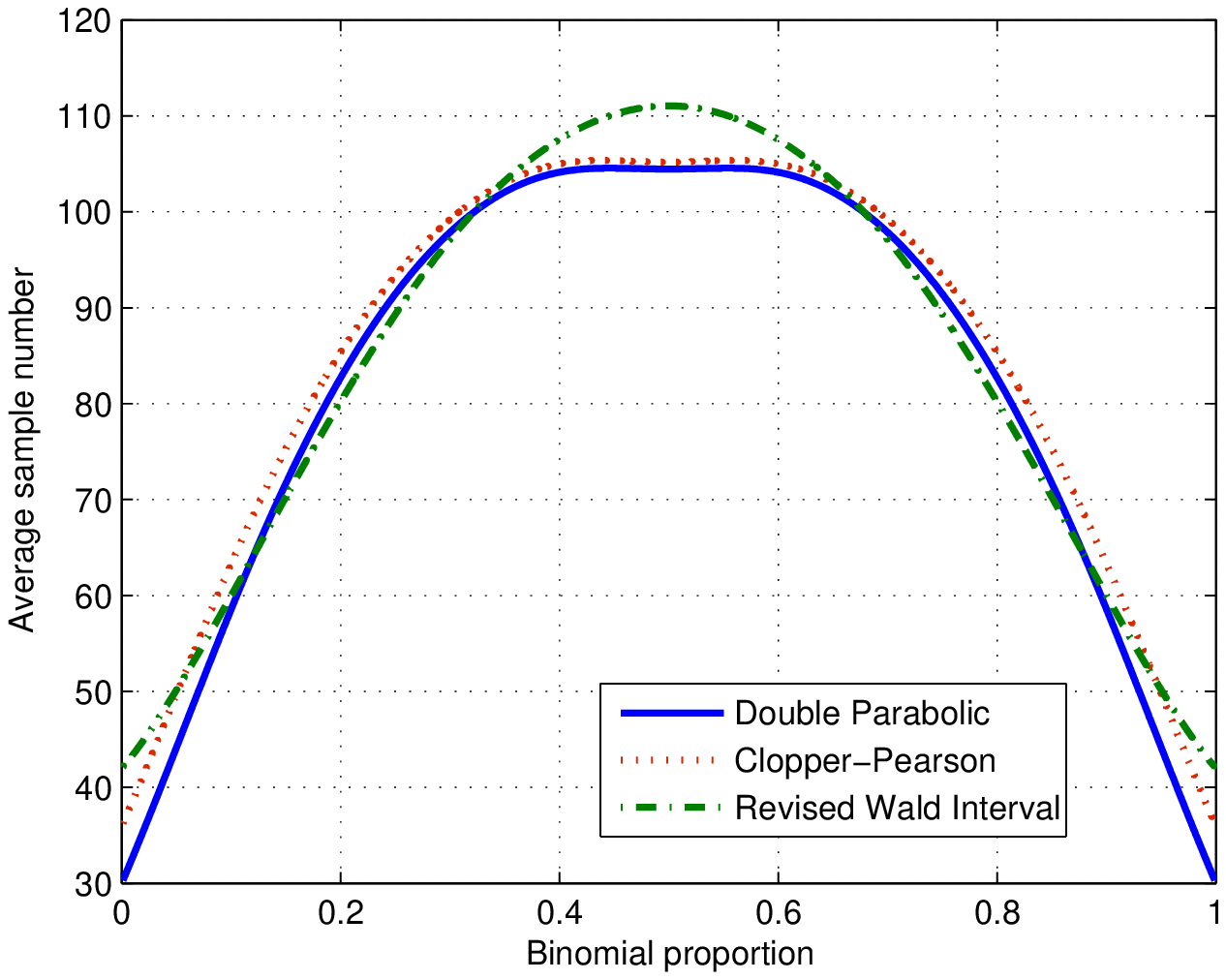} \includegraphics{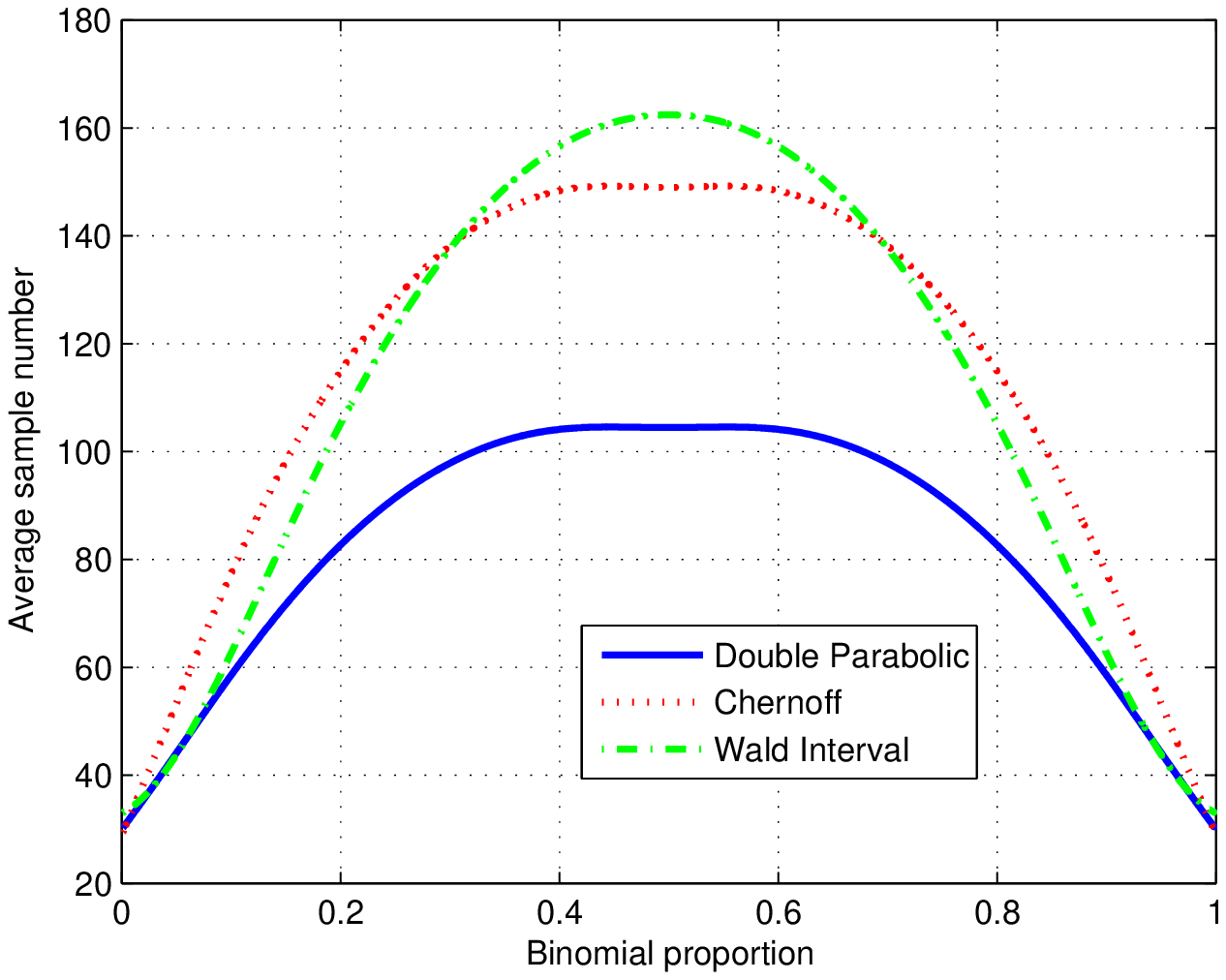} \caption{Comparison of average sample numbers} \label{fig_ASN_DP_CP_Frey_F5}
\end{figure}

We have investigated the impact of confidence intervals on the performance of fully sequential sampling schemes constructed from the inclusion
principle. We have observed that the stopping rule derived from Clopper-Pearson intervals generally outperforms the stopping rules derived from
other types of confidence intervals.  However, via appropriate choice of the dilation coefficient, the double-parabolic sampling scheme can
perform uniformly better than the stopping rule derived from Clopper-Pearson intervals.  To illustrate, consider the case that $\vep = 0.1$ and
$\de = 0.05$. For stopping rules derived from Clopper-Pearson intervals, Fishman's intervals, Wilson's intervals, and revised Wald intervals
with $a = 4$, we have determined that to guarantee the prescribed confidence level $1 - \de = 0.95$, it suffices to set the coverage tuning
parameter $\ze$ as $0.5, \; 1, \; 2.4$ and $0.37$, respectively.  For the stopping rule derived from Wald intervals, we have determined $\ze =
0.77$ to ensure the confidence level, under the condition that the minimum sample size is taken as {\small $\li \lc \f{1}{\vep} \ln \f{1}{\ze
\de}  \ri \rc$}. Recall that for the double-parabolic sampling scheme with $\ro = \f{3}{4}$, we have obtained $\ze = 2.4$ for purpose of
guaranteeing  the confidence level. The average sample numbers of these sampling schemes are shown in Figure \ref{fig_ASN_DP_CP_Frey_F5}. From
these plots, it can be seen that as compared to the stopping rule derived from Clopper-Pearson intervals, the stopping rule derived from the
revised Wald intervals performs better in the region of $p$ close to $0$ or $1$, but performs worse in the region of $p$ in the middle of $(0,
1)$. The performance of stopping rules from Fishman's intervals (i.e., from Chernoff bound) and Wald intervals are obviously inferior as
compared to that of the stopping rule derived from Clopper-Pearson intervals.  It can be observed that the double-parabolic sampling scheme
uniformly outperforms the stopping rule derived from Clopper-Pearson intervals.

\bsk \bsk  \bsk \bsk

\vspace*{2.5in}

\subsection{Estimation with High Confidence Level} \la{SecsubEHCL}

\begin{figure}  [here]
\includegraphics{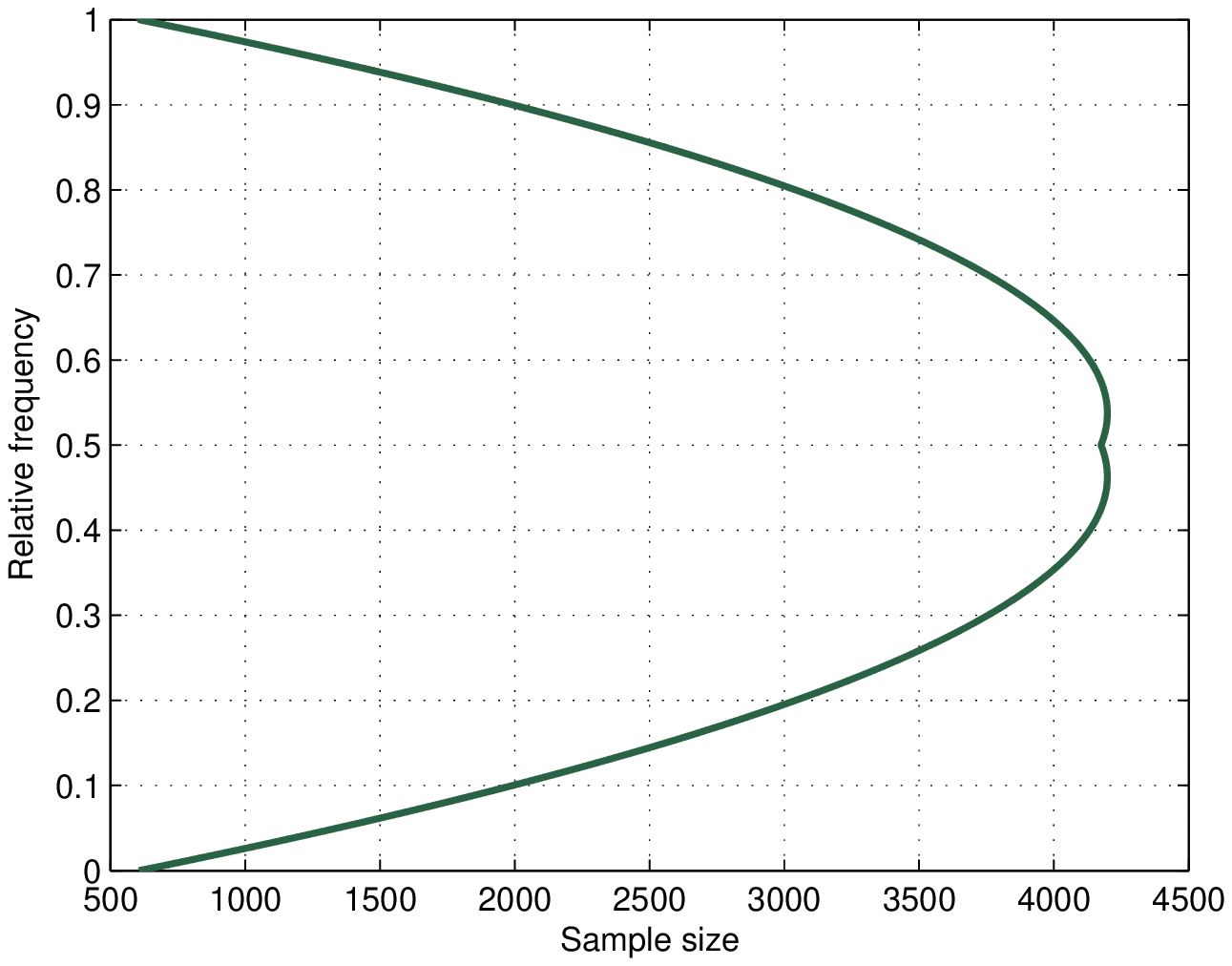} \includegraphics{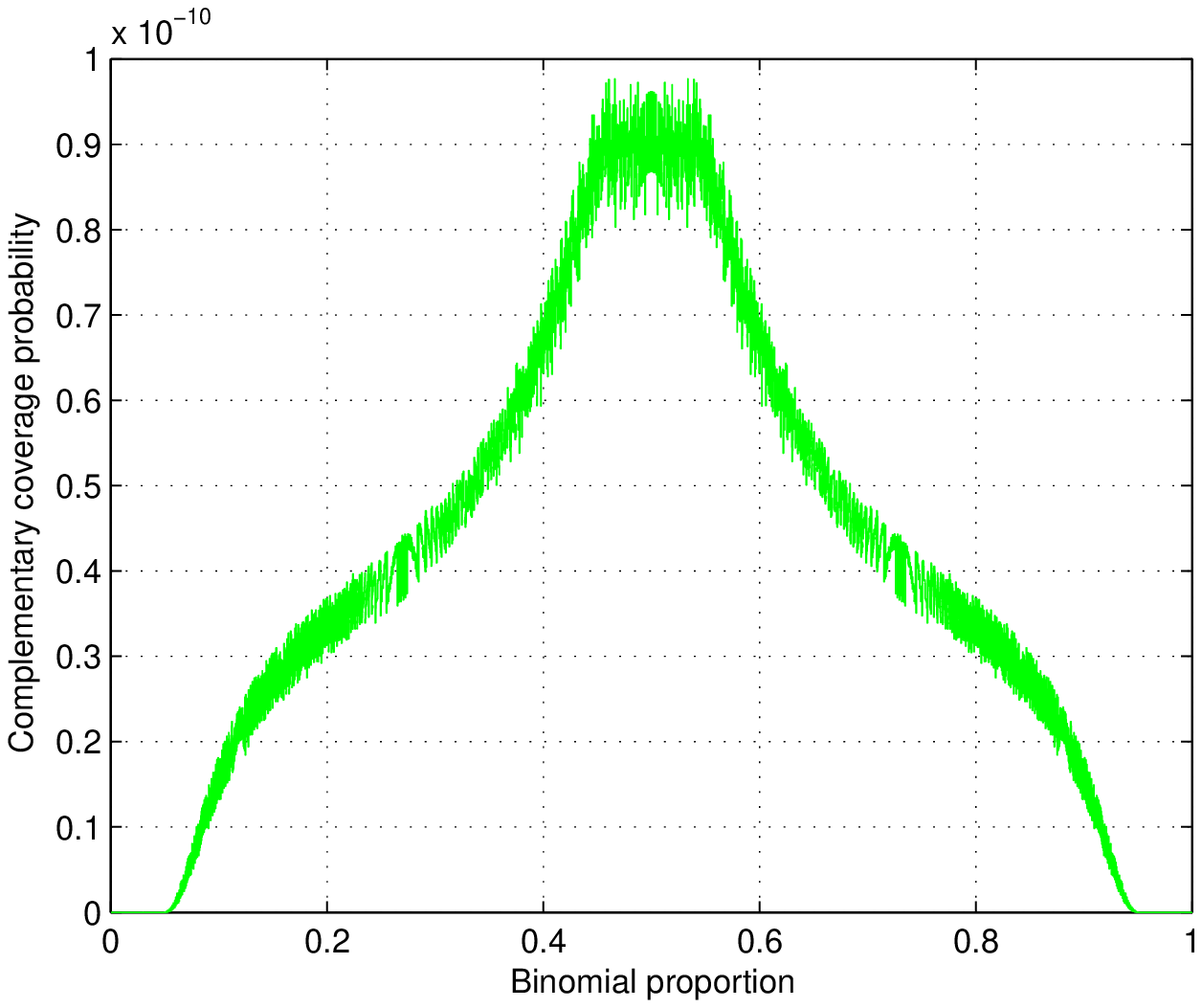} \caption{Double-parabolic sampling with $\vep = 0.05, \; \de = 10^{-10}, \; \ro = \f{3}{4}$ and $\ze = 7.65$}
\label{fig_ST_tiny_de}
\end{figure}

In some situations, we need to estimate a binomial proportion with a high confidence level.  For example, one might want to construct a sampling
scheme such that, for $\vep =  0.05$ and $\de = 10^{-10}$, the resultant sequential estimator $\wh{\bs{p}}$ satisfies $\Pr \{ | \wh{\bs{p}} - p
| < \vep \mid p \} > 1 - \de$ for any $p \in (0, 1)$. By working with the complementary coverage probability, we determined that it suffices to
let the dilation coefficient $\ro = \f{3}{4}$ and the coverage tuning parameter $\ze = 7.65$.   The stopping boundary and the function of
coverage probability with respect to the binomial proportion are displayed, respectively, in the left and right sides of Figure
\ref{fig_ST_tiny_de}. As addressed in Section \ref{SubsecWCCP},  it should be noted that it is impossible to obtain such a sampling scheme
without working with the complementary coverage probability.

\bsk \bsk  \bsk \bsk

 \vspace*{2.5in}

\section{Illustrative Examples for Clinical Trials} \la{SecsubIECT}

In this section, we shall illustrate the applications of our double-parabolic group sequential estimation method in clinical trials.

\begin{figure}[here]
\includegraphics{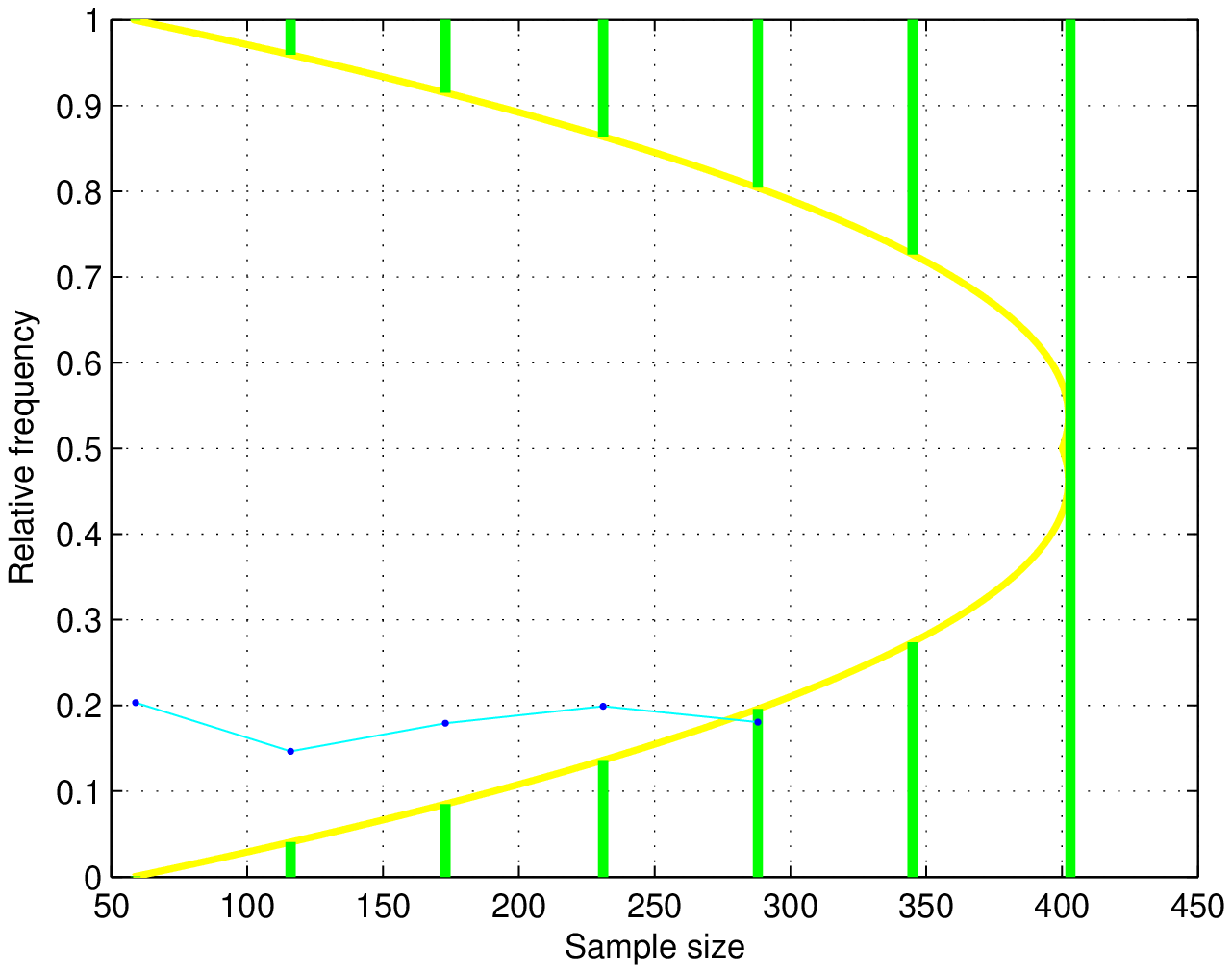} \includegraphics{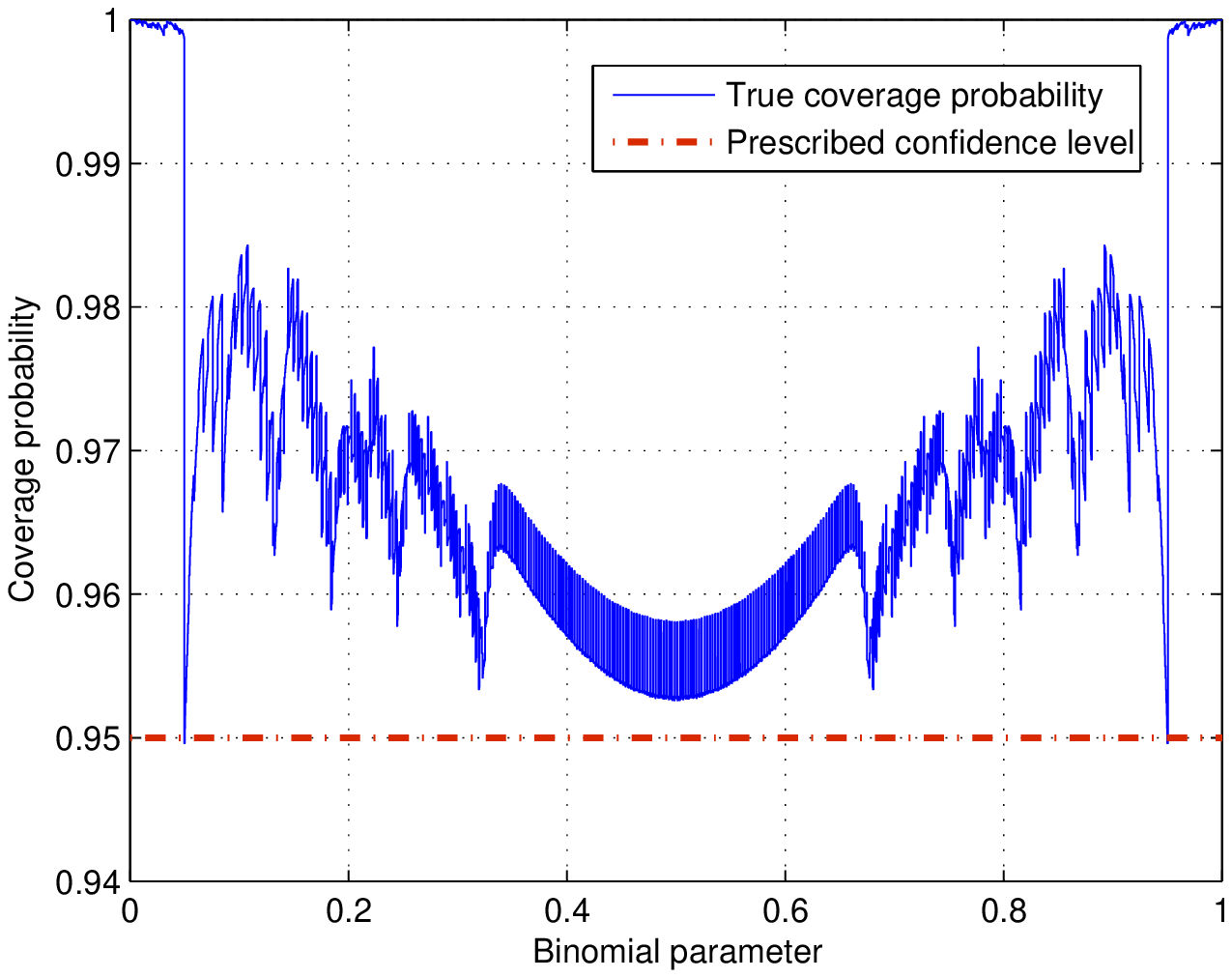} \caption{Double-parabolic sampling with $\vep = \de = 0.05, s = 7, \rho = \f{3}{4}$ and $\ze = 2.6759$}
\label{fig_ST_seven_0505_path}
\end{figure}

An example of our double-parabolic sampling scheme can be illustrated as follows.  Assume that $\vep = \de = 0.05$ is given and that the
sampling procedure is expected to have $7$ stages with sample sizes satisfying (\ref{equalss}).  Choosing $\rho = \f{3}{4}$, we have determined
that it suffices to take $\ze = 2.6759$ to guarantee that the coverage probability is no less than $1 - \de = 0.95$ for all $p \in (0, 1)$.
Accordingly, the sample sizes of this sampling scheme are calculated  as $59, 116, 173, 231, 288, 345$ and $403$. This sampling scheme, with a
sample path, is shown in the left side of Figure \ref{fig_ST_seven_0505_path}.   In this case, the stopping rule can be equivalently described
by virtue of Figure \ref{fig_ST_seven_0505_path} as: Continue sampling until $(\wh{\bs{p}}_\ell, n_\ell)$ hit a green line at some stage.    The
coverage probability is shown in the right side of Figure \ref{fig_ST_seven_0505_path}.

To apply this estimation method in a clinical trial for estimating the proportion $p$ of a binomial response with margin of error $0.05$ and
confidence level $95 \%$, we can have seven groups of patients with group sizes $59, 57, 57, 58, 57, 57$ and $58$.   In the first stage, we
conduct experiment with the $59$ patients of the first group.  We observe the relative frequency of response and record it as $\wh{\bs{p}}_1$.
Suppose there are $12$ patients having positive responses, then the relative frequency at the first stage is $\wh{\bs{p}}_1 = \f{12}{59} =
0.2034$. With the values of $(\wh{\bs{p}}_1, n_1) = (0.2034, 59)$, we check if the stopping rule is satisfied.  This is equivalent to see if the
point $(\wh{\bs{p}}_1, n_1)$ hit a green line at the first stage.  For such value of $(\wh{\bs{p}}_1, n_1)$, it can be seen that the stopping
condition is not fulfilled. So, we need to conduct the second stage of experiment with the $57$ patients of the second group. We observe the
response of these $57$ patients.  Suppose we observe that $5$ patients among this group have positive responses. Then, we add $5$ with $12$, the
number of positive responses before the second stage, to obtain $17$ positive responses among $n_2 = 59 + 57 = 116$ patients.  So, at the second
stage, we get the relative frequency  $\wh{\bs{p}}_2 = \f{17}{116} = 0.1466$. Since the stopping rule is not satisfied with the values of
$(\wh{\bs{p}}_2, n_2) = (0.1466, 116)$, we need to conduct the third stage of experiment with the $57$ patients of the third group. Suppose we
observe that $14$ patients among this group have positive responses. Then, we add $14$ with $17$, the number of positive responses before the
third stage, to get $31$ positive responses among $n_3 = 59 + 57 + 57 = 173$ patients. So, at the third stage, we get the relative frequency
$\wh{\bs{p}}_3 = \f{31}{173} = 0.1792$. Since the stopping rule is not satisfied with the values of $(\wh{\bs{p}}_3, n_3) = (0.1792, 173)$, we
need to conduct the fourth stage of experiment with the $58$ patients of the fourth group. Suppose we observe that $15$ patients among this
group have positive responses. Then, we add $15$ with $31$, the number of positive responses before the fourth stage, to get $46$ positive
responses among $n_4 = 59 + 57 + 57 + 58 = 231$ patients.  So, at the fourth stage, we get the relative frequency $\wh{\bs{p}}_4 = \f{46}{231} =
0.1991$. Since the stopping rule is not satisfied with the values of $(\wh{\bs{p}}_4, n_4) = (0.1991, 231)$, we need to conduct the fifth stage
of experiment with the $57$ patients of the fifth group. Suppose we observe that $6$ patients among this group have positive responses. Then, we
add $6$ with $46$, the number of positive responses before the fifth stage, to get $52$ positive responses among $n_5 = 59 + 57 + 57 + 58 + 57=
288$ patients.  So, at the fifth stage, we get the relative frequency $\wh{\bs{p}}_5 = \f{52}{288} = 0.1806$.  It can be seen that the stopping
rule is satisfied with the values of $(\wh{\bs{p}}_5, n_5) = (0.1806, 288)$. Therefore, we can terminate the sampling experiment and take
$\wh{\bs{p}} = \f{52}{288} = 0.1806$ as an estimate of the proportion of the whole population having positive responses.  With a $95 \%$
confidence level, one can believe that the difference between the true value of $p$ and its estimate $\wh{\bs{p}} = 0.1806$ is less than $0.05$.

\bsk  \bsk \bsk  \bsk \bsk

\vspace*{2.5in}

In this experiment, we only use $288$ samples to obtain the estimate for $p$. Except the roundoff error, there is no other source of error for
reporting statistical accuracy, since no asymptotic approximation is involved. As compared to fixed-sample-size procedure, we achieved a
substantial save of samples. To see this, one can check that using the rigorous formula (\ref{SSNCH}) gives a sample size $738$, which is overly
conservative. From the classical approximate formula (\ref{ASNAP}), the sample size is determined as $385$, which has been known to be
insufficient to guarantee the prescribed confidence level $95 \%$.  The exact method of \cite{ChenJP} shows that at least $391$ samples are
needed.  As compared to the best fixed sample size obtained by the method of \cite{ChenJP}, the reduction of sample sizes resulted from our
double-parabolic sampling scheme is $391 - 288 = 103$. It can be seen that the fixed-sample-size procedure wastes $\f{103}{288} = 35.76 \%$
samples as compared to our group sequential method, which is also an exact method.  This percentage may not be serious if it were a save of
number of simulation runs. However, as the number count is for patients, the reduction of samples is important for ethical and economical
reasons.   Using our group sequential method, the worst-case sample size is equal to $403$, which is only $12$ more than the minimum sample size
of fixed-sample procedure. However, a lot of samples can be saved in the average case.

\begin{figure}[here]
\includegraphics{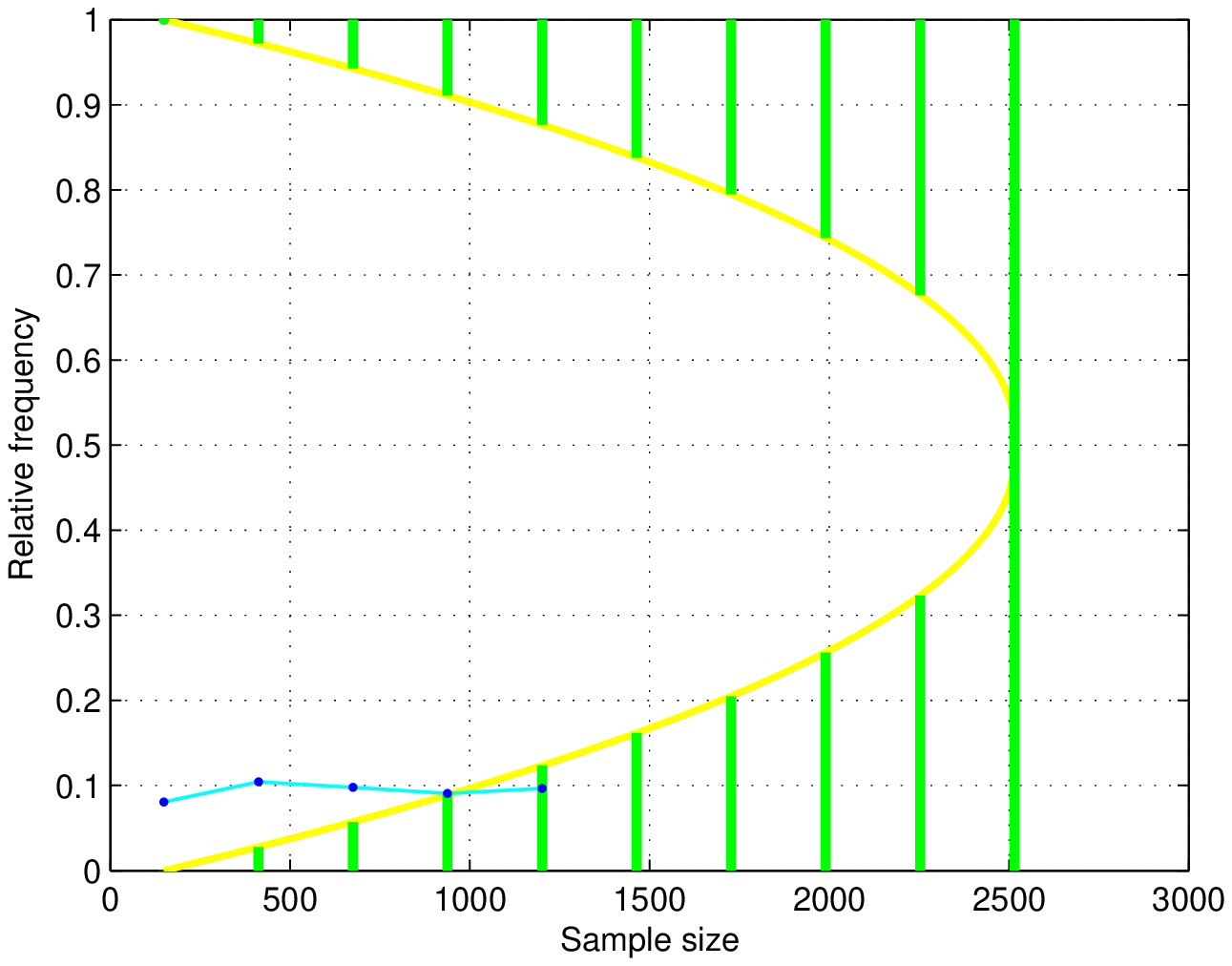} \includegraphics{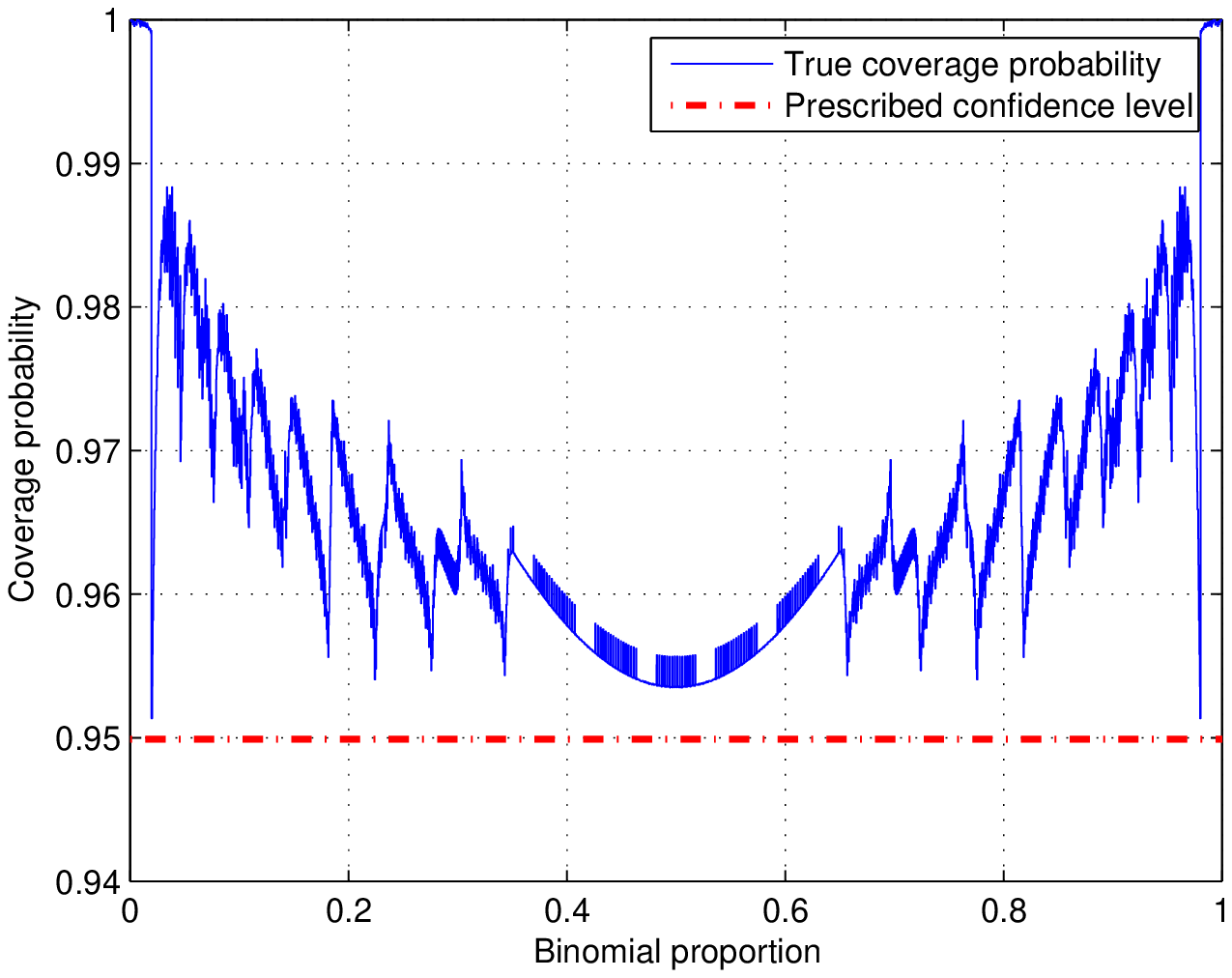} \caption{Double-parabolic sampling with $\vep = 0.02,  \; \de = 0.05, \; s = 10, \; \rho = \f{3}{4}$ and $\ze =
2.6725$} \label{fig_ST_DP_0205_m10}
\end{figure}

As $\vep$ or $\de$ become smaller, the reduction of samples is more significant.  For example, let $\vep = 0.02$ and $\de = 0.05$, we have a
double-parabolic sample scheme with $10$ stages.  The sampling scheme, with a sample path,  is shown in the left side of Figure
\ref{fig_ST_DP_0205_m10}. The coverage probability is shown in the right side of Figure \ref{fig_ST_DP_0205_m10}.

\bsk \bsk \bsk \bsk

\vspace*{2.5in}

\section{Conclusion}

In this paper, we have reviewed recent development of group sequential estimation methods for a binomial proportion.  We have illustrated the
inclusion principle and its applications to various stopping rules.  We have introduced computational techniques in the literature, which
suffice for determining parameters of stopping rules to guarantee desired confidence levels.  Moreover, we have proposed a new family of
sampling schemes with stopping boundary of double-parabolic shape, which are parameterized by the coverage tuning parameter and the dilation
coefficient. These parameters can be determined by the exact computational techniques to reduce the sampling cost, while ensuring prescribed
confidence levels. The new family of sampling schemes are extremely simple in structure and asymptotically optimal as the margin of error tends
to $0$.  We have established analytic bounds for the distribution and expectation of the sample number at the termination of the sampling
process.  We have obtained parameter values via the exact computational techniques for the proposed sampling schemes such that the confidence
levels are guaranteed and that the sampling schemes are generally more efficient as compared to existing ones.

\appendix

\sect{Proof of Theorem \ref{Wisleq} } \la{Wisleq_app}

Consider function $g (x, z) = \f{ (x - z)^2 }{x (1 - x)}$ for $x \in (0, 1)$ and $z \in [0, 1]$. It can be checked that  $\f{\pa g(x, z)}{\pa x}
= (x - z) [z (1 - x) + x (1 - z)] [x (1 - x)]^{-2}$,  which shows that for any fixed $z \in [0, 1]$,  $-g(x,z)$ is a unimodal function of $x \in
(0, 1)$, with a maximum attained at $x = z$. By such a property of $g(x,z)$ and the definition of Wilson's confidence intervals, we have {\small
\bee \li \{ \wh{\bs{p}}_\ell - \vep \leq L_\ell  \ri \} = \li \{ 0 < \wh{\bs{p}}_\ell - \vep \leq L_\ell  \ri \} \cup \{ \wh{\bs{p}}_\ell \leq
\vep \} = \li \{ 0 < \wh{\bs{p}}_\ell - \vep \leq L_\ell \leq \wh{\bs{p}}_\ell, \; g(L_\ell, \wh{\bs{p}}_\ell) = \f{ \mcal{Z}_{\ze \de}^2
}{n_\ell}  \ri \} \cup \{
\wh{\bs{p}}_\ell \leq \vep \} &  &\\
 = \li \{ \wh{\bs{p}}_\ell > \vep, \; \f{ \vep^2 } { (\wh{\bs{p}}_\ell - \vep) [1 - (\wh{\bs{p}}_\ell - \vep) ] } \geq  \f{ \mcal{Z}_{\ze
\de}^2 }{n_\ell} \ri \} \cup \{ \wh{\bs{p}}_\ell \leq  \vep \} \qqu \qqu \qqu \qqu \qqu \qqu \qqu \qu &  &
 \eee}
and {\small \bee &  & \li \{ \wh{\bs{p}}_\ell + \vep \geq U_\ell  \ri \} = \li \{ 1 > \wh{\bs{p}}_\ell + \vep \geq U_\ell \ri \} \cup \{
\wh{\bs{p}}_\ell + \vep \geq 1 \} = \li \{ 1 > \wh{\bs{p}}_\ell + \vep \geq U_\ell \geq \wh{\bs{p}}_\ell, \; g(U_\ell, \wh{\bs{p}}_\ell)
 = \f{ \mcal{Z}_{\ze \de}^2 }{n_\ell}  \ri \}\\
 &  &  \cup \{ \wh{\bs{p}}_\ell + \vep \geq 1 \}
  = \li \{ \wh{\bs{p}}_\ell < 1 - \vep, \; \f{ \vep^2 } { (\wh{\bs{p}}_\ell + \vep) [1 - (\wh{\bs{p}}_\ell + \vep) ] } \geq \f{ \mcal{Z}_{\ze
\de}^2 }{n_\ell} \ri \} \cup \{ \wh{\bs{p}}_\ell \geq  1 - \vep \}
 \eee}
for $\ell = 1, \cd, s$, where we have used the fact that $\{  \wh{\bs{p}}_\ell > \vep \} \subseteq \{ L_\ell > 0 \}, \;  \{  \wh{\bs{p}}_\ell <
1 - \vep \} \subseteq \{ U_\ell < 1 \}$ and $0 \leq L_\ell \leq \wh{\bs{p}}_\ell \leq U_\ell \leq 1$.  Recall that $0 < \vep < \f{1}{2}$.
It follows that {\small \bee &  & \li \{ \wh{\bs{p}}_\ell - \vep \leq L_\ell \leq U_\ell \leq  \wh{\bs{p}}_\ell + \vep \ri \}\\
&  & = \li \{ \vep < \wh{\bs{p}}_\ell < 1 - \vep, \; \f{ \vep^2 } { (\wh{\bs{p}}_\ell - \vep) [1 - (\wh{\bs{p}}_\ell - \vep) ] } \geq
\f{\mcal{Z}_{\ze \de}^2}{n_\ell}, \; \f{ \vep^2 } { (\wh{\bs{p}}_\ell + \vep) [1 - (\wh{\bs{p}}_\ell + \vep) ] } \geq
\f{ \mcal{Z}_{\ze \de}^2 }{n_\ell} \ri \} \\
&  & \bigcup \li \{ \wh{\bs{p}}_\ell \leq \vep, \; \f{ \vep^2 } { (\wh{\bs{p}}_\ell + \vep) [1 - (\wh{\bs{p}}_\ell + \vep) ] } \geq \f{
\mcal{Z}_{\ze \de}^2 }{n_\ell} \ri \}  \bigcup \li \{ \wh{\bs{p}}_\ell \geq 1 - \vep, \; \f{ \vep^2 } { (\wh{\bs{p}}_\ell - \vep)
[1 - (\wh{\bs{p}}_\ell - \vep) ] } \geq \f{ \mcal{Z}_{\ze \de}^2 }{n_\ell} \ri \} \\
&  & =  \li \{ \vep < \wh{\bs{p}}_\ell < 1 - \vep, \; \li ( \li | \wh{\bs{p}}_\ell - \f{1}{2} \ri | - \vep  \ri )^2 \geq \f{1}{4} - n_\ell \li (
\f{ \vep } {\mcal{Z}_{\ze \de} } \ri )^2  \ri \} \bigcup \\
&  & \li \{ \wh{\bs{p}}_\ell \leq \vep, \; \li ( \li | \wh{\bs{p}}_\ell - \f{1}{2} \ri | - \vep  \ri )^2 \geq \f{1}{4} - n_\ell \li ( \f{ \vep }
{\mcal{Z}_{\ze \de} } \ri )^2  \ri \} \bigcup \li \{ \wh{\bs{p}}_\ell \geq 1 - \vep, \; \li ( \li | \wh{\bs{p}}_\ell - \f{1}{2} \ri | - \vep \ri
)^2 \geq \f{1}{4} - n_\ell \li
( \f{ \vep } {\mcal{Z}_{\ze \de} } \ri )^2  \ri \}\\
&  & =  \li \{  \li ( \li | \wh{\bs{p}}_\ell - \f{1}{2} \ri | - \vep  \ri )^2 \geq \f{1}{4} - n_\ell \li ( \f{ \vep } {\mcal{Z}_{\ze \de} } \ri
)^2  \ri \}  \eee} for $\ell = 1, \cd, s$.  This completes the proof of the theorem.

\sect{Proof of Theorem \ref{UniformControl}} \la{UniformControl_app}

By the assumption that $n_s \geq \f{ 1 }{ 2 \vep^2} \ln \f{1}{\ze \de}$, we have  $\f{1}{4} + \f{ \vep^2 n_s }{2 \ln (\ze \de)} \leq 0$ and
consequently, {\small $\Pr  \{  (  | \wh{\bs{p}}_s - \f{1}{2}  | - \ro  \vep )^2 \geq \f{1}{4} + \f{ \vep^2 n_s }{2 \ln (\ze \de)} \} = 1$}. It
follows from the definition of the sampling scheme that the sampling process must stop at or before the $s$-th stage. In other words, $\Pr \{
\bs{l} \leq s \} = 1$.  This allows one to write \bel \Pr \{ | \wh{\bs{p}} - p | \geq \vep \mid p \} & = & \sum_{\ell = 1}^s \Pr \{  |
\wh{\bs{p}} - p | \geq \vep, \; \bs{l} = \ell \mid p \}
= \sum_{\ell = 1}^s \Pr \{  | \wh{\bs{p}}_\ell - p | \geq \vep, \; \bs{l} = \ell \mid p \} \nonumber\\
&  \leq & \sum_{\ell = 1}^s \Pr \{  | \wh{\bs{p}}_\ell - p | \geq \vep \mid p \} \la{crita} \eel for $p \in (0, 1)$.  By virtue of the
well-known Chernoff-Hoeffding bound \cite{Chernoff, Hoeffding}, we have \be \la{critb} \Pr \{  | \wh{\bs{p}}_\ell - p | \geq \vep \mid p \} \leq
2 \exp ( - 2 n_\ell \vep^2 ) \ee for $\ell = 1, \cd, s$. Making use of (\ref{crita}), (\ref{critb}) and the fact that $n_1 \geq 2 \ro (
\f{1}{\vep} - \ro) \ln \f{1}{\ze \de}$ as can be seen from (\ref{add338}), we have \bee \Pr \{ | \wh{\bs{p}} - p | \geq \vep \mid p \} & \leq &
2 \sum_{\ell = 1}^s \exp ( - 2 n_\ell \vep^2 )
\leq 2 \sum_{m = n_1}^\iy \exp ( - 2 m \vep^2 )   =   \f{ 2 \exp ( - 2 n_1 \vep^2 )  } { 1 - \exp ( - 2  \vep^2 ) }\\
&  \leq  & \f{ 2 \exp \li ( - 2 \vep^2 \times 2 \ro ( \f{1}{\vep} - \ro ) \ln \f{1}{\ze \de} \ri ) } { 1 - \exp ( - 2  \vep^2 ) } = \f{2 \exp
\li ( 4 \vep \ro ( 1 - \ro \vep) \ln (\ze \de)  \ri )  } { 1 - \exp ( - 2  \vep^2 ) }  \eee for any $p \in (0, 1)$.  Therefore, to guarantee
that $\Pr \{ | \wh{\bs{p}} - p | < \vep \mid p \} \geq 1 - \de$ for any $p \in (0, 1)$, it is sufficient to choose $\ze$ such that $2 \exp \li (
4 \vep \ro ( 1 - \ro \vep) \ln (\ze \de) \ri )   \leq \de[ 1 - \exp ( - 2  \vep^2 ) ]$.  This inequality can be written as $4 \vep \ro ( 1 - \ro
\vep) \ln (\ze \de) \leq \ln \f{\de}{2} + \ln \li [  1 - \exp ( - 2  \vep^2 ) \ri ]$ or equivalently, {\small $\ze \leq \f{1}{\de} \exp \li (
\f{ \ln \f{\de}{2} + \ln \li [  1 - \exp ( - 2  \vep^2 ) \ri ]}{ 4 \vep \ro ( 1 - \ro \vep) } \ri )$}.  The proof of the theorem is thus
completed.

\sect{Proof of Theorem \ref{AspOptimality}} \la{AspOptimality_app}

First, we need to show that $\Pr \{ \lim_{\vep \to 0}  \f{ \mbf{n} } { N ( p, \vep, \de, \ze ) } = 1 \mid p \} = 1$ for any $p \in (0, 1)$.
Clearly, the sample number $\mbf{n}$ is a random number dependent on $\vep$. Note that for any $\om \in \Om$, the sequences $\{ \ovl{X}_{
\mbf{n} (\om) } (\om) \}_{\vep \in (0, 1) }$ and $\{ \ovl{X}_{ \mbf{n} (\om) - 1} (\om) \}_{\vep \in (0, 1) }$ are subsets of $\{ \ovl{X}_m
(\om)  \}_{m = 1}^\iy$. By the strong law of large numbers, for almost every $\om \in \Om$, the sequence $\{ \ovl{X}_m (\om)  \}_{m = 1}^\iy$
converges to $p$. Since every subsequence of a convergent sequence must converge, it follows that the sequences $\{ \ovl{X}_{ \mbf{n} (\om) }
(\om) \}_{\vep \in (0, 1) }$ and $\{ \ovl{X}_{ \mbf{n} (\om) - 1 } (\om) \}_{\vep \in (0, 1) }$ converge to $p$ as $\vep \to 0$ provided that
$\mbf{n} (\om) \to \iy$ as $\vep \to 0$. Since it is certain that $\mbf{n} \geq 2 \ro ( \f{1}{\vep}  - \ro ) \ln \f{1}{\ze \de}  \to \iy$ as
$\vep \to 0$, we have that $\li \{ \lim_{\vep \to 0}  \f{ \mbf{n} - 1 }{ \mbf{n} } = 1 \ri \}$ is a sure event.  It follows that $B = \{
\lim_{\vep \to 0} \ovl{X}_{\mbf{n} - 1} = p, \; \lim_{\vep \to 0}  \ovl{X}_{\mbf{n} } = p, \; \lim_{\vep \to 0} \f{ \mbf{n} - 1 }{ \mbf{n} } = 1
\}$  is an almost sure event.   By the definition of the sampling scheme, we have that
\[
A = \li \{ \li (  \li | \ovl{X}_{\mbf{n} - 1}  - \f{1}{2} \ri | - \ro  \vep \ri )^2 < \f{1}{4} + \f{ \vep^2 (\mbf{n} - 1) }{2 \ln (\ze \de)},
\qu \li ( \li | \ovl{X}_{\mbf{n}} - \f{1}{2} \ri | - \ro  \vep \ri )^2 \geq \f{1}{4} + \f{ \vep^2 \mbf{n} }{2 \ln (\ze \de)} \ri \}
\]
is a sure event.  Hence, $A \cap B$ is an almost sure event.  Define $C = \li \{  \lim_{\vep \to 0} \f{ \mbf{n}  }{ N ( p, \vep, \de, \ze )  } =
1 \ri \}$. We need to show that $C$ is an almost sure event.  For this purpose,  we let $\om \in A \cap B$ and expect to show that $\om \in C$.
As a consequence of $\om \in A \cap B$, {\small \[ \f{ \mbf{n} (\om)   }{  N ( p, \vep, \de, \ze )  } < \f{ \mbf{n} (\om) } { \mbf{n} (\om)  - 1
} \f{ \li [ \f{1}{4} -  \li ( \li | \ovl{X}_{\mbf{n} (\om)  - 1} (\om) - \f{1}{2} \ri | - \ro  \vep \ri )^2 \ri ]}{p ( 1 - p)}, \qu \lim_{\vep
\to 0} \ovl{X}_{\mbf{n} (\om) - 1} (\om)  = p, \qu \lim_{\vep \to 0} \f{ \mbf{n} (\om) - 1 }{ \mbf{n} (\om)  } = 1.
\]}
By the continuity of the function $\li | x - \f{1}{2} \ri | - \ro \vep$ with respect to $x$ and $\vep$, we have \be \la{ineqsup} \limsup_{\vep
\to 0} \f{ \mbf{n} (\om) }{ N ( p, \vep, \de, \ze ) } \leq  \lim_{\vep \to 0} \f{ \mbf{n} (\om) } { \mbf{n} (\om) - 1 } \times \f{ \li [
\f{1}{4} - \li ( \li | \lim_{\vep \to 0} \ovl{X}_{\mbf{n} (\om)  - 1} (\om) - \f{1}{2} \ri | - \lim_{\vep \to 0} \ro  \vep \ri )^2 \ri ]}{p ( 1
- p)} = 1. \qqu \ee On the other hand, as a consequence of $\om \in A \cap B$,
\[
\f{ \mbf{n} (\om) }{ N ( p, \vep, \de, \ze )  } \geq \f{ \li [ \f{1}{4} -  \li (  \li | \ovl{X}_{\mbf{n} (\om) } (\om) - \f{1}{2} \ri | - \ro
\vep \ri )^2 \ri ]}{p ( 1 - p)}, \qqu \qqu \lim_{\vep \to 0}  \ovl{X}_{\mbf{n} (\om) } (\om) = p.
\]
Making use of the continuity of the function $\li | x - \f{1}{2} \ri | - \ro \vep$ with respect to $x$ and $\vep$, we have \be \la{ineqinf}
\liminf_{\vep \to 0} \f{ \mbf{n}  (\om) }{ N ( p, \vep, \de, \ze )  }  \geq \f{ \li [ \f{1}{4} -  \li ( \li | \lim_{\vep \to 0} \ovl{X}_{\mbf{n}
(\om) } (\om) - \f{1}{2} \ri | - \lim_{\vep \to 0} \ro  \vep \ri )^2 \ri ]}{p ( 1 - p)} = 1. \ee  Combining (\ref{ineqsup}) and (\ref{ineqinf})
yields $\lim_{\vep \to 0} \f{ \mbf{n} (\om) }{ N ( p, \vep, \de, \ze ) } = 1$ and thus $A \cap B \subseteq C$.  This implies that $C$ is an
almost sure event and thus  {\small $\Pr \li \{ \lim_{\vep \to 0} \f{ \mbf{n} }{ N ( p, \vep, \de, \ze ) }  = 1  \mid p \ri \} = 1$} for $p \in
(0, 1)$.

Next, we need to show that $\lim_{\vep \to 0} \Pr \{ | \wh{\bs{p}} - p | < \vep \mid p \} = 2 \Phi \li ( \sq{ 2 \ln \f{1}{\ze \de} } \ri )  - 1$
for any $p \in (0, 1)$.  For simplicity of notations, let $\si = \sq{ p ( 1 - p) }$ and $a = \sq{ 2 \ln \f{1}{\ze \de} }$.   Note that $\Pr \{ |
\wh{\bs{p}} - p | < \vep \mid p \} = \Pr \{ | \ovl{X}_\mbf{n} - p | < \vep \mid p \} = \Pr \{ \sq{\mbf{n}}| \ovl{X}_\mbf{n} - p | \sh \si < \vep
\sq{\mbf{n}} \sh \si \}$. Clearly, for any $\eta \in (0, a)$, {\small \bel \Pr \{ \sq{\mbf{n}}| \ovl{X}_\mbf{n} -  p | \sh \si < \vep
\sq{\mbf{n}} \sh \si \} & \leq & \Pr \{ \sq{\mbf{n}}| \ovl{X}_\mbf{n} -  p | \sh \si < \vep \sq{\mbf{n}} \sh \si, \;  \vep
\sq{\mbf{n}} \sh \si \in [a - \eta, a + \eta] \} \nonumber\\
&  &  + \Pr \{  \vep \sq{\mbf{n}} \sh \si \notin [a - \eta, a + \eta] \} \nonumber\\
&  \leq  & \Pr \{ \sq{\mbf{n}}| \ovl{X}_\mbf{n} -  p | \sh \si < a + \eta, \;  \vep
\sq{\mbf{n}} \sh \si \in [a - \eta, a + \eta] \} \nonumber\\
&  &  + \Pr \{  \vep \sq{\mbf{n}} \sh \si \notin [a - \eta, a + \eta] \} \nonumber\\
&  \leq  & \Pr \{ \sq{\mbf{n}}| \ovl{X}_\mbf{n} -  p | \sh \si < a + \eta \}
 + \Pr \{  \vep \sq{\mbf{n}} \sh \si \notin [a - \eta, a +
\eta] \} \qqu \qqu \la{ineqA} \eel}  and \bel &  & \Pr \{ \sq{\mbf{n}}| \ovl{X}_\mbf{n} -  p | \sh \si < \vep \sq{\mbf{n}} \sh \si \} \geq  \Pr
\{ \sq{\mbf{n}}| \ovl{X}_\mbf{n} -  p | \sh \si
< \vep \sq{\mbf{n}} \sh \si, \;  \vep \sq{\mbf{n}} \sh \si \in [a - \eta, a + \eta] \} \nonumber\\
&   & \geq  \Pr \{ \sq{\mbf{n}}| \ovl{X}_\mbf{n} -  p | \sh \si <  a - \eta, \;  \vep \sq{\mbf{n}} \sh \si \in [a - \eta, a + \eta] \} \nonumber\\
&  & \geq \Pr \{ \sq{\mbf{n}}| \ovl{X}_\mbf{n} -  p | \sh \si < a - \eta \}
 - \Pr \{  \vep \sq{\mbf{n}} \sh \si \notin [a - \eta, a +
\eta] \}. \la{ineqB}  \eel Recall that we have established that $\mbf{n} \sh N ( p, \vep, \de, \ze ) \to 1$  almost surely as $\vep \to 0$. This
implies that $\vep \sq{\mbf{n}} \sh \si \to a$ and $\mbf{n} \sh N ( p, \vep, \de, \ze ) \to 1$ in probability as $\vep$ tends to zero.  It
follows from Anscombe's random central limit theorem \cite{Ascombe} that as $\vep$ tends to zero, $\sq{\mbf{n}} ( \ovl{X}_\mbf{n} - p ) \sh \si$
converges in distribution to a Gaussian random variable with zero mean and unit variance.  Hence, from (\ref{ineqA}), \bee &  & \limsup_{\vep
\to 0}  \Pr \{ \sq{\mbf{n}}|
\ovl{X}_\mbf{n} -  p | \sh \si < \vep \sq{\mbf{n}} \sh \si \}\\
&  & \leq \lim_{\vep \to 0}  \Pr \{ \sq{\mbf{n}}| \ovl{X}_\mbf{n} -  p | \sh \si < a + \eta \}
 + \lim_{\vep \to 0} \Pr \{  \vep \sq{\mbf{n}} \sh \si \notin [a - \eta, a +
\eta] \} = 2 \Phi (a + \eta) - 1 \eee  and from (\ref{ineqB}), \bee &  & \liminf_{\vep \to 0} \Pr \{ \sq{\mbf{n}}| \ovl{X}_\mbf{n} -  p | \sh \si < \vep \sq{\mbf{n}} \sh \si \}\\
&  & \geq \lim_{\vep \to 0} \Pr \{ \sq{\mbf{n}}| \ovl{X}_\mbf{n} -  p | \sh \si < a - \eta \}
 - \lim_{\vep \to 0} \Pr \{  \vep \sq{\mbf{n}} \sh \si \notin [a - \eta, a +
\eta] \} = 2 \Phi (a - \eta) - 1.  \eee  Since this argument holds for arbitrarily small $\eta \in (0, a)$,  it must be true that
\[
\liminf_{\vep \to 0} \Pr \{ \sq{\mbf{n}}| \ovl{X}_\mbf{n} -  p | \sh \si < \vep \sq{\mbf{n}} \sh \si  \} = \limsup_{\vep \to 0}  \Pr \{
\sq{\mbf{n}}| \ovl{X}_\mbf{n} -  p | \sh \si < \vep \sq{\mbf{n}} \sh \si \} = 2 \Phi (a) - 1.
\]
So, {\small $\lim_{\vep \to 0} \Pr \{ | \wh{\bs{p}} - p | < \vep \mid p \} = \lim_{\vep \to 0}  \Pr \{ \sq{\mbf{n}}| \ovl{X}_\mbf{n} -  p | \sh
\si < \vep \sq{\mbf{n}} \sh \si \} = 2 \Phi (a) - 1 = 2 \Phi \li ( \sq{ 2 \ln \f{1}{\ze \de} } \ri )  - 1$} for any $p \in (0, 1)$.

Now, we focus our attention to show that $\lim_{\vep \to 0} \f{ \bb{E} [ \mbf{n} ] } { N ( p, \vep, \de, \ze ) } = 1$ for any $p \in (0, 1)$.
For this purpose, it suffices to show that \be \la{suf}
 1 - \eta \leq \liminf_{\vep \to 0} \f{ \bb{E}
[ \mbf{n} ] } { N ( p, \vep, \de, \ze ) } \leq \limsup_{\vep \to 0} \f{ \bb{E} [ \mbf{n} ]  } { N ( p, \vep, \de, \ze ) } \leq 1 + \eta, \qqu
\fa p \in (0, 1) \ee for any $\eta \in (0, 1)$.  For simplicity of notations, we abbreviate $N ( p, \vep, \de, \ze )$ as $N$ in the sequel.
Since we have established $\Pr \{ \lim_{\vep \to 0} \f{ \mbf{n} } { N ( p, \vep, \de, \ze ) }  = 1 \} = 1$, we can conclude that \be
\lim_{\vep \to 0} \Pr \{ (1 - \eta) N  \leq  \mbf{n} \leq (1 + \eta) N \} = 1.  \la{limeq6} \\
\ee  Noting that {\small \bee \bb{E} [ \mbf{n} ] = \sum_{m = 0}^\iy m \Pr \{ \mbf{n} = m \}  \geq \sum_{ (1 - \eta) N \leq m \leq (1 + \eta) N }
m \Pr \{ \mbf{n} = m \} \geq (1 - \eta) N  \sum_{ (1 - \eta) N \leq m \leq (1 + \eta) N   } \Pr \{ \mbf{n} = m \},  \eee}
we have \be \la{ieq99}
\bb{E} [ \mbf{n} ]  \geq (1 - \eta) N \Pr \{ (1 - \eta) N \leq \mbf{n} \leq (1 + \eta) N \}. \ee Combining (\ref{limeq6}) and (\ref{ieq99})
yields
\[
\liminf_{\vep \to 0} \f{ \bb{E} [ \mbf{n} ]  } { N ( p, \vep, \de, \ze ) } \geq (1 - \eta) \lim_{\vep \to 0} \Pr \{ (1 - \eta) N  \leq \mbf{n}
\leq (1 + \eta) N \} = 1 - \eta.
\]
On the other hand, using $\bb{E} [ \mbf{n} ] = \sum_{m = 0}^\iy \Pr \{ \mbf{n} > m  \}$, we can write {\small \bee \bb{E} [ \mbf{n} ]  = \sum_{0
\leq m < (1 + \eta) N } \Pr \{ \mbf{n}
> m \}  + \sum_{m \geq (1 + \eta) N } \Pr \{ \mbf{n} > m \} \leq  \lc (1 + \eta) N  \rc + \sum_{m \geq (1 + \eta) N } \Pr \{ \mbf{n} > m \}.
\eee}  Since {\small $\limsup_{\vep \to 0}  \f{ \lc (1 + \eta) N  \rc  } { N ( p, \vep, \de, \ze ) } =  1 + \eta$},  for the purpose of
establishing $\limsup_{\vep \to 0} \f{ \bb{E} [ \mbf{n} ] } { N ( p, \vep, \de, \ze ) } \leq 1 + \eta$, it remains to show that
\[
\limsup_{\vep \to 0} \f{ \sum_{m \geq (1 + \eta) N} \Pr \{ \mbf{n} > m \} } { N ( p, \vep, \de, \ze ) } = 0.
\]
Consider functions $f(x) = \f{1}{4} - \li ( \li | x - \f{1}{2} \ri | - \ro  \vep \ri
 )^2$ and $g(x) = x ( 1 - x)$ for $x \in [0, 1]$.  Note that
 \bee
\li |  f(x) - g (x) \ri |  = \li |  \li ( x - \f{1}{2} \ri )^2  -  \li ( \li | x - \f{1}{2} \ri | - \ro  \vep \ri )^2  \ri | = \ro \vep \li | |
2 x - 1 | - \ro  \vep \ri | \leq \ro \vep ( 1 + \ro \vep)
 \eee
for all $x \in [0, 1]$. For $p \in (0, 1)$, there exists a positive number $\ga < \min \{ p, 1 - p \}$ such that $| g(x) - g(p) | < \f{\eta}{2}
p ( 1 - p)$ for any $x \in (p - \ga, p + \ga)$, since $g(x)$ is a continuous function of $x$.  From now on, let $\vep > 0$ be sufficiently small
such that $\ro \vep ( 1 + \ro \vep) < \f{\eta}{2} p ( 1 - p)$.  Then, \[ f(x)  \leq g (x)  +  \ro \vep ( 1 + \ro \vep) < g(p) + \f{\eta}{2} p (
1 - p) +  \ro \vep ( 1 + \ro \vep) < (1 + \eta) p ( 1 - p) \] for all $x \in (p - \ga, p + \ga)$.  This implies that \be \la{inc98}
 \{ \ovl{X}_m \in (p - \ga, p + \ga) \} \subseteq \li \{ (1 + \eta) p ( 1 -
p) \geq \f{1}{4} - \li ( \li | \ovl{X}_m - \f{1}{2} \ri | - \ro  \vep \ri )^2 \ri \}
 \ee for all $m > 0$. Taking complementary events on both sides of (\ref{inc98}) leads to
\[
\li  \{ (1 + \eta) p ( 1 - p) < \f{1}{4} -  \li (  \li | \ovl{X}_m - \f{1}{2} \ri | - \ro  \vep \ri )^2 \ri \} \subseteq \{ \ovl{X}_m \notin (p
- \ga, p + \ga) \}
 \]  for all $m > 0$.  Since $(1 + \eta) p
 ( 1 - p) = \f{ (1 + \eta) N \vep^2 }{ 2 \ln \f{1}{\ze \de} }
 \leq \f{ m \vep^2 }{ 2 \ln \f{1}{\ze \de} }$
  for all $m \geq (1 + \eta) N$, it follows that \[
\li \{ \f{ m \vep^2 }{ 2 \ln \f{1}{\ze \de} }  < \f{1}{4} - \li ( \li | \ovl{X}_m - \f{1}{2} \ri | - \ro  \vep \ri )^2  \ri \} \subseteq \{
\ovl{X}_m \notin (p - \ga, p + \ga) \}
\]
for all $m \geq (1 + \eta) N$.  Therefore, we have shown that if $\vep$ is sufficiently small, then there exists a number $\ga > 0$ such that
\[
\{  \mbf{n} > m \} \subseteq  \li \{  \li ( \li | \ovl{X}_m - \f{1}{2} \ri | - \ro  \vep \ri )^2 < \f{1}{4} + \f{m \vep^2}{2 \ln (\ze \de)} \ri
\} \subseteq \{ \ovl{X}_m \notin (p - \ga, p + \ga) \}
\]
for all $m \geq (1 + \eta) N$.  Using this inclusion relationship and the Chernoff-Hoeffding bound \cite{Chernoff, Hoeffding}, we have \be
\la{CH98} \Pr \{  \mbf{n}
> m \} \leq \Pr \{ \ovl{X}_m \notin (p - \ga, p + \ga) \} \leq 2 \exp ( - 2 m \ga^2 ) \ee for all $m \geq (1 + \eta) N$ provided that $\vep > 0$
is sufficiently small.  Letting  $k = \lc (1 + \eta) N \rc$ and using (\ref{CH98}), we have \bee &  &  \sum_{m \geq (1 + \eta) N} \Pr \{ \mbf{n}
> m \} = \sum_{m \geq k }  \Pr \{  \mbf{n} > m \} \leq \sum_{m \geq k} 2 \exp ( - 2 m \ga^2 )  = \f{2  \exp ( - 2 k \ga^2 ) } {  1 - \exp ( - 2
\ga^2 ) } \eee  provided that $\vep$ is sufficiently small.     Consequently,
 \[
\limsup_{\vep \to 0} \f{ \sum_{m \geq (1 + \eta) N} \Pr \{ \mbf{n} > m \} } { N ( p, \vep, \de, \ze ) } \leq \limsup_{\vep \to 0} \f{2}{N}
\f{\exp ( - 2 k \ga^2 ) } {  1 - \exp ( - 2 \ga^2 ) }  = 0,
 \]
since $k \to \iy$ and $N \to \iy$ as $\vep \to 0$.  So, we have established (\ref{suf}).  Since the argument holds for arbitrarily small $\eta
> 0$, it must be true that $\lim_{\vep \to 0} \f{ \bb{E} [ \mbf{n} ] } { N ( p, \vep, \de, \ze ) } = 1$ for any $p \in (0, 1)$.  This completes the proof of the theorem.

\sect{Proof of Theorem \ref{bounds} } \la{App_bounds}

Recall that $\bs{l}$ denotes the index of stage at the termination of the sampling process.  Observing that \bee n_s - n_1 \; \Pr \{
\bs{l} = 1 \} & = & n_s \; \Pr \{ \bs{l} \leq s \} - n_1 \; \Pr \{ \bs{l} \leq 1 \} \\
& = & \sum_{\ell = 2}^s \li ( n_\ell \; \Pr \{ \bs{l} \leq \ell \} - n_{\ell - 1} \; \Pr \{ \bs{l} < \ell \} \ri )\\
& = & \sum_{\ell = 2}^s n_\ell \; ( \Pr \{ \bs{l} \leq \ell \} - \Pr \{ \bs{l} < \ell \} ) + \sum_{\ell = 2}^s
(n_\ell - n_{\ell -1} ) \; \Pr \{ \bs{l} < \ell \}\\
& = & \sum_{\ell = 2}^s n_\ell \; \Pr \{ \bs{l} = \ell \}  + \sum_{\ell = 1}^{s-1} (n_{\ell + 1} - n_\ell) \; \Pr \{ \bs{l} \leq \ell \}, \eee
we have $n_s - \sum_{\ell = 1}^s n_\ell \; \Pr \{ \bs{l} = \ell \}  = \sum_{\ell = 1}^{s-1} \; (n_{\ell + 1} - n_\ell) \; \Pr \{ \bs{l} \leq
\ell \}$.  Making use of this result and the fact $n_s = n_1 + \sum_{\ell = 1}^{s-1} \; (n_{\ell + 1} - n_\ell)$, we have \bel \bb{E} [
\mathbf{n} ] & = & \sum_{\ell = 1}^s n_\ell \; \Pr \{ \bs{l} = \ell  \}  =  n_s - \li ( n_s - \sum_{\ell = 1}^s n_\ell \; \Pr \{
\bs{l} = \ell \} \ri ) \nonumber \\
& = & n_1 + \sum_{\ell = 1}^{s-1} \; (n_{\ell + 1} - n_\ell) - \sum_{\ell = 1}^{s-1} \; (n_{\ell + 1} - n_\ell) \; \Pr \{ \bs{l} \leq \ell \} \nonumber\\
& = &  n_1 + \sum_{\ell = 1}^{\tau - 1} \; (n_{\ell + 1} - n_{\ell}) \; \Pr \{ \bs{l} > \ell \} + \sum_{\ell = \tau}^{s - 1} \; (n_{\ell + 1} -
n_{\ell}) \; \Pr \{ \bs{l} > \ell \}. \la{best} \eel  By the definition of the stopping rule, we have \bel \{  \bs{l} > \ell  \} & \subseteq &
\li \{  \li (  \li | \wh{\bs{p}}_\ell - \f{1}{2} \ri | - \ro  \vep \ri
)^2 < \f{1}{4} + \f{ \vep^2 n_\ell}{2 \ln (\ze \de)}  \ri \} \nonumber\\
&  = & \li \{ \ro \vep - \sq{\f{1}{4} + \f{ \vep^2 n_\ell } { 2 \ln (\ze \de) }} < \li | \wh{\bs{p}}_\ell - \f{1}{2} \ri | < \ro \vep +
\sq{\f{1}{4} + \f{ \vep^2 n_\ell } { 2 \ln (\ze \de) }} \ri \} \nonumber\\
& = & \li \{ \ro \vep - \sq{\f{1}{4} + \f{ \vep^2 n_\ell } { 2 \ln (\ze \de) }} < \f{1}{2} - \wh{\bs{p}}_\ell < \ro \vep +
\sq{\f{1}{4} + \f{ \vep^2 n_\ell } { 2 \ln (\ze \de) }}, \; \; \wh{\bs{p}}_\ell \leq \f{1}{2} \ri \} \nonumber\\
&  & \bigcup \li \{ \ro \vep - \sq{\f{1}{4} + \f{ \vep^2 n_\ell } { 2 \ln (\ze \de) }} < \wh{\bs{p}}_\ell - \f{1}{2} < \ro \vep + \sq{\f{1}{4} +
\f{ \vep^2 n_\ell } { 2 \ln (\ze \de) }}, \;\; \wh{\bs{p}}_\ell > \f{1}{2} \ri \} \nonumber\\
& \subseteq & \{ a_\ell < \wh{\bs{p}}_\ell < b_\ell \} \cup \{ 1 - b_\ell < \wh{\bs{p}}_\ell < 1 - a_\ell \} \la{better} \eel for $1 \leq \ell <
s$,  where $b_\ell = \f{1}{2} - \ro \vep + \sq{\f{1}{4} + \f{ \vep^2 n_\ell } { 2 \ln (\ze \de) }}$ for $\ell = 1, \cd, s - 1$.  By the
assumption that $\vep$ and $\ro$ are non-negative,  we have $1 - b_\ell - a_\ell = 2 \ro \vep \geq 0$ for $\ell = 1, \cd, s - 1$.  It follows
from (\ref{better}) that $\{  \bs{l} > \ell  \} \subseteq  \{ \wh{\bs{p}}_\ell >  a_\ell \}$ for $\ell = 1, \cd, s - 1$. By the definition of
$\tau$, we have $p < a_\ell$ for $\tau \leq \ell < s$.  Making use of this fact, the inclusion relationship $\{  \bs{l}
> \ell  \} \subseteq  \{ \wh{\bs{p}}_\ell >  a_\ell \}, \; \ell = 1, \cd, s - 1$,
and Chernoff-Hoeffding bound \cite{Chernoff, Hoeffding}, we have \be \la{bestgood} \Pr \{ \mbf{n} > n_\ell
\mid p \} = \Pr \{ \bs{l} > \ell  \mid p \} \leq \Pr \{ \wh{\bs{p}}_\ell > a_\ell \mid p \} \leq \exp (n_\ell \mscr{M} (a_\ell,p)) \ee for $\tau
\leq \ell < s$. It follows from (\ref{best}) and (\ref{bestgood}) that \bee \bb{E} [ \mbf{n} ] & \leq &  n_1 + \sum_{\ell = 1}^{\tau - 1} \;
(n_{\ell + 1} - n_{\ell}) + \sum_{\ell = \tau}^{s -
1} \; (n_{\ell + 1} - n_{\ell}) \; \Pr \{ \bs{l} > \ell \}\\
& =  & n_\tau + \sum_{\ell = \tau}^{s-1} ( n_{\ell + 1} - n_\ell ) \; \Pr \{ \bs{l} > \ell \} \leq  n_\tau + \sum_{\ell = \tau}^{s-1} ( n_{\ell
+ 1} - n_\ell ) \; \exp (n_\ell \mscr{M} (a_\ell,p)). \eee  This completes the proof of the theorem.

\end{document}